\documentclass[11pt]{amsart}
\usepackage{tocvsec2}
\usepackage{amsmath,amsfonts,latexsym,amssymb,amsthm,mathrsfs,array}

\usepackage[letterpaper]{geometry}
\usepackage[usenames]{color}
\usepackage[colorlinks=true,linkcolor=webgreen,filecolor=webbrown,citecolor=webgreen]{hyperref}

\definecolor{webgreen}{rgb}{0,.5,0}
\definecolor{webbrown}{rgb}{.6,0,0}

\newtheorem{dfn}{Definition}[section]
\newcommand{\bdfn}{\begin{dfn}\rm}
	\newcommand{\edfn}{\end{dfn}}
\newtheorem{thm}[dfn]{Theorem}
\newcommand{\bthm}{\begin{thm}}
	\newcommand{\ethm}{\end{thm}}
\newtheorem{lmma}[dfn]{Lemma}                   
\newcommand{\blmma}{\begin{lmma}}                   
	\newcommand{\elmma}{\end{lmma}}                   
\newtheorem{ppsn}[dfn]{Proposition}
\newcommand{\bppsn}{\begin{ppsn}}
	\newcommand{\eppsn}{\end{ppsn}}
\newtheorem{crlre}[dfn]{Corollary}
\newcommand{\bcrlre}{\begin{crlre}} 
	\newcommand{\ecrlre}{\end{crlre}}
\newtheorem{rmk}[dfn]{Remark}
\newcommand{\brmk}{\begin{rmk}\rm} 
	\newcommand{\ermk}{\end{rmk}}

\numberwithin{equation}{section}

\title[Rank-one geometry and mixed complexes in Cartan type Lie algebras]
{Rank-one geometry and mixed complexes in representations of Cartan type Lie algebras on a torus}
\author{S. Eswara Rao}
\author{Souvik Pal}         
\address{S. Eswara Rao, School of Mathematics, Tata Institute of Fundamental Research,
Homi Bhabha Road, Colaba, Mumbai 400005, India.}                  
\email{sena98672@gmail.com, senapati@math.tifr.res.in}
\address{Souvik Pal, Department of Sciences and Humanities, CHRIST (Deemed to be University), Mysore Road, Bangalore 560 074, India.}     
\email{pal.souvik90@gmail.com, souvik.pal@christuniversity.in}

\date{}

\begin{document}
\let\thefootnote\relax\footnote{Data sharing is not applicable to this article as no datasets were generated and analyzed during the current study.}



\maketitle
\begin{abstract}
In this paper, we develop a unified theory of reducibility and indecomposability for Shen--Larsson modules over the Witt, special and Hamiltonian type Lie algebras on a torus. Our approach is based on a \textit{rank-one mechanism} governing irreducible submodules, Loewy filtrations, rank reduction, uniseriality and mixed complex structures. We first provide a uniform intrinsic characterization of the trivial and fundamental representations of $\mathfrak{gl}_N, \mathfrak{sl}_N,  \mathfrak{sp}_{2n}$ in terms of quadratic relations satisfied by rank-one elements of these matrix Lie algebras and utilize it to determine the irreducibility of Shen--Larsson modules over $\mathcal{W}_N, \mathcal{S}_N, \mathcal{H}_{2n}$. Using the rank-one operators arising from these relations, we then construct rank-reducing operators corresponding to distinguished lattice directions and apply them to show that the submodule structure of the reducible Shen--Larsson modules over $\mathcal{W}_N, \mathcal{S}_N,  \mathcal{H}_{2n}$ attached to the fundamental representations of  $\mathfrak{gl}_N, \mathfrak{sl}_N,  \mathfrak{sp}_{2n}$ respectively are generically uniserial. In the Hamiltonian case, we show that the submodules of these reducible Shen--Larsson modules come from kernels and images of differentials of the de Rham and Koszul-type complexes. These differentials anti-commute and thus endow the tensor field modules with a mixed complex structure, which also admit a natural interpretation formally analogous to the de Rham differential and co-differential type operator appearing in symplectic Hodge theory. In particular, we provide complete answers to the questions recently posed by Pei--Sheng--Tang--Zhao [\textit{J.\ Inst.\ Math.\ Jussieu}\ 2023] concerning the structure of Shen--Larsson modules over $\mathcal{H}_{2n}$. 
\end{abstract} 

\settocdepth{section}
\tableofcontents

\section{Introduction}
\subsection{Background and Motivation} The representation theory of infinite-dimensional Lie algebras exhibits profound connections across algebra, geometry, mathematical physics and homological algebra, with varied applications in symplectic dynamics and  quantum field theory.
In this article, we deal with Cartan type Lie algebras on a torus, which include the Hamiltonian Lie algebra $\mathcal{H}_{2n}$ of symplectic (complex) polynomial vector fields on a $2n$-dimensional torus. This is the space of polynomial vector fields $X = \sum_{i=1}^{2n} f_i(\partial/\partial{t_i})$ satisfying $\mathcal{L}_X\omega =0$, where $\omega = \sum_{i=1}^{n}dt_i \wedge dt_{n+i}$ is the usual symplectic $2$-form, $f_i \in \mathbb{C}[t_1^{\pm1}, \ldots, t_{2n}^{\pm1}]$ and $\mathcal{L}_X$ is the Lie derivative of $\omega$ with respect to $X$.

The Lie algebra $\mathcal{H}_N^{+}$ of Hamiltonian vector fields on $\mathbb{A}^{N}$ is one of the four members of the Cartan type Lie algebras, which also consists of  the Witt series $\mathcal{W}_N^{+}$, the special series $\mathcal{S}_N^{+}$ and the contact series $\mathcal{K}_N^{+}$, with their coefficients belonging to $\mathbb{C}[t_1, \ldots, t_N]$ ($N=2n$ for the Hamiltonian type). In a series of two papers \cite{R, AR}, Rudakov initiated the study of irreducible modules over Lie algebras of type $\widehat{\mathcal{W}}_N^{+}, \widehat{\mathcal{S}}_N^{+}$ and $\widehat{\mathcal{H}}_N^{+}$, where the coefficients lie in $\mathbb{C}[[t_1, \ldots, t_N]]$.  Subsequently, Shen developed his method of mixed products to construct irreducible modules over $\mathcal{W}_N^{+},  \mathcal{S}_N^{+},  \mathcal{H}_N^{+}$ and also for the finite-dimensional counterpart over algebraically closed fields of positive characteristic \cite{GS1}. 

A natural progression is to consider these Cartan type Lie algebras in the context of vector fields on an $N$-dimensional torus $\mathbb{T}^N$, by allowing the  coefficients to lie in $\mathbb{C}[t_1^{\pm1}, \ldots, t_N^{\pm1}]$ and henceforth develop their representation theory. The first step in this direction was taken by Larsson \cite{L}, where he constructed a large class of modules over $\mathcal{W}_N$,  the Lie algebra of polynomial vector fields on $\mathbb{T}^N$. These modules appeared in the work of Larsson in connection to conformal field theory, while Shen's technique of mixed products was motivated by modular representation theory (also see \cite{AF,B}). We shall refer to the modules constructed in this way as Shen--Larsson modules. 

Over the last few decades, many authors have contributed to the study of Shen--Larsson modules. For Lie algebras of type $\mathcal{W}, \ \mathcal{S}$ and $\mathcal{H}$, there exists a natural functor, known as the Shen--Larsson functor, from the category of representations of finite-dimensional simple (or reductive) Lie algebras to the category of representations of the corresponding Cartan type Lie algebras. A fundamental result of Billig--Futorny \cite{BF} shows that all irreducible Harish-Chandra modules over $\mathcal{W}_N$ come from quotients of Shen--Larsson modules, which underlines the central role of these modules in the representation theory of Cartan type Lie algebras.

\subsection{Main Purpose and Challenges} Although the irreducibility criteria and partial submodule descriptions are known for some of these Cartan type Lie algebras (see \cite{R1,GZ, LZ, LGW,T1} and the references therein), a unified structural mechanism explaining irreducibility, indecomposability and the submodule structure is not available.  Our main purpose is to develop such a framework based on a \textit{rank-one mechanism} that uniformly controls the structure of irreducible modules, canonical filtrations and indecomposability properties. In particular, we provide a complete description of the submodules of Shen--Larsson modules over $\mathcal{H}_{N}$, thereby resolving the questions recently recorded in \cite[Remark 4.11]{PSTZ} related to their structure. It is also worth emphasizing here that understanding the structure of Shen--Larsson modules over $\mathcal{H}_{N}$ turns out to be far more delicate and challenging than in the case of either $\mathcal{W}_N$ or $\mathcal{S}_N$, mainly because of the following two reasons, due to which most of the 
existing techniques fail in this setup. 

\begin{itemize}
\item The non-zero $\mathbb{Z}^{N}$-graded components of $\mathcal{H}_{N}$ are extremely thin, i.e. they are $1$-dimensional.
\item There do not exist simple explicit realizations of the fundamental representations of the symplectic Lie algebra $\mathfrak{sp}_{N}$, except for the first fundamental representation.
\end{itemize}


\subsection{Rank-one Characterization of Fundamental Representations}
The starting point of our analysis is a uniform intrinsic characterization detecting the trivial and fundamental modules over the matrix Lie algebras $\mathfrak{gl}_N, \mathfrak{sl}_N$ and $\mathfrak{sp}_{N}$ by means of rank-one operators. For $\mathfrak{g} = \mathfrak{gl}_N, \mathfrak{sl}_N$ or $\mathfrak{sp}_{N}$, let us denote the variety of rank-one elements of  $\mathfrak{g}$ by $\mathcal{R}_1(\mathfrak{g})$, which also spans $\mathfrak{g}$. Then we have the following result from Corollary \ref{CSp}, Corollary \ref{CW} and Corollary \ref{CS} (see Corollary \ref{Nilpotent} and Remark \ref{Orbit} for the connection with \textit{minimal nilpotent orbits}).
\bthm
Let $\mathfrak{g} = \mathfrak{gl}_N, \mathfrak{sl}_N$ or $\mathfrak{sp}_{N}$  and $V$ be any non-zero irreducible module (need not be a weight module) over $\mathfrak{g}$. Then the following are equivalent.
\begin{enumerate}
\item $x^2 - \mathrm{tr}(x)x \in \mathrm{Ann}_{\mathrm{U}(\mathfrak{g})}V \ \forall \ x \in \mathcal{R}_1(\mathfrak{g})$.
\item $V$ is the trivial representation or a fundamental representation of $\mathfrak{g}$.
\end{enumerate}
\ethm
\noindent This rank-one criterion forms the conceptual foundation of our paper, which is initially used to show that the Shen--Larsson module $V \otimes A_{N}$  over $\mathcal{W}_N, \mathcal{S}_N, \mathcal{H}_{N}$ associated to an \textit{irreducible} $\mathfrak{g}$-module $V$, is again irreducible, unless $V$ is the \textit{trivial} or a \textit{fundamental} module over $\mathfrak{gl}_N, \mathfrak{sl}_N,\mathfrak{sp}_{N}$ respectively. This mechanism thus gives a uniform conceptual proof of the irreducibility of Shen--Larsson modules over these three Cartan type Lie algebras (see Theorem \ref{Criterion}, Theorem \ref{TW} and Theorem \ref{TS}).

\subsection{Rank Reduction and Submodule Structure}
Motivated by the above perspective, we construct families of \textit{rank-reducing operators} associated with distinguished lattice directions $k+\beta$ for any $k \in \mathbb{Z}^N$ and $\beta \in \mathbb{C}^N$ satisfying $k+\beta \neq 0$. In case of $\mathcal{W}_N$ and $\mathcal{S}_N$, these operators are induced by the rank-one elements of $\mathfrak{gl}_N$ and $\mathfrak{sl}_N$ respectively (see (\ref{W-invariant}) and Remark \ref{Rank-reducing}),
while for $\mathcal{H}_N$, they are governed by the symplectic rank-one elements (see (\ref{H-invariant})).
A salient feature of these operators is that they eliminate the distinguished direction determined by the vector $k+\beta$ and reduce the action to its transverse complement, so that each $\mathbb{Z}^N$-graded component of the associated Shen--Larsson module naturally acquires the structure of a lower rank module over the matrix algebras: a $\mathfrak{gl}_{N-1}$-module for $\mathcal{W}_N$, an $\mathfrak{sl}_{N-1}$-module for $\mathcal{S}_N$ and an $\mathfrak{sp}_{N-2}$-module for $\mathcal{H}_N$. These rank-reducing operators provide the mechanism for separating vectors inside submodules of reducible Shen--Larsson modules.

Using the above rank-reduction principle, we determine the complete submodule structure of Shen--Larsson modules corresponding to the fundamental $\mathfrak{g}$-modules, which we call \textit{exceptional modules}. 
These modules exhibit a sharp dichotomy depending on the parameter $\beta \in \mathbb{C}^N$. 

\begin{thm}
If $\beta \notin \mathbb{Z}^N$, then the Shen--Larsson modules associated to the fundamental modules are \textit{uniserial} across $\mathcal{W}, \mathcal{S}, \mathcal{H}$ types and their Loewy filtrations \textit{exhaust} all possible submodules. 
\end{thm}

\noindent On the other hand, if $\beta \in \mathbb{Z}^N$, a distinguished homogeneous component becomes trivial and so every subspace of this component forms a submodule. Apart from these trivial summands, no additional submodules occur (see Theorem \ref{Main}, Corollary \ref{CW} and Remark \ref{RS}).

\subsection{Loewy Filtrations in the Hamiltonian Case} In case of $\mathcal{H}_N$, the exceptional Shen--Larsson modules admit canonical semi-simple filtrations 
\begin{align*}
(0)  \subseteq W_{min}^{\beta}\big(V(\omega_p)\big)& \subseteq W_{int}^{\beta}\big(V(\omega_p)\big) \subseteq W_{max}^{\beta}\big(V(\omega_p)\big) \subseteq L_{\mathcal{H}}^{\beta}\big(V(\omega_p)\big); \\
\text{where} \ W_{min}^{\beta}\big(V(\omega_p)\big) & = \mathrm{Im}\big[(T^{\beta} \circ \pi^\beta)|_{L_{\mathcal{H}}^{\beta}\big(V(\omega_p)\big)}\big], \ W_{int}^{\beta}\big(V(\omega_p)\big) = \mathrm{Ker}\big[T^{\beta}|_{L_{\mathcal{H}}^{\beta}\big(V(\omega_p)\big)}\big], \\ W_{max}^{\beta}\big(V(\omega_p)\big) & = \mathrm{Ker}\big[(\pi^{\beta} \circ T^{\beta})|_{L_{\mathcal{H}}^{\beta}\big(V(\omega_p)\big)}\big], \ \text{for} \ \beta \notin \mathbb{Z}^N.
\end{align*}
If $\beta \in \mathbb{Z}^N$, then we obtain the same filtration by removing $V(\omega_p) \otimes \mathbb{C}t^{-\beta}$ from the above two kernels. We prove that the above filtration coincides with a Loewy filtration of $L_{\mathcal{H}}^{\beta}\big(V(\omega_p)\big)$ for $1 \leqslant p \leqslant n$. The boundary cases exhibit additional degeneracies: for $p=1$, the intermediate module equals the maximal module, and for $p=n$, the intermediate and minimal modules coincide (see Lemma \ref{Equal}). Consequently, $L_{\mathcal{H}}^{\beta}\big(V(\omega_p)\big)$ has \textit{Loewy length} $4$ for $1<p<n$ and \textit{Loewy length} $3$ for $p \in \{1,n \}$. Here $T^{\beta}$ and $\pi^{\beta}$ denote the \textit{Koszul-type} and \textit{de Rham-type} differential operators arising from the natural exterior product constructions. More precisely, on each homogeneous component indexed by $k \in \mathbb{Z}^{N}, \ T^{\beta}$ matches with \textit{contraction} by the symplectic dual vector $\overline{k+\beta}$ and  $\pi^{\beta}$ is the \textit{exterior multiplication} by $k+\beta$ (see Lemma \ref{Intermediate} and (\ref{Exterior})). The \textit{intermediate} submodule $W_{int}^{\beta}\big(V(\omega_p)\big)$ thus consists of those $p$-forms transverse to the distinguished direction $k+\beta$ (see Lemma \ref{Lemma} for another realization of $W_{int}^{\beta}\big(V(\omega_p)\big)$), while the \textit{minimal} submodule $W_{min}^{\beta}\big(V(\omega_p)\big)$ is generated by the $p$-forms containing this direction. Moreover, $W_{min}^{\beta}\big(V(\omega_p)\big)$ is also realized using the action of the \textit{rank-one} operators $(k+\beta)(\overline{k+\beta})^T$ on $V(\omega_p)$, where $k$ varies over $\mathbb{Z}^N$ (see Lemma \ref{Similar}).

 Another structural feature of these differential maps is that their composition $f_p^{\beta} := \pi_{p-1}^{\beta} \circ T_p^{\beta}$ gives rise to a canonical \textit{square-zero} endomorphism of $F^{\beta}(\Lambda^{p}\mathbb{C}^N)$, which also remains invariant under $L_{\mathcal{H}}^{\beta}\big(V(\omega_p)\big)$. This operator is induced directly by the symplectic rank-one operators governing the reducing property and is compatible with the mixed complex structure. Explicitly,
\begin{align*}
 f_p^{\beta}(v \otimes t^k) = (k+\beta)(\overline{k+\beta})^Tv \otimes t^k, \ k \in \mathbb{Z}^N,
\end{align*}
which shows that the \textit{maximal submodule} $W_{max}^{\beta}\big(V(\omega_p)\big)$ is intrinsically determined by the associated rank-one operator. Geometrically, $f_p^{\beta}$ measures the component of a form lying in the distinguished rank-one direction and reinserts this direction through exterior multiplication.

\subsection{Special Features of the Hamiltonian Case} A striking difference in the Hamiltonian case is the appearance of the additional Koszul-like differential $T^{\beta}$, which has no analogue in the Witt and special settings. For $\mathcal{W}_N$ and $\mathcal{S}_N$, the submodules of exceptional Shen--Larsson modules arise solely from the de Rham complex associated to $\pi^{\beta}$ (see Corollary \ref{CW} and Remark \ref{RS}). On the other hand, the differential operators $\pi^{\beta}$ and $T^{\beta}$ satisfy the \textit{anti-commutation} relation $T^{\beta}\pi^{\beta} + \pi^{\beta}T^{\beta} =0$, which endows the modules of \textit{tensor fields} on $\mathbb{T}^N$ with the structure of a \textit{mixed complex} in the sense of Kassel \cite{Ka} (see \S \ref{SS7.1}). This mixed complex structure also admits a natural interpretation analogous to the interplay between the de Rham and symplectic co-differentials that appear in symplectic Hodge theory (see \S \ref{Hodge}). 
  
A further distinctive characteristic of the Hamiltonian case is that 
\begin{align*}
\mathrm{Soc}L_{\mathcal{H}}^{\beta}\big(V(\omega_p)\big) \cong \mathrm{Head}L_{\mathcal{H}}^{\beta}\big(V(\omega_p)\big) \cong W_{min}^{\beta}\big(V(\omega_p)\big) \  \text{for} \ \beta \notin \mathbb{Z}^N, 
\end{align*}
while for $\beta \in \mathbb{Z}^{N}$, they differ by finitely many \textit{trivial summands} (see Lemma \ref{JH}, Theorem \ref{Main} and Corollary \ref{C6.8}). This phenomenon is absent in the case of $\mathcal{W}_N$ and $\mathcal{S}_N$, where the \textit{unique} irreducible submodule of the exceptional modules \textit{does not} simultaneously appear as both head and socle.


\subsection{Properties of the Shen--Larsson Functor} We establish uniform functorial properties of the Shen--Larsson construction for the Witt, special and Hamiltonian families. More precisely:
\bthm The Shen--Larsson functor sends every \textit{non-trivial irreducible} $\mathfrak{g}$-module $V$ to the \textit{indecomposable} module $V \otimes A_N$ over the associated Cartan type Lie algebra. All proper submodules of $V \otimes A_N$ are also likewise indecomposable upto trivial summands, across all three families.
\ethm

In contrast, if $V$ is trivial and $\beta \in \mathbb{Z}^N$, then $F_{\mathcal{W}}^{\beta}(V)$ is either \textit{irreducible} or \textit{uniserial} of Loewy length $2$ with exactly one proper submodule; $F_{\mathcal{S}}^{\beta}(V)$ and $L_{\mathcal{H}}^{\beta}(V)$ are \textit{completely reducible}, admitting exactly two Loewy series of length $2$, which exhaust all the proper submodules. If $\beta \notin \mathbb{Z}^N$, then $F_{\mathcal{W}}^{\beta}(V), F_{\mathcal{S}}^{\beta}(V)$ and $L_{\mathcal{H}}^{\beta}(V)$ are all \textit{irreducible} (see Corollary \ref{C6.8}, Corollary \ref{CW} and Remark \ref{RS}).


\subsection{Homology, Realizations and Connection to EALAs} We further compute the homology and cohomology associated with the mixed complex, which shed light on their exactness properties and structural behavior (see Proposition \ref{de Rham} and (\ref{HC})). We also provide various realizations of the Hamiltonian Lie algebra in terms of the space of differential $1$-forms (see \S \ref{Realization}). Finally, we connect these differential $1$-forms to the Hamiltonian extended affine Lie algebras, thereby situating our representation-theoretic results within a broader framework (see \S \ref{EALAs}).

\smallskip

\noindent \textbf{Organization of the paper.} In Section \ref{Notation}, we introduce preliminary notations and terminologies before proceeding to investigate the irreducibility of Shen--Larsson modules over $\mathcal{H}_N$ in the next section. In Section \ref{Submodules}, we present three classes of $\mathcal{H}_N$-submodules of $F^{\beta}(\Lambda^p \mathbb{C}^N)$ and consequently of $L_{\mathcal{H}}^{\beta}\big(V(\omega_p)\big)$, which we call the \textit{minimal}, \textit{intermediate} and \textit{maximal} submodules, all of whose realizations are novel, to the best of our knowledge. We further prove that, unlike $L_{\mathcal{H}}^{\beta}\big(V(\omega_p)\big)$, the Shen--Larsson module $F^{\beta}(\Lambda^p \mathbb{C}^N)$ is \textit{not} uniserial  over $\mathcal{H}_N$ (see Lemma \ref{Inclusion}). In Section \ref{S6}, we classify the submodules of exceptional modules over $\mathcal{H}_N$, while in the following section, we connect our results to homological algebra and symplectic Hodge theory. In Section \ref{S8}, we recover the irreducibility results for Shen--Larsson modules over $\mathcal{W}_N$ and $\mathcal{S}_N$ by following our unified rank-one approach, which substantially simplifies the existing proofs in the literature. In the final section, we determine the submodule structure of exceptional Shen--Larsson modules over $\mathcal{W}_N$ and $\mathcal{S}_N$.
 
\brmk
An earlier version of this work, titled ``Representations of Hamiltonian vector fields on a torus" was posted on arXiv in April 2025 \cite{RP}. In that preprint, we independently established the irreducibility criterion for Shen--Larsson modules over $\mathcal{H}_N$, concurrently with Futorny--Tantubay \cite{FT1} (posted in March 2025). Our approach is different in nature, and the present version is a substantial development of \cite{RP}, introducing a rank-one/orbit-theoretic mechanism  that treats the Witt, special and Hamiltonian families uniformly within a single conceptual framework. In particular, the rank-reducing operators, Loewy filtrations, mixed complex structures and algebraic version of the classical $d\delta$-lemma - which form the bulk of this paper- are not addressed in \cite{FT}.
\ermk


\section{Notations and Preliminaries}\label{Notation}
Throughout this paper, unless otherwise stated, all the vector spaces, algebras, matrices and tensor products are over the field of complex numbers $\mathbb{C}$. We shall denote the set of integers, natural numbers, non-negative integers and non-zero complex numbers by $\mathbb{Z}$, $\mathbb{N}$, $\mathbb{Z_{+}}$ and $\mathbb{C}^{\times}$ respectively. 

\smallskip

Fix any $N \in \mathbb{N}$ and consider the column vector space $\mathbb{C}^N$ of $N \times 1$ matrices with the standard basis $\{e_1, \ldots, e_N \}$. Denote by $(\cdot | \cdot)$  the usual bilinear form on $\mathbb{C}^N$  given by $(u|v) = u^Tv \in \mathbb{C} \ \forall \ u, v \in \mathbb{C}^N$, where $u^T$ is the matrix transpose. Let $A_N = \mathbb{C}[t_1^{\pm 1}, \ldots, t_N^{\pm 1}]$ be the algebra of Laurent polynomials in $N$ variables. For $k = (k_1, \ldots, k_N) \in \mathbb{Z}^N$, write $t^k = t_1^{k_1} \ldots t_N^{k_N}$, which is a typical element of $A_N$.  
 
\subsection{Lie algebras of type $\mathcal{W}_N$ and $\mathcal{S}_N$} Consider the Lie algebra of derivations of $A_N$, popularly known as the Witt algebra, which we shall denote by $\mathcal{W}_N$. Setting $d_i = t_i(\partial/\partial t_i)$ for each $1 \leqslant i \leqslant N$, it is then easy to check that $\mathcal{W}_N = \text{span} \{t^rd_i \ | \ r \in \mathbb{Z}^N, \ 1 \leqslant  i \leqslant N \}$. For $u \in \mathbb{C}^N$ and $r \in \mathbb{Z}^N$, we  define $D(u,r) := \sum_{i=1}^{N}u_it^rd_i$. The Lie bracket operation on $\mathcal{W}_N$ is given by 
\begin{align*}
[D(u,r), D(v,s)] = D\big((u|s)v - (v|r)u, r+s \big) \ \forall \ u, v \in \mathbb{C}^N, \ r, s \in \mathbb{Z}^N.
\end{align*}
This Lie algebra has an interesting subalgebra, namely the divergence-zero vector fields, which we shall denote by $\mathcal{S}_N$. More precisely, it was shown in \cite{T1} that
\begin{align*}
\mathcal{S}_N = \text{span} \{D(u,r) \ | \ (u|r) = 0, \ u \in \mathbb{C}^N, \ r \in \mathbb{Z}^N \},
\end{align*}
whence it follows that $\mathcal{S}_N = \sum_{i=1}^{N}\mathbb{C}d_i \ltimes \big(\bigoplus_{0 \neq r \in \mathbb{Z}^N} (S_N)_r \big)$.

\subsection{Hamiltonian Lie algebra} The Hamiltonian Lie algebra or the algebra of Hamiltonian vector fields on a torus, is known to exist only when $N$ is even. So take $N = 2n$. We closely follow \cite{R2,T2} to define this Lie algebra and denote it by $\mathcal{H}_N$. Put
$\overline{r} := (r_{n+1}, \ldots, r_{2n}, -r_1, -r_2, \ldots, -r_n)$ for each $r = (r_1, \ldots, r_n, r_{n+1}, \ldots, r_N) \in \mathbb{Z}^N$ and define
$h_r := \sum_{i=1}^{n} (r_{n+i}t^rd_i - r_it^rd_{n+i})$.
One can check that $[h_r,h_s] = (\overline{r}|s)h_{r+s}, \ [h_r, h_{-r}] = 0, \ h_0 = 0$.
As earlier, $[D(u,0), h_r] = (u|r)h_r$. Then the corresponding Hamiltonian Lie algebra is defined by  
\begin{align*}
\mathcal{H}_{N} := \mathrm{span} \{h_r,  d_i \ | \ 0 \neq r \in \mathbb{Z}^{N}, \ 1 \leqslant i \leqslant N \}, 
\end{align*}
where the bracket operations are induced from $\mathcal{S}_{N}$. 

\brmk \label{Simple}
\
\begin{enumerate}
\item We also take $\overline{w} := (w_{n+1}, \ldots, w_{2n}, -w_1,-w_2,\ldots,-w_n)$ for $w = (w_1, \ldots , w_N) \in \mathbb{C}^N$. Note that $\mathcal{H}_N = \text{span}\{D(Jr, r), d_i \ | \ r \in \mathbb{Z}^N, \ 1 \leqslant i \leqslant N \}$, where 
\begin{equation*} J =
\begin{pmatrix}
\begin{matrix} 0_{n \times n} \end{matrix} & \mathrm{Id}_{n \times n} \\ -\mathrm{Id}_{n \times n} & \begin{matrix} 0_{n \times n} \end{matrix}
\end{pmatrix}
\end{equation*}

\item $\mathcal{H}_{N}$ is a subalgebra of $\mathcal{S}_{N}$ and in particular, $\mathcal{H}_2$ coincides with $\mathcal{S}_2$.
\item $\mathcal{W}_N, \ \mathcal{S}_N^{\prime} = \bigoplus_{0 \neq r \in \mathbb{Z}^N} \big(\mathcal{S}_N \big)_r$ and $\mathcal{H}_N^{\prime} = \bigoplus_{0 \neq r \in \mathbb{Z}^N} \big(\mathcal{H}_N \big)_r$ are all known to be simple.
\item $D = \sum_{i=1}^{N}\mathbb{C}d_i$ serves as a Cartan subalgebra of $\mathcal{W}_N, \ \mathcal{S}_N$ and $\mathcal{H}_N$.
\end{enumerate}
\ermk

\subsection{Weight modules} \label{Harish-Chandra}
 $V$ is said to be a weight module over $\mathcal{W}_N$, $\mathcal{S}_N$ or $\mathcal{H}_N$ if the action of $D$ on $V$ is diagonalizable. More precisely, we have
\begin{align*}
V = \bigoplus_{\mu \in D^*} V_{\mu}, \ \text{where} \  V_{\mu} = \{v \in V \ | \ d_iv=\mu(d_i)v \ \forall \ 1 \leqslant i \leqslant N \}.
\end{align*}
In this case, the colllection $P_D(V):= \{\mu \in D^* \ | \ V_{\mu} \neq (0) \}$ is known as the set of all \textit{weights} of $V$ (with respect to $D$) and $V_\mu$ is called the \textit{weight space} relative to the weight $\mu$. Moreover, if $\mathrm{dim}V_{\mu} < \infty \ \forall \ \mu \in P_D(V)$, then $V$ is usually referred to as a \textit{Harish-Chandra module}.

\subsection{Finite-dimensional modules over $\mathfrak{gl}_N$ and $\mathfrak{sl}_N$} Let $\mathfrak{gl}_N$ be the Lie algebra of $N \times N$ matrices 
and $\mathfrak{sl}_N$ denotes the subalgebra of traceless matrices. For $1 \leqslant i,j \leqslant N, \ E_{i,j}$ stands for the elementary matrix with $(i,j)$-th entry $1$ and $0$ elsewhere. Then $\mathfrak{h}_{\mathfrak{gl}_N} = \oplus_{i=1}^{N} \mathbb{C}E_{i,i}$ is a Cartan subalgebra of $\mathfrak{gl}_N$. Put $\mathfrak{h}_{\mathfrak{sl}_N} := \mathfrak{h}_{\mathfrak{gl}_N} \cap \mathfrak{sl}_n$, which is a Cartan subalgebra of $\mathfrak{sl}_N$. For each $1 \leqslant i \leqslant N$, define $\epsilon_i \in \mathfrak{h}_{\mathfrak{gl}_N}^*$ by setting $\epsilon_i(\sum_{j=i}^{N}a_j E_{j,j}) = a_i$, where $a_j \in \mathbb{C}$. Consider  $P_{\mathfrak{gl}_n}^{+} = \{\lambda \in \mathfrak{h}_{\mathfrak{gl}_N}^* | \ \lambda (E_{i,i} - E_{i+1,i+1}) \in \mathbb{Z}_{\geqslant 0} \ \forall \ 1 \leqslant i < N \}$ and $\overline{\delta_p} = \sum_{j=1}^{p} \epsilon_j \in \mathfrak{h}_{\mathfrak{gl}_N}^* \ \forall \ 1 \leqslant p \leqslant N$. Set $\overline{\delta_0} := 0$ and $\delta_p := \overline{\delta_p}|_{\mathfrak{h}_{\mathfrak{sl}_N}^*} \ \forall \ 0 \leqslant p \leqslant N$. Note that $\delta_0 = \delta_N = 0$. It is well-known that any finite-dimensional irreducible $\mathfrak{gl}_N$-module ($\mathfrak{sl}_N$-module) is isomorphic to a highest weight module $V(\overline{\lambda})$ (respectively $V(\lambda)$) for a unique $\overline{\lambda} \in P_{\mathfrak{gl}_N}^{+}$ (respectively $\lambda \in P_{\mathfrak{sl}_N}^{+}$).  

Recall that $\mathfrak{gl}_N = \mathfrak{sl}_N \oplus \mathbb{C}\mathrm{Id}$. For any $\lambda \in P_{\mathfrak{sl}_N}^{+}$, we can realize $V(\lambda)$ as a $\mathfrak{gl}_N$-module by simply acting the identity matrix by some $b \in \mathbb{C}$. We shall denote this irreducible $\mathfrak{gl}_N$-module by $V(\lambda,b)$. Then the natural action of $\mathfrak{gl}_N$ (or $\mathfrak{sl}_N$) on $\mathbb{C}^N$ is isomorphic to $V(\delta_1,1)$ (respectively $V(\delta_1)$) and the exterior product $\Lambda^p\mathbb{C}^N = \mathbb{C}^N \wedge \ldots \wedge \mathbb{C}^N$ ($p$ times), where $1 \leqslant p \leqslant N$, is again a $\mathfrak{gl}_N$-module ($\mathfrak{sl}_N$-module), under the action 
\begin{align*}
X(v_1 \wedge \ldots \wedge v_p) = \sum_{i=1}^{p}v_1 \wedge \ldots \wedge v_{i-1} \wedge Xv_i \wedge v_{i+1} \wedge \ldots \wedge v_p \ \forall \ X \in \mathfrak{gl}_N \ (\text{respectively} \ \mathfrak{sl}_N).    
\end{align*}
For $1 \leqslant p < N$, these are known as the \textit{fundamental representations} of $\mathfrak{sl}_N$, with
\begin{align*}
\Lambda^p\mathbb{C}^N \cong V(\delta_p,p) \ \text{as} \ \mathfrak{gl}_n\text{-modules} \ \text{and} \ \Lambda^p\mathbb{C}^N \cong V(\delta_p) \ \text{as} \ \mathfrak{sl}_n\text{-modules}. \\
\text{Also,} \ \Lambda^0 \mathbb{C}^N \cong V(\delta_0) \ \text{as} \ \mathfrak{sl}_n\text{-modules}, \ \text{where we take} \ \Lambda^0 \mathbb{C}^N= \mathbb{C}. 
\end{align*}

\brmk
We shall also refer to $V(\overline{\delta_p}) = V(\delta_p,p)$, where $1 \leqslant p < N$, as the \textit{fundamental representations} of $\mathfrak{gl}_N$.
\ermk

\subsection{Structure of symplectic Lie algebra} \label{Structure}
Consider the non-degenerate and skew-symmetric bilinear form $B$ on $\mathbb{C}^N$ given by $B(u,v) = u^TJv$, where $J$ is the skew-symmetric matrix recorded in Remark \ref{Simple}. Then the \textit{symplectic Lie algebra} (of type $C_n$) is defined by
\begin{align*}
\mathfrak{sp}_{N} := \{ A \in \mathfrak{gl}_N \ | \ A^TJ + JA = 0 \},
\end{align*}
with $\mathfrak{h}_{\mathfrak{sp}_{N}} = \sum_{i=1}^{n} \mathbb{C}H_i$ being a Cartan subalgebra of $\mathfrak{sp}_N$, where $H_i = E_{i,i} - E_{n+i,n+i} \ \forall \ 1 \leqslant i \leqslant n$. The roots of $\mathfrak{sp}_N$ with respect to $\mathfrak{h}_{\mathfrak{sp}_{N}}$ consist of short roots $\{\epsilon_i - \epsilon_j \ | \ 1 \leqslant i \neq j \leqslant n \} \cup \{\pm(\epsilon_i + \epsilon_j) \ | \ 1 \leqslant i \neq j \leqslant n \}$ as well as long roots $\{\pm 2\epsilon_i \ | \ 1 \leqslant i \leqslant n \}$, where $\epsilon_j \in \mathfrak{h}_{\mathfrak{sp}_{N}}^{*}$ satisfies $\epsilon_j(H_i) = \delta_{i,j}$. Put $X_{i,j} = E_{i,j} - E_{n+j,n+i}, \ Y_{i,j} = E_{i,n+j} + E_{j,n+i}, \ Z_{i,j} =  E_{n+j,i} + E_{n+i,j} \ \forall \ 1 \leqslant i \neq j \leqslant n$ and $U_i = E_{i,n+i}, \ V_i = E_{n+i,i} \ \forall \ 1 \leqslant i \leqslant n$, whence $\{X_{i,j}, Y_{i,j}, Z_{i,j} \ | \ 1 \leqslant i \neq j \leqslant n \} \cup \{U_i, V_i \ | \ 1 \leqslant i \leqslant n \}$ is the set of all root vectors of $\mathfrak{sp}_{N}$.  
In particular, $\text{span} \{U_i, V_i, H_i \} \cong \mathfrak{sl}_2 \ \forall \ 1 \leqslant i \leqslant n$.

\brmk \label{basis}
The symplectic Lie algebra can be defined for any basis $\{v_1, w_1,  \ldots, v_n, w_n \}$ of $\mathbb{C}^N$ and any non-degenerate skew-symmetric bilinear form with $B(v_i,w_i) = c \ \forall \ 1 \leqslant i \leqslant n$ and $B(v_i,w_j) = 0 \ \forall \ 1 \leqslant i \neq j \leqslant n$ for some $c \in \mathbb{C}^{\times}$ \cite{FH}. We shall refer to this basis as a \textit{symplectic basis} of $\mathbb{C}^N$. It is well known that any $v_1, w_1 \in \mathbb{C}^N$ satisfying $B(v_1,w_1) \neq 0$ can be extended to a symplectic basis of $\mathbb{C}^N$ \cite{HK}. In \S \ref{Structure}, the symplectic basis is given by $\{e_1, e_{n+1}, \ldots, e_n, e_{2n}\}$. For the remainder of this paper, the bilinear form on $\mathbb{C}^N$ will be taken to be $B(v,w) = (\overline{v}|w) \ \forall \ v, w \in \mathbb{C}^N$. 
\ermk

\subsection{Finite-dimensional $\mathfrak{sp}_N$-modules} 
Let $\alpha_i^{\vee} = H_i-H_{i+1} \ \forall \ 1 \leqslant i < n$ and $\alpha_n^{\vee} = H_n$, which are the simple co-roots of $\mathfrak{sp}_N$, where $N=2n$. Put $P_{\mathfrak{sp}_N}^{+} = \{\lambda \in \mathfrak{h}_{\mathfrak{sp}_N}^* \ | \ \lambda(\alpha_i^{\vee}) \in \mathbb{Z}_{\geqslant 0} \ \forall \ 1 \leqslant i \leqslant n \}$. It is a standard fact that every finite-dimensional irreducible $\mathfrak{sp}_N$-module is isomorphic to a highest weight module $V(\lambda)$ for a unique $\lambda \in P_{\mathfrak{sp}_N}^{+}$. The \textit{fundamental weights} of $\mathfrak{sp}_{N}$ are given by $\omega_p := \sum_{i=1}^{p} \epsilon_i$, for $1 \leqslant p \leqslant n$, where $\epsilon_i \in \mathfrak{h}_{\mathfrak{sp}_N}^*$ is defined by setting $\epsilon_i(H_j) = \delta_{j,i} \ \forall \ 1 \leqslant i,j \leqslant n$.  

Clearly, $\mathfrak{sp}_N$ is a subalgebra of $\mathfrak{sl}_N$ and the natural action of $\mathfrak{sl}_N$ on $\mathbb{C}^N$, when restricted to $\mathfrak{sp}_N$, again remains irreducible and becomes isomorphic to the fundamental module $V(\omega_1)$. However, the other fundamental representations $\Lambda^p\mathbb{C}^N$ of $\mathfrak{sl}_N$, for $2 \leqslant p \leqslant n$, are \textit{not} irreducible over $\mathfrak{sp}_N$.In order to describe these fundamental representations of $\mathfrak{sp}_N$, let us consider the \textit{contraction map} 
	\begin{eqnarray} \label{contraction}
	\theta_p : \Lambda^p \mathbb{C}^N & \longrightarrow & \Lambda^{p-2} \mathbb{C}^N, \ 2 \leqslant p \leqslant N
	\\
	v_1 \wedge \ldots \wedge v_p &\longmapsto & \sum_{i,j=1, \ i<j}^{p} (-1)^{i+j-1}(\overline{v_i}|v_j)\big(v_1 \wedge \ldots \wedge \widehat{v_i} \wedge \ldots \wedge \widehat{v_j} \wedge \ldots \wedge v_p \big),
\end{eqnarray}
where the notation $\widehat{v_i}$ means that the term $v_i$ is omitted from the expression. Then $\theta_p$ is an \textit{onto} $\mathfrak{sp}_N$-module map if $2 \leqslant p \leqslant n$ and the \textit{fundamental representations} of $\mathfrak{sp}_N$ are given by $V(\omega_1) := \mathbb{C}^N$ and $V(\omega_p):= \mathrm{Ker}\theta_p$, with $\mathrm{dim}V(\omega_p) = {N \choose p} - {N \choose p-2}$. Moreover, it is also well-known that  $\Lambda^p\mathbb{C}^N \cong V(\omega_p) \oplus V(\omega_{p-2}) \oplus \ldots \oplus V(\omega_1)$ (or $V(0)$), if $p$ is odd (respectively even) \cite{C,FH,GW}. 
\brmk \label{isomorphism}
\
\begin{enumerate}
\item $\theta_{n+1}$ is an \textit{isomorphism} of $\mathfrak{sp}_N$-modules \cite{E}. 
\item Every finite-dimensional irreducible $\mathfrak{sp}_N$-module $V$ is isomorphic to its (irreducible) dual module $V^*$, which gives $\Lambda^{N-p}(\mathbb{C}^N) \cong \big(\Lambda^p\mathbb{C}^N)^* \cong \Lambda^p(\mathbb{C}^N)^* \cong \Lambda^p \mathbb{C}^N \ \forall \ 1 \leqslant p \leqslant n$ \cite{H}.  
\end{enumerate}
\ermk    

\subsection{Shen--Larsson modules over $\mathcal{W}_N$ and $\mathcal{S}_N$} \label{SL}
 For $\beta \in \mathbb{C}^N$ and a $\mathfrak
 {gl}_N$-module $V$, consider $F^{\beta}_{\mathcal{W}}(V) = V \otimes A_N$. Then $F^{\beta}_{\mathcal{W}}(V)$ is a $\mathcal{W}_N$-module under the action
\begin{align*}
D(u,r)(v \otimes t^k) = (u|k+\beta)v \otimes t^{k+r} + (ru^T)v \otimes t^{k+r} \ \forall \ k \in \mathbb{Z}^N.
\end{align*}
We shall refer to this module as the Shen--Larsson module over $\mathcal{W}_N$. It was shown in \cite{R1} that for each $1 \leqslant p < N, \ F^{\beta}_{\mathcal{W}}\big(V(\overline{\delta_p})\big) = F^{\beta}_{\mathcal{W}}\big(V(\delta_p,p)\big)$ contains an \textit{irreducible} $\mathcal{W}_N$-module 
\begin{align*}
W^{\beta}(\Lambda^p\mathbb{C}^N) = \bigoplus_{k \in \mathbb{Z}^N}W^{k,\beta}(\Lambda^p\mathbb{C}^N) \otimes \mathbb{C}t^k, \  W^{k,\beta}(\Lambda^p\mathbb{C}^N) = \mathbb{C}(k + \beta) \wedge \mathbb{C}^N \ldots \wedge \mathbb{C}^N \subseteq \Lambda^p\mathbb{C}^N. 
\end{align*}
For any $\mathfrak{sl}_N$-module $V$, take the same action of $\mathcal{S}_N^{\prime}$ (see Remark \ref{Simple}) on $F^{\beta}_{\mathcal{S}}(V) = V\otimes A_N$ and put
\begin{align*}
d_i(v \otimes t^k) = (k_i + \alpha_i)(v \otimes t^k) \ \forall \ k \in \mathbb{Z}^N, \ 1 \leqslant i \leqslant N \ \text{and some} \ \alpha \in \mathbb{C}^N.
\end{align*}
This gives the Shen--Larsson module over $\mathcal{S}_N$. In this case, $F^{\beta}_{\mathcal{S}}\big(V(\delta_p)\big)$ also consists of the \textit{irreducible} $\mathcal{S}_N$-module $W^{\beta}(\Lambda^p\mathbb{C}^N)$ \cite{T1}. For $\beta \in \mathbb{Z}^N$, $W^{\beta}(\Lambda^p\mathbb{C}^N)$ is contained in the slightly bigger submodule $W^{\beta}(\Lambda^p\mathbb{C}^N) \oplus \big(\Lambda^p\mathbb{C}^N \otimes \mathbb{C}t^{-\beta} \big)$ over both $\mathcal{W}_N$ and $\mathcal{S}_N$, which we shall denote by $\widehat{W}^{\beta}(\Lambda^p\mathbb{C}^N)$.

\subsection{Shen--Larsson modules over $\mathcal{H}_N$} \label{type H}
For  $\beta \in \mathbb{C}^N$ and any $\mathfrak{sp}_N$-module $V$, define $L^{\beta}_{\mathcal{H}}(V) = V \otimes A_N = \bigoplus_{k \in \mathbb{Z}^N} L^{k,\beta}_{\mathcal{H}}(V) \otimes \mathbb{C}t^k$. Then $\mathcal{H}_N$ acts on $L^{\beta}_{\mathcal{H}}(V)$ via the action
\begin{align*}
h_r(v \otimes t^k) = (\overline{r}|k+\beta)v \otimes t^{k+r} + (r\overline{r}^T)v \otimes t^{k+r}; \\
d_i(v \otimes t^k) = (k_i + \alpha_i)(v \otimes t^k) \ \forall \ k \in \mathbb{Z}^N, \ 1 \leqslant i \leqslant N \ \text{and some} \ \alpha \in \mathbb{C}^N.
\end{align*}
We shall refer to this module as the Shen--Larsson module over $\mathcal{H}_N$.

\brmk \label{Alpha}
\
\begin{enumerate}
\item We shall denote $\Lambda^p \mathbb{C}^N \otimes A_N$ by simply $F^{\beta}(\Lambda^p\mathbb{C}^N), \ 1 \leqslant p \leqslant N$. 
\item Unless explicitly stated, we shall assume that the (zero) degree derivations always act by a grade-shift of $\alpha \in \mathbb{C}^N$ on every Shen--Larsson module over $\mathcal{H}_N$ and $\mathcal{S}_N$.
\end{enumerate}  
\ermk
\blmma \cite{LZ,MZ} \label{L1}
Let $\mathfrak{g}$ be a finite-dimensional simple Lie algebra and $V$ be an irreducible $\mathfrak{g}$-module (need not be a weight module). Then:
\begin{enumerate}
\item Every root vector of $\mathfrak{g}$ either acts injectively or locally nilpotently on $V$.
\item $V$ is finite-dimensional if and only if all the root vectors of $\mathfrak{g}$ act locally nilpotently on $V$. 
\end{enumerate} 
\elmma

\subsection{Loewy Filtration of a Module} 
Let $M$ be a module of finite length over $L$. A filtration of $L$-modules on $M$ is said to be \textit{semi-simple} if the successive quotients are completely reducible. A semi-simple filtration on $M$  of minimal length is known as a \textit{Loewy filtration} (or a \textit{Loewy series}) of $M$. The completely reducible quotients appearing in this filtration are called the \textit{Loewy layers} and the number of Loewy layers is said to be the \textit{Loewy length} of $M$ (which depends only on $M$). 

The \textit{socle} of $M$, which we shall denote by $\mathrm{Soc}M$, is the \textit{sum} of all \textit{irreducible} submodules of $M$. The \textit{radical} of $M$, denoted by $\mathrm{Rad}M$, is the \textit{intersection} of all \textit{maximal} submodules of $M$ and the \textit{head} of $M$ is given by $\mathrm{Head}M := M/\mathrm{Rad}M$.

\subsection{Uniserial Module} 
A module $M$ over $L$ is said to be \textit{uniserial} if its submodule lattice is totally ordered by inclusion. In this case, $M$ has a \textit{unique} composition series (and Loewy series).

\subsection{Annihilator Ideal of a Module} For an $L$-module $V$, the annihilator ideal of $V$ is defined as $\mathrm{Ann}V := \{u \in \mathrm{U}(L) \ | \ u.v = 0 \ \forall \ v \in V \}$, where $\mathrm{U}(L)$ is the universal enveloping algebra of $L$.

\section{Irreducibility criterion for Shen--Larsson modules over $\mathcal{H}_N$}\label{S3}
In this section, we give a sufficient condition for the Shen--Larsson modules over $\mathcal{H}_N$ to be irreducible. We start this section with a simple lemma, which follows with an application of the Vandermonde argument and will be used repeatedly throughout this paper.
\blmma \label{L2}
Let $V$ be a vector space and $W$ be its subspace. For any $M \in \mathbb{N}$, consider the polynomial $P(X_1,X_2, \ldots, X_M)$ in $M$ variables with coefficients in $V$. If $P(X_1, \ldots, X_M) \in W$ for every $X_1, X_2, \ldots, X_M \in \mathbb{Z}$, then each coefficient of $P(X_1,X_2, \ldots, X_M)$ lies in $W$.
\elmma

We now introduce the key lemma of this section, which also gives us a characterization of the trivial and fundamental representations of $\mathfrak{sp}_N$. The analogous results for $\mathfrak{gl}_N$ and $\mathfrak{sl}_N$ are also provided in Lemma \ref{LW} and Lemma \ref{LS} respectively (see Corollary \ref{Nilpotent} and Remark \ref{Orbit} for connections with minimal nilpotent orbits). 

\blmma \label{L3}
Let $V$ be an irreducible $\mathfrak{sp}_{N}$-module (not necessarily a weight module).  Then:  
\begin{enumerate}
\item $\mathfrak{sp}_N =  \mathrm{span}\{r\overline{r}^T \ | \ r \in \mathbb{Z}^N \}  = \mathrm{span} \{u\overline{u}^T \ | \ u \in \mathbb{C}^N \}$. 
\item $(u\overline{v}^T + v\overline{u}^T) \in \mathfrak{sp}_N \ \forall \ u, v \in \mathbb{C}^N$. 
\item $\mathcal{J}_{\mathcal{H}}(V) := \{v \in V \ | \ (r \overline{r}^{T})^2v = 0 \ \forall \ r \in \mathbb{Z}^{N} \}$ is an $\mathfrak{sp}_{N}$-submodule of $V$.
\item $\mathcal{J}_{\mathcal{H}}(V) = V$ if and only if $V \cong V(\omega_k)$ for some $0 \leqslant k \leqslant n$, where we take $\omega_0 = 0$. 
\end{enumerate}
\elmma

\begin{proof}
(1) It can be easily verified that $u\overline{u}^T \in \mathfrak{sp}_N \ \forall \ u \in \mathbb{C}^N$. On the other hand, we have $X_{i,j} = (e_i + e_{n+j}) \overline{(e_i + e_{n+j})}^T - e_i\overline{e_i}^T - e_{n+j}\overline{e_{n+j}}^T, \ Y_{i,j} = (e_i + e_j) \overline{(e_i + e_j)}^T - e_i\overline{e_i}^T - e_{j}\overline{e_{j}}^T, \ Z_{i,j} = (e_{n+i} + e_{n+j}) \overline{(e_{n+i} + e_{n+j})}^T - e_{n+i}\overline{e_{n+i}}^T - e_{n+j}\overline{e_{n+j}}^T \ \forall \ 1 \leqslant i \neq j \leqslant n$ and $U_i = -e_i\overline{e_i}^T, \ V_i = e_{n+i}\overline{e_{n+i}}^T, \ H_i = (e_i + e_{n+i}) \overline{(e_i + e_{n+i})}^T - e_i\overline{e_i}^T - e_{n+i}\overline{e_{n+i}}^T \ \forall \ 1 \leqslant i \leqslant n$, which proves (1).\\
(2) $u\overline{v}^T + v\overline{u}^T = (u+v)(\overline{u+v})^T - u\overline{u}^T - v\overline{v}^T \ \forall \ u, v \in \mathbb{C}^N$ and so we are done using (1).\\
(3) \textbf{Claim.} $\big(v\overline{v}^T(u\overline{v}^T + v\overline{u}^T) + (u\overline{v}^T + v\overline{u}^T)v\overline{v}^T\big)$ acts trivially on $\mathcal{J}_{\mathcal{H}}(V)$ for each $u, v \in \mathbb{C}^N$. \\
Pick any $w \in \mathcal{J}_{\mathcal{H}}(V)$. By applying (1), we then have
\begin{align*}
0 = \big((u+v)\overline{(u+v)}^T)^2w = (u\overline{u}^T)(u\overline{v}^T + v\overline{u}^T)w + (u\overline{v}^T + v\overline{u}^T)(u\overline{u}^T)w  
+ (v\overline{v}^T)(u\overline{v}^T + v\overline{u}^T)w \\
 + (u\overline{v}^T + v\overline{u}^T)(v\overline{v}^T)w \ \forall \ u, v \in \mathbb{C}^N.   
\end{align*}
Now replace $u$ by $pu$ and $v$ by $qv$ for $p, q \in \mathbb{Z}$ in the above equation. Henceforth, comparing the coefficient of $pq^3$ and using Lemma \ref{L2}, the claim follows. \\
Again, for any $r, s \in \mathbb{Z}^N$, consider
\begin{align*}
(s\overline{s}^T)^2(r\overline{r}^T)w = (s\overline{s}^T)[s\overline{s}^T, r\overline{r}^T] w + [s\overline{s}^T, r\overline{r}^T](s\overline{s}^T)w = (\overline{s}|r) \big(s\overline{s}^T(s\overline{r}^T + r\overline{s}^T) + (s\overline{r}^T + r\overline{s}^T)s\overline{s}^T\big)w 
\end{align*}
The desired result is now an immediate consequence of (1) and the above claim. \\
(4) $(u\overline{u}^T)^2$ acts trivially on ${\Lambda}^k\mathbb{C}^N \ (0 \leqslant k \leqslant N)$ for all $u \in \mathbb{C}^N$ , which proves the ``if part". \\
Conversely, assume that $\mathcal{J}_{\mathcal{H}}(V) = V$. We first show that $V$ must be finite-dimensional. \\
By hypothesis, $U_i^2w = V_i^2w = 0  \ \forall \ w \in V$. Consequently, for each $1 \leqslant i \leqslant n$, there exists $0 \neq w_i, w_i^{\prime} \in V$ such that $(V_i)w_i = (U_i)w_i^{\prime} = 0$. \\
\textbf{Claim 1.} $X_{i,j}$ does not act injectively on $V$ for any $1 \leqslant i \neq j \leqslant n$. \\
Taking $r = e_i + e_{n+j} \in \mathbb{Z}^{N}$, we get $(X_{i,j})^2w_j = (r\overline{r}^T)^2w_j = 0$, as $[U_i, V_j] = 0$. Hence the claim. \\
\textbf{Claim 2.} $Y_{i,j}$ does not act injectively on $V$ for any $1 \leqslant i \neq j \leqslant n$. \\
Taking $r = e_i + e_j \in \mathbb{Z}^{N}$, we obtain $(Y_{i,j})^2w_i^{\prime} = (r\overline{r}^T)^2w_i^{\prime} = 0$, since $[U_i, U_j] = 0$. Hence the claim. \\
\textbf{Claim 3.} $Z_{i,j}$ does not act injectively on $V$ for any $1 \leqslant i \neq j \leqslant n$. \\
Taking $r = e_{n+i} + e_{n+j} \in \mathbb{Z}^{N}$, we get $(Z_{i,j})^2w_i = (r\overline{r}^T)^2w_i = 0$, as $[V_i, V_j] = 0$. Hence the claim.  \\
We can now directly invoke Lemma \ref{L1} to conclude that $V$ is finite-dimensional. \\
Assume that $V \ncong V(\omega_k)$ for any $0 \leqslant k \leqslant n$. Then $V \cong V(\lambda)$ for some $\lambda \in P_{\mathfrak{sp}_{N}}^{+} \setminus \{\omega_k\}_{k=0}^{n}$. So there exist $a_1, \ldots, a_k \in \mathbb{Z}_{\geqslant 0}$ such that $\lambda = \sum_{k=1}^{n} a_k \omega_k$ with $a_i + a_j \geqslant 2$ for some $1 \leqslant i \leqslant j \leqslant n$. Note that $\omega_i(H_i) = 1 = \omega_j(H_i)$, which implies that $\lambda(H_i) \geqslant 2$. Now pick $0 \neq v_{\lambda} \in V(\lambda)_{\lambda}$ and consider the $S_i$-module $\text{span} \{(E_{n+i,i})^rv_{\lambda} \ | \ r \in \mathbb{Z}_{\geqslant 0}  \}$, where $S_i := \text{span} \{U_i, V_i, H_i \} \cong \mathfrak{sl}_2$. Consequently, by $\mathfrak{sl}_2$-theory, it immediately follows that $(E_{n+i,i})^2v_{\lambda} \neq 0$ and thus $v_{\lambda} \notin \mathcal{J}_{\mathcal{H}}(V(\lambda))$. 
\end{proof}

\bcrlre \label{CSp}
Let $V$ be any irreducible $\mathfrak{sp}_{N}$-module. Then:
\begin{enumerate}
\item The rank-one operators of $\mathfrak{sp}_{N}$ are given by $\mathcal{R}_1(\mathfrak{sp}_{N}) := \{0 \neq u\overline{u}^T \ | \ u \in \mathbb{C}^{N} \}$.
\item $\mathfrak{sp}_{N}$ is spanned by $\mathcal{R}_1(\mathfrak{sp}_{N})$.
\item $x^2  \in \mathrm{Ann}_{\mathrm{U}(\mathfrak{sp}_{N})}V$ for all $x \in \mathcal{R}_1(\mathfrak{sp}_{N})$ if and only if $V$ is the trivial or a fundamental representation of $\mathfrak{sp}_{N}$.
\end{enumerate}
\ecrlre

\bthm \label{Criterion}
Let $V$ be any irreducible $\mathfrak{sp}_N$-module (need not be a weight module) and $\beta \in \mathbb{C}^N$. Then $L^{\beta}_{\mathcal{H}}(V)$ is irreducible over $\mathcal{H}_N$ if $V \ncong V(\omega_k)$ for any $0 \leqslant k \leqslant n$.
\ethm

\begin{proof} Assume that $V \ncong V(\omega_k)$ for any $0 \leqslant k \leqslant n$. Let $W$ be a non-zero $\mathcal{H}_N$-submodule of $L_{\mathcal{H}}^{\beta}(V)$. Taking $W_m = \{v \in V \ | \ v \otimes t^m \in W \}$, we obtain $W = \bigoplus_{m \in \mathbb{Z}^N} W_m \otimes \mathbb{C}t^m$, since $W$ is $\mathbb{Z}^N$-graded. Set $\widetilde{W} = \bigcap_{m \in \mathbb{Z}^N} W_m$. We now claim that $\widetilde{W} \neq (0)$. \\
By Lemma \ref{L3}, $\mathcal{J}_{\mathcal{H}}(V) = (0)$ and thus for any $0 \neq v \in V$, there exists $r \in \mathbb{Z}^N$ such that $(r\overline{r}^T)^2v \neq 0$. Fix any non-zero $w \in W_m$ for some $m \in \mathbb{Z}^N$. For any $s \in \mathbb{Z}^N$, consider
\begin{align*}
h_{-r}h_{r+s}(w \otimes t^m) = -(\overline{r}|m+r+s+\beta)(\overline{r+s}|m+\beta)w\otimes t^{m+s} + (\overline{r+s}|m+\beta)(r\overline{r}^T)w \otimes t^{m+s} \\
- (\overline{r}|m+s+\beta)(r+s)(\overline{r+s})^Tw \otimes t^{m+s} + (r\overline{r})^T(r+s)(\overline{r+s})^Tw \otimes t^{m+s} \in W_{m+s} \otimes \mathbb{C}t^{m+s}.
\end{align*} 
Replacing $r$ by $pr$ for $p \in \mathbb{Z}$ in the above equation and comparing the coefficients of $p^4$, we can then apply Lemma \ref{L2} to conclude that $(r\overline{r}^T)^2w \in W_{m+s} \ \forall \ s \in \mathbb{Z}^N$. Hence the claim follows. \\ 
Again, for any $m, s \in \mathbb{Z}^N$ and $v \in \widetilde{W}$, we have
\begin{align*}
h_m(v \otimes t^{s-m}) = (\overline{m}|s-m+\beta)v \otimes t^s + (m\overline{m}^T)v \otimes t^s,
\end{align*}
which gives $(m\overline{m}^T)v \in \widetilde{W}$. This proves $\widetilde{W} = V$ using Lemma \ref{L3}, and thus the result follows. 
\end{proof}

\brmk
The construction of $\widetilde{W}$ in Theorem \ref{Criterion} is inspired by \cite{GZ}.
\ermk

\section{Submodules of exceptional Shen--Larsson modules}\label{Submodules}
In this section, we introduce three different classes of $\mathcal{H}_N$-submodules of $L^{\beta}_{\mathcal{H}}\big(V(\omega_p)\big), \ 0 \leqslant p \leqslant n$ and study the inclusions between them.

\subsection{Minimal submodule.} 
For $\beta \in \mathbb{C}^N$ and $1 \leqslant p < N$, set
\begin{align*}
W_{min}^{\beta}(\Lambda^p\mathbb{C}^N) = \bigoplus_{k \in \mathbb{Z}^N}W_{min}^{k,\beta}(\Lambda^p\mathbb{C}^N) \otimes \mathbb{C}t^k, \  W_{min}^{k,\beta}(\Lambda^p\mathbb{C}^N) = \text{span}\{(k+\beta)(\overline{k+\beta})^Tv \ | \ v \in \Lambda^p\mathbb{C}^N \}.
\end{align*}
\bppsn
$W_{min}^{\beta}(\Lambda^p\mathbb{C}^N)$ is an $\mathcal{H}_N$-module.
\eppsn
\begin{proof}
Let $(k+\beta)(\overline{k+\beta})^Tv \otimes t^k \in W_{min}^{k,\beta}(\Lambda^p\mathbb{C}^N)$. For any $r \in \mathbb{Z}^N$, we have
\begin{align*}
h_r\big((k+\beta)(\overline{k+\beta})^Tv \otimes t^k\big) = (\overline{r}|k+\beta)(k+\beta)(\overline{k+\beta})^Tv \otimes t^{k+r} + r\overline{r}^T(k+\beta)(\overline{k+\beta})^Tv \otimes t^{k+r}.
\end{align*} 
So it suffices to  show that $(\overline{r}|k+\beta)(k+\beta)(\overline{k+\beta})^Tv + r\overline{r}^T(k+\beta)(\overline{k+\beta})^Tv \in W_{min}^{(k+r),\beta}(\Lambda^p\mathbb{C}^N)$. To this end, let us consider 
\begin{align*}
(\overline{r}|k+\beta)((k+\beta+r)(\overline{k+\beta+r})^Tv \big) + (k+\beta+r)(\overline{k+\beta+r})^Tr\overline{r}^Tv \in W_{min}^{(k+r),\beta}(\Lambda^p\mathbb{C}^N) \\
=(\overline{r}|k+\beta)\big((k+\beta)(\overline{k+\beta})^T + r\overline{r}^T + (k+\beta)\overline{r}^T + r(\overline{k+\beta})^T \big)v + (k +\beta)(\overline{k+\beta})^Tr\overline{r}^Tv \\
+ (r\overline{r}^T)^2v + \big((k+\beta)\overline{r}^T +r(\overline{k+\beta})^T\big)r\overline{r}^Tv \ldots\ldots (*) 
\end{align*}
Now by Claim of Lemma \ref{L3}, we obtain 
\begin{align*}
0 = r\overline{r}^T\big((k+\beta)\overline{r}^T + r(\overline{k+\beta})^T\big)v + \big((k+\beta)\overline{r}^T + r\overline({k+\beta})^T\big)r\overline{r}^Tv \\ 
= [r\overline{r}^T,(k+\beta)\overline{r}^T + r(\overline{k+\beta})^T]v  + 2\big((k+\beta)\overline{r}^T + r(\overline{k+\beta})^T\big)r\overline{r}^Tv \\
= 2\big((\overline{r}|k+\beta)r\overline{r}^Tv + \big(r(\overline{k+\beta})^T + (k+\beta)\overline{r}^T \big)r\overline{r}^Tv \big), 
\end{align*}
which thereby reduces $(*)$ to
\begin{align*}
(\overline{r}|k+\beta)(k+\beta)(\overline{k+\beta})^Tv + (\overline{r}|k+\beta)\big((k+\beta)\overline{r}^T + r(\overline{k+\beta})^T \big)v + (k +\beta)(\overline{k+\beta})^Tr\overline{r}^Tv \\ 
= (\overline{r}|k+\beta)(k+\beta)(\overline{k+\beta})^Tv + (\overline{r}|k+\beta)\big((k+\beta)\overline{r}^T + r(\overline{k+\beta})^T\big)r\overline{r}^Tv \\
+ [(k+\beta)(\overline{k+\beta})^T, r\overline{r}^T]v 
+ r\overline{r}^T(k+\beta)(\overline{k+\beta})^Tv \\ 
= (\overline{r}|k+\beta)(k+\beta)(\overline{k+\beta})^Tv  + r\overline{r}^T(k+\beta)(\overline{k+\beta})^Tv \in W_{min}^{(k+r),\beta}(\Lambda^p\mathbb{C}^N)  
\end{align*}
and gives us the desired result.
\end{proof}

\noindent For $\beta \in \mathbb{C}^N$ and $1 \leqslant p \leqslant n$, define 
\begin{align*}
W_{min}^{\beta} \big(V(\omega_p) \big) = \bigoplus_{k \in \mathbb{Z}^N} W_{min}^{k,\beta}(V(\omega_p)) \otimes \mathbb{C}t^k, \ \text{where} \ W_{min}^{k,\beta} \big(V(\omega_p)\big) = W_{min}^{k,\beta}(\Lambda^p\mathbb{C}^N)  \cap V(\omega_p). 
\end{align*} 
We shall refer to $W_{min}^{\beta}\big(V(\omega_p) \big)$ as the \textit{minimal module} for $L_{\mathcal{H}}^{\beta}\big(V(\omega_p)\big)$. For $\beta \in \mathbb{Z}^N$, $W_{min}^{\beta}\big(V(\omega_p)\big)$ is contained in the $\mathcal{H}_N$-module $W_{min}^{\beta}\big(V(\omega_p)\big) \oplus \big(V(\omega_p) \otimes \mathbb{C}t^{-\beta} \big)$, which we shall denote by $\widehat{W}_{min}^{\beta} \big(V(\omega_p) \big)$. \\
From \S \ref{type H}, recall that $W^{\beta}(\Lambda^p\mathbb{C}^N)$ is a $\mathcal{W}_N$-module and hence an $\mathcal{H}_N$-module. Furthermore, set
\begin{align*}
W^{\beta} \big(V(\omega_p) \big) = \bigoplus_{k \in \mathbb{Z}^N} W^{k,\beta}(V(\omega_p)) \otimes \mathbb{C}t^k, \ \text{where} \ W^{k,\beta} \big(V(\omega_p)\big) = W^{k,\beta}(\Lambda^p\mathbb{C}^N)  \cap V(\omega_p). 
\end{align*} 
For $\beta \in \mathbb{Z}^N$, $W^{\beta}\big(V(\omega_p)\big)$ is again contained in the $\mathcal{H}_N$-submodule $W^{\beta}\big(V(\omega_p)\big) \oplus \big(V(\omega_p) \otimes \mathbb{C}t^{-\beta} \big)$. We shall denote this module by $\widehat{W}^{\beta} \big(V(\omega_p) \big)$. 

\subsection{Extension of $\mathfrak{sp}_N$-module maps} For an $\mathfrak{sp}_N$-module map $E : V \longrightarrow W$ and $\beta \in \mathbb{C}^N$, set
\begin{eqnarray*}
  \widetilde{E} : L_{\mathcal{H}}^{\beta}(V) & \longrightarrow & L_{\mathcal{H}}^{\beta}(W) \\
  v \otimes t^k &\longmapsto & E(v) \otimes t^k, \ k \in \mathbb{Z}^N, \ v \in V,
\end{eqnarray*}
where $L_{\mathcal{H}}^{\beta}(V)$ and $L_{\mathcal{H}}^{\beta}(W)$ are constructed following \S 2.7 (also see Remark \ref{Alpha}).  

\bppsn \label{extension}
$\widetilde{E}$ is an $\mathcal{H}_N$-module homomorphism.
\eppsn

\begin{proof}
Using Lemma \ref{L3} and the fact that $E$ is an $\mathfrak{sp}_N$-module homomorphism, we have 
\begin{eqnarray*}
h_r \widetilde{E}(v \otimes t^k ) &=& (\overline{r}|k+\beta)\widetilde{E}(v \otimes t^{k}) + r\overline{r}^T \widetilde{E}(v \otimes t^{k}) \\
&=& (\overline{r}|k+\beta)E(v) \otimes t^{k+r}  + r\overline{r}^T \big(E(v)\big) \otimes t^{k+r} \\
&=& (\overline{r}|k+\beta)E(v) \otimes t^{k+r}  + E\big((r\overline{r}^T)v\big) \otimes t^{k+r} \\
&=&E\big((\overline{r}|k+\beta)v + (r\overline{r}^T)(v) \big) \otimes t^{k+r} \\
 &=& \widetilde{E}\big(h_r(v \otimes t^k) \big) \ \forall \ v \in V \ \text{and} \ r, k \in \mathbb{Z}^N.
\end{eqnarray*}
This shows that $\widetilde{E}$ is an $\mathcal{H}_N$-module map. 
\end{proof}

\brmk
Applying Lemma \ref{LW} (Lemma \ref{LS}) reveals that any $\mathfrak{gl}_N$-module map ($\mathfrak{sl}_N$-module map) can be similarly extended to a $\mathcal{W}_N$-module map (respectively to an $\mathcal{S}_N$-module map) for the corresponding Shen--Larsson module over $\mathcal{W}_N$ (respectively $\mathcal{S}_N$).
\ermk

\subsection{Intermediate submodule.} \label{Intermediate}
For $\beta \in \mathbb{C}^N$ and $1 \leqslant p \leqslant N$, define
\begin{eqnarray*}
  \widetilde{\theta}_p^{\beta} : F^{\beta}(\Lambda^{p}\mathbb{C}^N) & \longrightarrow & F^{\beta}(\Lambda^{p-2}\mathbb{C}^N), \ p \neq 1 \\
  v_1 \wedge \ldots \wedge v_p \otimes t^k &\longmapsto & \theta_p(v_1\wedge \ldots \wedge v_p) \otimes t^k, \ k \in \mathbb{Z}^N \\
  	T_p^{\beta} : F^{\beta}(\Lambda^{p}\mathbb{C}^N) & \longrightarrow & F^{\beta}(\Lambda^{p-1}\mathbb{C}^N) \\
  v_1 \wedge \ldots \wedge v_p \otimes t^k &\longmapsto & \sum_{i=1}^{p} (-1)^{i}(\overline{k+\beta}|v_i) \big(v_1 \wedge \ldots \widehat{v_i} \ldots \wedge v_p \big) \otimes t^k, \ k \in \mathbb{Z}^N,
\end{eqnarray*}
where $\theta_p$ is the contraction map as defined in (\ref{contraction}) 
and the notation $\widehat{v_i}$ means that $v_i$ is missing.

\bppsn \label{homomorphisms}
$\widetilde{\theta}_p^{\beta}$ and $T_p^{\beta}$ are $\mathcal{H}_N$-module homomorphisms.
\eppsn

\begin{proof}
From Proposition \ref{extension}, it follows that $\widetilde{\theta}_p^{\beta}$ is an $\mathcal{H}_N$-module map, as $\theta_p$ is an $\mathfrak{sp}_N$-module homomorphism. We now deal with the map  $T_p^{\beta}$. The well-definedness of $T_p^{\beta}$ is immediate from the universal property of exterior products, as $\sum_{i=1}^{p} (-1)^{i}(\overline{k+\beta}|v_i)v_1 \wedge \ldots \widehat{v_i} \ldots \wedge v_p = 0$, whenever $v_i=v_j$ for some $1 \leqslant i \neq j \leqslant p$. Consider\\
$h_r T_p^{\beta}\big((v_1 \wedge \ldots \wedge v_p) \otimes t^k \big)$
\begin{eqnarray*}
&=& (\overline{r}|k+\beta)T_p^{\beta}\big((v_1 \wedge \ldots \wedge v_p) \otimes t^{k+r} \big) + r\overline{r}^TT_p^{\beta}\big((v_1 \wedge \ldots \wedge v_p) \otimes t^{k+r} \big) \\
&=& (\overline{r}|k+\beta) \big(\sum_{i=1}^{p} (-1)^{i}(\overline{k+\beta}|v_i) (v_1 \wedge \ldots \widehat{v_i} \ldots \wedge v_p) \otimes t^{k+r} \big)\\ 
&&+ r\overline{r}^T \big(\sum_{i=1}^{p} (-1)^{i}(\overline{k+\beta}|v_i) (v_1 \wedge \ldots  \widehat{v_i} \ldots \wedge v_p ) \otimes t^{k+r} \big)\\
&=& (\overline{r}|k+\beta) \big(\sum_{i=1}^{p} (-1)^{i}(\overline{k+\beta+r}|v_i) (v_1 \wedge \ldots \widehat{v_i} \ldots \wedge v_p) \otimes t^{k+r} \big) \\
&&+ \sum_{i=1}^{p}(-1)^{i}(\overline{k+\beta}|v_i) \big(\sum_{j \neq i}(\overline{r}|v_j)(v_1 \wedge \ldots \wedge v_{j-1} \wedge r \wedge v_{j+1} \wedge \ldots \widehat{v_i} \ldots \wedge v_p) \otimes t^{k+r}\big) \\
&& - (\overline{r}|k+\beta)\big(\sum_{i=1}^{p}(-1)^{i}(\overline{r}|v_i)(v_1 \wedge \ldots \widehat{v_i} \ldots \wedge v_p) \otimes t^{k+r} \big) 
\end{eqnarray*}
\begin{eqnarray*}
&=& (\overline{r}|k+\beta) \big(\sum_{i=1}^{p} (-1)^{i}(\overline{k+\beta +r}|v_i) (v_1 \wedge \ldots \widehat{v_i} \ldots \wedge v_p) \otimes t^{k+r} \big) \\ 
&& + \sum_{i=1}^{p}(-1)^{i}\big(\sum_{j \neq i}(\overline{k+\beta +r}|v_i)(\overline{r}|v_j)(v_1 \wedge \ldots \wedge v_{j-1} \wedge r \wedge v_{j+1} \wedge \ldots \widehat{v_i} \ldots \wedge v_p) \otimes t^{k+r}\big) \\
&&- \sum_{i=1}^{p}(-1)^{i}\big(\sum_{j \neq i}(\overline{r}|v_i)(\overline{r}|v_j)(v_1 \wedge \ldots \wedge v_{j-1} \wedge r \wedge v_{j+1} \wedge \ldots \widehat{v_i} \ldots \wedge v_p) \otimes t^{k+r}\big) \\
&&- (\overline{r}|k+\beta) \big(\sum_{i=1}^{p} (-1)^{i}(\overline{r}|v_i) (v_1 \wedge \ldots \widehat{v_i} \ldots \wedge v_p) \otimes t^{k+r} \big) \\
&=& (\overline{r}|k+\beta) \big(\sum_{i=1}^{p} (-1)^{i}(\overline{k+\beta +r}|v_i) (v_1 \wedge \ldots \widehat{v_i} \ldots \wedge v_p) \otimes t^{k+r} \big) \\
&& + \sum_{i=1}^{p}(-1)^{i}\big(\sum_{j \neq i}(\overline{k+\beta +r}|v_i)(\overline{r}|v_j)(v_1 \wedge \ldots \wedge v_{j-1} \wedge r \wedge v_{j+1} \wedge \ldots \widehat{v_i} \ldots \wedge v_p) \otimes t^{k+r}\big) \\
&&- (\overline{r}|k+\beta) \big(\sum_{i=1}^{p} (-1)^{i}(\overline{r}|v_i) (v_1 \wedge \ldots \widehat{v_i} \ldots \wedge v_p) \otimes t^{k+r} \big).
\end{eqnarray*}
Using $(\overline{k+\beta +r}|r) = -(\overline{r}|k+\beta)$ and interchanging $i$ and $j$ in the first and last summands gives 
\begin{eqnarray*}
(\overline{r}|k+\beta) \big(\sum_{j=1}^{p} (-1)^{j}(\overline{k+\beta +r}|v_j) (v_1 \wedge \ldots \widehat{v_j} \ldots \wedge v_p) \otimes t^{k+r} \big) + \\
  \sum_{j=1}^{p}(\overline{r}|v_j)\big(\sum_{i \neq j}(-1)^{i}(\overline{k+\beta +r}|v_i)(v_1 \wedge \ldots \wedge v_{j-1} \wedge r \wedge v_{j+1} \wedge \ldots \widehat{v_i} \ldots \wedge v_p) \otimes t^{k+r}\big) + \\
\big(\sum_{j=1}^{p}(-1)^{j}(\overline{r}|v_j)(\overline{k+\beta +r}|r) (v_1 \wedge \ldots  \widehat{v_j} \ldots \wedge v_p) \otimes t^{k+r} \big) \\
= (\overline{r}|k+\beta) \big(\sum_{j=1}^{p} (-1)^{j}(\overline{k+\beta +r}|v_j) (v_1 \wedge \ldots \widehat{v_j} \ldots \wedge v_p) \otimes t^{k+r} \big) + \\ \sum_{j=1}^{p}(-1)^{j}(\overline{r}|v_j)\big(\sum_{i=1}^{j-1}  (-1)^{i}(\overline{k+\beta+r}|v_i)(v_1 \wedge \ldots \widehat{v_i}  \ldots \wedge v_{j-1} \wedge r \wedge v_{j+1} \wedge \ldots \wedge v_p) \otimes t^{k+r}\big) \\
+ \sum_{j=1}^{p}(\overline{r}|v_j)(\overline{k+\beta+r}|r)(v_1 \wedge \ldots \wedge v_{j-1} \wedge v_{j+1} \wedge \ldots \wedge v_p) \otimes t^{k+r} + \\ \sum_{j=1}^{p}(-1)^{j}(\overline{r}|v_j)\big(\sum_{i=j+1}^{p}  (-1)^{i}(\overline{k+\beta+r}|v_i)(v_1 \wedge \ldots \wedge v_{j-1} \wedge r \wedge v_{j+1} \wedge \ldots \widehat{v_i} \ldots \wedge v_p) \otimes t^{k+r}\big) \\
= T_p^{\beta}\big(h_r(v_1 \wedge \ldots \wedge v_p \otimes t^k) \big) \ \forall \ r, k \in \mathbb{Z}^N,
\end{eqnarray*}
which proves that $T_p^{\beta}$ is an $\mathcal{H}_N$-module homomorphism, as required.
\end{proof}

For $\beta \in \mathbb{C}^N$ and $1 \leqslant p \leqslant N$, define
\begin{align*}
W_{int}^{\beta}(\Lambda^{p-1}\mathbb{C}^N) := \bigoplus_{k \in \mathbb{Z}^N}W_{int}^{k,\beta}(\Lambda^{p-1}\mathbb{C}^N) \otimes \mathbb{C}t^k, \ \text{where}\\
 W_{int}^{k,\beta}(\Lambda^{p-1}\mathbb{C}^N) = \text{span} \{\sum_{i=1}^{p} (-1)^{i}(\overline{k+\beta}|v_i)v_1 \wedge \ldots \widehat{v_i} \ldots \wedge v_p \ | \ v_1 \wedge \ldots \wedge v_p \in \Lambda^p\mathbb{C}^N \} \subseteq \Lambda^{p-1}\mathbb{C}^N. 
\end{align*}  
By Proposition \ref{homomorphisms}, $W_{int}^{\beta}(\Lambda^{p-1}\mathbb{C}^N)$ is an $\mathcal{H}_N$-module. Furthermore, for $1 \leqslant p \leqslant n+1$, set   
\begin{align*}
W_{int}^{\beta} \big(V(\omega_{p-1}) \big) = \bigoplus_{k \in \mathbb{Z}^N} W_{int}^{k,\beta}(V(\omega_{p-1})) \otimes \mathbb{C}t^k, \ W_{int}^{k,\beta} \big(V(\omega_{p-1})\big) = W_{int}^{k,\beta}(\Lambda^{p-1}\mathbb{C}^N)  \cap V(\omega_{p-1}). 
\end{align*} 
We shall refer to $W_{int}^{\beta} \big(V(\omega_{p-1})\big)$ as the \textit{intermediate module} for $L_{\mathcal{H}}^{\beta}\big(V(\omega_{p-1})\big)$. For $\beta \in \mathbb{Z}^N$, $W^{\beta}_{int}\big(V(\omega_{p-1})\big)$ is again contained in the $\mathcal{H}_N$-submodule $W^{\beta}_{int}\big(V(\omega_{p-1})\big) \oplus \big(V(\omega_{p-1}) \otimes \mathbb{C}t^{-\beta} \big)$. We shall denote this module by $\widehat{W}_{int}^{\beta} \big(V(\omega_{p-1})\big)$.  

\subsection{Maximal submodule.} For $\beta \in \mathbb{C}^N$ and $1 \leqslant p < N$, define
\begin{align*}
\widehat{W}_{max}^{\beta}(\Lambda^p\mathbb{C}^N) = \bigoplus_{k \in \mathbb{Z}^N}W_{max}^{k,\beta}(\Lambda^p\mathbb{C}^N) \otimes \mathbb{C}t^k, \  \widehat{W}_{max}^{k,\beta}(\Lambda^p\mathbb{C}^N) = \{v \in \Lambda^p\mathbb{C}^N \ | \ (k+\beta)(\overline{k+\beta})^Tv = 0 \}.
\end{align*}  

\bppsn
$\widehat{W}_{max}^{\beta}(\Lambda^p\mathbb{C}^N)$ is an $\mathcal{H}_N$-module.
\eppsn
\begin{proof}
Let $v \in \widehat{W}_{max}^{k,\beta}(\Lambda^p\mathbb{C}^N)$, which implies that $(k+\beta)(\overline{k+\beta})^Tv = 0$. We need to show that $(k + \beta + r)(\overline{k + \beta + r})^T\big(h_r(v \otimes t^k)\big) \in \widehat{W}_{max}^{(k+r),\beta}(\Lambda^p\mathbb{C}^N)$ for any $r \in \mathbb{Z}^N$. To prove this, we consider 
\begin{align*}
(k + \beta + r)(\overline{k + \beta + r})^T\big((\overline{r}|k+\beta)v + r\overline{r}^Tv\big) \\
= (\overline{r}|k+\beta)\big(r\overline{r}^Tv + \{(k+\beta)\overline{r}^T + r(\overline{k+\beta})^T \}v \big) + \big((k+\beta)(\overline{k+\beta})^T+ (k+\beta)\overline{r}^T + r(\overline{k+\beta})^T\big)r\overline{r}^Tv \\ 
= (\overline{r}|k+\beta)\big(r\overline{r}^Tv + \{(k+\beta)\overline{r}^T + r(\overline{k+\beta})^T\}v \big) + r\overline{r}^T(k+\beta)(\overline{k+\beta})^Tv \\
 +  (\overline{k+\beta}|r)\big((k+\beta)\overline{r}^T 
 + r(\overline{k+\beta})^T\big)v + \big((k+\beta)\overline{r}^T 
 + r(\overline{k+\beta})^T\big)r\overline{r}^Tv \\
= (\overline{r}|k+\beta)r\overline{r}^Tv + \big((k+\beta)\overline{r}^T + r(\overline{k+\beta})^T\big)r\overline{r}^Tv.    
\end{align*}
Now by Claim of Lemma \ref{L3}, we have 
\begin{align*}
0 = r\overline{r}^T\big((k+\beta)\overline{r}^T + r(\overline{k+\beta})^T\big)v + \big((k+\beta)\overline{r}^T + r(\overline{k
+\beta})^T\big)r\overline{r}^Tv  \\
= [r\overline{r}^T,(k+\beta)\overline{r}^T + r(\overline{k+\beta})^T]v  + 2\big((k+\beta)\overline{r}^T + r(\overline{k+\beta})^T\big)r\overline{r}^Tv \\ 
= 2\big((\overline{r}|k+\beta)r\overline{r}^Tv + \{(k+\beta)\overline{r}^T + r(\overline{k+\beta})^T \}r\overline{r}^Tv\big), 
\end{align*}
whence we finally get $(k + \beta + r)(\overline{k + \beta + r})^T\big((\overline{r}|k+\beta)w + r\overline{r}^Tv\big) = 0$ and the result follows.
\end{proof}

For $\beta \in \mathbb{C}^N$ and $1 \leqslant p \leqslant n$, define 
\begin{align*}
\widehat{W}_{max}^{\beta} \big(V(\omega_p) \big) = \bigoplus_{k \in \mathbb{Z}^N} \widehat{W}_{max}^{k,\beta}(V(\omega_p)) \otimes \mathbb{C}t^k, \ \text{where} \ \widehat{W}_{max}^{k,\beta} \big(V(\omega_p)\big) = \widehat{W}_{max}^{k,\beta}(\Lambda^p\mathbb{C}^N)  \cap V(\omega_p). 
\end{align*} 
We shall refer to $\widehat{W}_{max}^{\beta}\big(V(\omega_p) \big)$ as the \textit{maximal module} associated with $L_{\mathcal{H}}^{\beta}\big(V(\omega_p)\big)$. For $\beta \in \mathbb{Z}^N$, $\widehat{W}_{max}^{\beta}\big(V(\omega_p)\big)$ contains the $\mathcal{H}_N$-submodule $\widehat{W}_{max}^{\beta}\big(V(\omega_p)\big) \setminus \big(V(\omega_p) \otimes \mathbb{C}t^{-\beta} \big)$. We shall denote this smaller submodule by $W_{max}^{\beta}\big(V(\omega_p)\big)$. Also, we shall take $\widehat{W}_{max}^{\beta}\big(V(\omega_p) \big) = W_{max}^{\beta}\big(V(\omega_p)\big)$ if $\beta \notin \mathbb{Z}^N$.

\blmma \label{Inclusion}
Let $\beta \in \mathbb{C}^N$ and $1 \leqslant p < N$.
\
\begin{enumerate}
\item $W_{min}^{\beta}(\Lambda^p\mathbb{C}^N) \subseteq W^{\beta}(\Lambda^p\mathbb{C}^N) \subseteq W_{max}^{\beta}(\Lambda^p\mathbb{C}^N) \subseteq F^{\beta}(\Lambda^{p}\mathbb{C}^N)$. 
\item $W_{min}^{\beta}(\Lambda^p\mathbb{C}^N) \subseteq W_{int}^{\beta}(\Lambda^p\mathbb{C}^N) \subseteq W_{max}^{\beta}(\Lambda^p\mathbb{C}^N) \subseteq F^{\beta}(\Lambda^{p}\mathbb{C}^N)$.      
\item $W_{int}^{\beta}(\Lambda^0\mathbb{C}^N) \subseteq F^{\beta}(\Lambda^{0}\mathbb{C}^N)$, where $\Lambda^0 \mathbb{C}^N = \mathbb{C}$. 
\end{enumerate}
\elmma
\begin{proof}
(1) and (3) directly follow from the definitions.\\
(2) Fix $k \in \mathbb{Z}^N$ such that $k+\beta \neq 0$. Now, for any $v_1 \wedge \ldots \wedge v_p \in \Lambda^p\mathbb{C}^N$, we have
\begin{align*}
(k+\beta)(\overline{k +\beta})^T\big(v_1 \wedge \ldots \wedge v_p \otimes t^k \big) \\
= \sum_{i=1}^{p}(-1)^{i-1}(\overline{k+\beta}|v_i)\big((k + \beta) \wedge \ldots \wedge v_{i-1} \wedge \widehat{v_i} \wedge v_{i+1} \wedge \ldots \wedge v_p \otimes t^k \big) \\
= T_{p+1}^{\beta}\big(k+ \beta \wedge v_1 \wedge v_2 \ldots \wedge v_p \otimes t^k \big),
\end{align*}
which implies that $W_{min}^{\beta}(\Lambda^p\mathbb{C}^N) \subseteq W_{int}^{\beta}(\Lambda^p\mathbb{C}^N)$.\\
\textbf{Claim.} $(k+\beta)(\overline{k +\beta})^T(\sum_{i=1}^{p} (-1)^{i}(\overline{k+\beta}|v_i) \big(v_1 \wedge \ldots \widehat{v_i} \ldots \wedge v_p) \big) = 0$. \\
We shall prove this claim by induction. For $p= 1$ and $p=2$, the claim follows trivially. Let us now assume that the claim holds good for $p-1$ and consider
\begin{align*}
(k+\beta)(\overline{k +\beta})^T(\sum_{i=1}^{p} (-1)^{i}(\overline{k+\beta}|v_i) \big(v_1 \wedge \ldots \widehat{v_i} \ldots \wedge v_p) \big) \\
= (k+\beta)(\overline{k +\beta})^T(\sum_{i=1}^{p-1} (-1)^{i}(\overline{k+\beta}|v_i) \big(v_1 \wedge \ldots \widehat{v_i} \ldots \wedge v_{p-1})\wedge v_p \big) \\
+ \sum_{i=1}^{p-1} (-1)^{i}(\overline{k+\beta}|v_i) \big(v_1 \wedge \ldots \widehat{v_i} \ldots \wedge v_{p-1}\wedge (k+\beta)(\overline{k +\beta})^Tv_p \big)\\
+ (k+\beta)(\overline{k +\beta})^T(-1)^p(\overline{k+\beta}|v_p)(v_1 \wedge \ldots \wedge v_{p-1}) 
\end{align*}
\begin{align*}
= \sum_{i=1}^{p-1} (-1)^{i}(\overline{k+\beta}|v_i)(\overline{k+\beta}|v_p)\big(v_1 \wedge \ldots \widehat{v_i} \ldots \wedge v_{p-1}\wedge k+\beta \big) \\  
+ (-1)^p (\overline{k+\beta}|v_p) \sum_{i=1}^{p-1}(\overline{k+\beta}|v_i)\big(v_1 \wedge \ldots \wedge v_{i-1} \wedge k+\beta \wedge v_{i+1} \wedge \ldots v_{p-1}\big) = 0,
\end{align*}
by applying our induction hypothesis. Hence the claim follows, which consequently shows that $W_{int}^{\beta}(\Lambda^p\mathbb{C}^N) \subseteq W_{max}^{\beta}(\Lambda^p\mathbb{C}^N)$ and so we get the desired inclusion.
\end{proof}

\brmk
Similar to the constructions of $\widehat{W}_{min}^{\beta}\big(V(\omega_p)\big), \widehat{W}_{int}^{\beta}\big(V(\omega_p)\big)$ and $\widehat{W}_{max}^{\beta}\big(V(\omega_p)\big)$ of $L_{\mathcal{H}}^{\beta}\big(V(\omega_p)\big)$ for $1 \leqslant p \leqslant n$, we also have the $\mathcal{H}_N$-submodules $\widehat{W}_{min}^{\beta}(\Lambda^p\mathbb{C}^N), \widehat{W}_{int}^{\beta}(\Lambda^p\mathbb{C}^N)$ and $\widehat{W}_{max}^{\beta}(\Lambda^p\mathbb{C}^N)$ of $F^{\beta}(\Lambda^{p}\mathbb{C}^N)$ for $1 \leqslant p \leqslant N$ and $\beta \in \mathbb{Z}^N$.    
\ermk

\section{Loewy filtration of exceptional Shen--Larsson modules}
In this section, we explicitly construct Loewy filtrations of $L_{\mathcal{H}}^{\beta}\big(V(\omega_p)\big), \ 0 \leqslant p \leqslant n$. To achieve this, we shall consider the following symplectic basis (see Remark \ref{basis}) of $L_{\mathcal{H}}^{k,\beta}\big(V(\omega_1)\big)$ for any $k \in \mathbb{Z}^N$ with $k+\beta \neq 0$, which will be used repeatedly for the remaining portion of this paper.

Pick $w_1^{(k)} \in \mathbb{Z}^N$ such that $(\overline{k+\beta}|w_1^{(k)}) \neq 0$. Then we can extend $\{k+\beta, w_1^{(k)} \}$ to a symplectic basis of $\mathbb{C}^N$, say $\{k+\beta, w_1^{(k)}, w_2^{(k)}, \ldots, w_{N-2}^{(k)}, w_{N-1}^{(k)} \}$. More precisely, $(\overline{k+\beta}|w_1^{(k)}) = (\overline{w_{2i}^{(k)}}|w_{2i+1}^{(k)}) = c$ (say) for all $1 \leqslant i \leqslant n-1$ and $0$ otherwise. Without loss of generality, we may take $c = 1$.

\blmma \label{Similar}
 $W_{min}^{\beta}\big(V(\omega_p) \big) = W^{\beta} \big(V(\omega_p) \big) \ \forall \ 1 \leqslant p \leqslant n$ and $\beta \in \mathbb{C}^N$.
\elmma

\begin{proof}
We clearly have  $W_{min}^{\beta}(\Lambda^p\mathbb{C}^N) \subseteq W^{\beta}(\Lambda^p\mathbb{C}^N)$, which gives $W_{min}^{\beta}\big(V(\omega_p) \big) \subseteq W^{\beta} \big(V(\omega_p) \big)$. Let 
\begin{align*}
v = \sum_{i_1 < \ldots < i_{p-1}} a_{i_1, \ldots, i_{p-1}} \big(k + \beta \wedge w_{i_1}^{(k)} \wedge \ldots \wedge w_{i_{p-1}}^{(k)}\big) + a_{i_{1}^{\prime} \ldots , i_{p-1}^{\prime}}\big(k + \beta \wedge w_{i_1^{\prime}}^{(k)} \wedge \ldots \wedge w_{i_{p-1}^{\prime}}^{(k)}\big) 
\end{align*}
be a non-zero element of $W^{k,\beta}\big(V(\omega_p)\big)$, where $k \in \mathbb{Z}^N$ such that $k+\beta \neq 0$ and the set containining all the summands of $v$ are taken to be linearly independent. \\
\textbf{Claim.} None of the above summands of $v$ contains $w_1^{(k)}$. \\
If not, then there exists a summand of $v$ of the form $a_{1, j_1, \ldots, j_{p-2}}\big(k + \beta \wedge w_1^{(k)} \wedge w_{j_1}^{(k)} \wedge \ldots \wedge w_{j_{p-2}}^{(k)} \big)$ for some $a_{1, j_1, \ldots, j_{p-2}} \neq 0$ and $j_m \neq 1 \ \forall \ 1 \leqslant m \leqslant p-2$. Without loss of generality, we can take $a_{i_{1}^{\prime} \ldots , i_{p-1}^{\prime}} \neq 0$ and $k + \beta \wedge w_{i_1^{\prime}}^{(k)} \wedge \ldots \wedge w_{i_{p-1}^{\prime}}^{(k)} =  
k + \beta \wedge w_1^{(k)} \wedge w_{j_1}^{(k)} \wedge \ldots \wedge w_{j_{p-2}}^{(k)}$. By hypothesis, 
\begin{align*}
0 = \theta_p(v) = \theta_p\big(a_{i_{1}^{\prime} \ldots , i_{p-1}^{\prime}}(k + \beta \wedge w_1^{(k)} \wedge w_{j_1}^{(k)} \wedge \ldots \wedge w_{j_{p-2}}^{(k)}) + \ldots \ldots \ldots \big) \\
= a_{i_{1}^{\prime} \ldots , i_{p-1}^{\prime}}(w_{j_1}^{(k)} \wedge \ldots \wedge w_{j_{p-2}}^{(k)}) + \ldots \ldots \ldots,
\end{align*}
which implies that $w_{j_1}^{(k)} \wedge \ldots \wedge w_{j_{p-2}}^{(k)}$ occurs in $\theta_p \big(\sum_{i_1 < \ldots < i_{p-1}} a_{i_1, \ldots, i_{p-1}} \big(k + \beta \wedge w_{i_1}^{(k)} \wedge \ldots \wedge w_{i_{p-1}}^{(k)})\big)$. But this is not possible, as $k+\beta$ is present in every summand of $v$. Hence the claim follows. \\
Finally, consider
\begin{align*}
w = \sum_{i_1 < \ldots < i_{p-1}} a_{i_1, \ldots, i_{p-1}} \big(w_1^{(k)} \wedge w_{i_1}^{(k)} \wedge \ldots \wedge w_{i_{p-1}}^{(k)}\big) + a_{i_{1}^{\prime} \ldots , i_{p-1}^{\prime}}\big(w_1^{(k)} \wedge w_{i_1^{\prime}}^{(k)} \wedge \ldots \wedge w_{i_{p-1}^{\prime}}^{(k)}\big),
\end{align*}
which again belongs to $V(\omega_p)$, since $v \in V(\omega_p)$. But then $(k+\beta)(\overline{k+\beta})^Tw = v$, whence we finally obtain $W^{\beta} \big(V(\omega_p) \big) \subseteq W_{min}^{\beta}\big(V(\omega_p) \big)$ and thus the lemma follows.	
\end{proof}

\noindent For $\beta \in \mathbb{C}^N$ and $1 \leqslant p < N$, define
\begin{eqnarray} \label{Exterior}
	\pi_p^{\beta} : F^{\beta}(\Lambda^{p}\mathbb{C}^N) & \longrightarrow & F^{\beta}(\Lambda^{p+1}\mathbb{C}^N) \\
	v_1 \wedge \ldots \wedge v_p \otimes t^k &\longmapsto & (k+\beta)\wedge v_1 \wedge \ldots \wedge v_p \otimes t^k, \ k \in \mathbb{Z}^N.  
\end{eqnarray}
It is well-known that $\pi_p^{\beta}$ is a $\mathcal{W}_N$-module map \cite{R1} and hence an $\mathcal{H}_N$-module map by restriction. 

\blmma \label{JH}
Let $1 \leqslant p < N$ and $\beta \in \mathbb{C}^N$.
\
\begin{enumerate}
	\item $F^{\beta}(\Lambda^{p}\mathbb{C}^N)/\widehat{W}_{max}^{\beta}(\Lambda^p\mathbb{C}^N) \cong W_{min}^{\beta}(\Lambda^p\mathbb{C}^N)$.
	\item $L^{\beta}_{\mathcal{H}}\big(V(\omega_p)\big)/\widehat{W}_{max}^{\beta}\big(V(\omega_p)\big) \cong W_{min}^{\beta}\big(V(\omega_p)\big) , \ 1 \leqslant p \leqslant n$. 
\end{enumerate}
\elmma

\begin{proof}
	Define 
	\begin{eqnarray*}
		f_p^{\beta} : F^{\beta}(\Lambda^{p}\mathbb{C}^N) & \longrightarrow & W_{min}^{\beta}(\Lambda^{p}\mathbb{C}^N)
		\\
		v_1 \wedge \ldots \wedge v_p \otimes t^k &\longmapsto & (k+\beta)(\overline{k+\beta})^T\big(v_1 \wedge v_2 \wedge \ldots \wedge v_{p}) \otimes t^k, \ k \in \mathbb{Z}^N. 
	\end{eqnarray*}
	It is easy to see that $f_p^{\beta} = T_{p+1}^{\beta} \circ \pi_p^{\beta}$, whence by Proposition \ref{homomorphisms}, it follows that $f_p^{\beta}$ is an $\mathcal{H}_N$-module homomorphism. Also, $\text{Ker}f_p^{\beta} = \widehat{W}_{max}^{\beta}(\Lambda^p\mathbb{C}^N)$ and $\text{Im}f_p^{\beta} = W_{min}^{\beta}(\Lambda^p\mathbb{C}^N)$, which proves (1). \\
	Let us now denote the restriction of $f_p^{\beta}$ to $L_{\mathcal{H}}^{\beta}\big(V(\omega_p)\big)$ by $\widetilde{f}_p^{\beta}$. Then $\text{Ker}\widetilde{f}_p^{\beta} = \widehat{W}_{max}^{\beta}\big(V(\omega_p)\big)$ and $\text{Im}\widetilde{f}_p^{\beta} = W_{min}^{\beta}\big(V(\omega_p)\big)$ (since $(k+\beta)(\overline{k+\beta})^T$ commutes with $\theta_p$), from which (2) follows.
\end{proof} 

\bthm \label{Irreducible}
$W_{min}^{\beta}\big(V(\omega_p)\big)$ is irreducible over $\mathcal{H}_N \ \forall \ 1 \leqslant p \leqslant n$ and $\beta \in \mathbb{C}^N$.
\ethm

\noindent We need some preparation to prove this theorem. The next lemma is easy to see.
\blmma \label{Basis}
Let $k \in \mathbb{Z}^N$ and $\beta \in \mathbb{C}^N$ with $k+\beta \neq 0$. Then for $1\leqslant p \leqslant n$, there
exists a basis $\{u_1^{(k)},\ldots, u_a^{(k)}, x_1^{(k)},\ldots, x_b^{(k)}, y_1^{(k)},\ldots, y_c^{(k)}, z_1^{(k)},\ldots, z_d^{(k)}\}$ of $L_{\mathcal{H}}^{k,\beta}\big(V(\omega_p)\big)$ such that $\{u_1^{(k)},\ldots, u_a^{(k)}\}$ is a basis of 
$W_{min}^{k,\beta}\big(V(\omega_p) \big), \{u_1^{(k)},\ldots, u_a^{(k)}, x_1^{(k)},\ldots, x_b^{(k)}\}$ is a basis of $W_{int}^{k,\beta}\big(V(\omega_p)\big)$ and \\$\{u_1^{(k)},\ldots, u_a^{(k)}, x_1^{(k)},\ldots, x_b^{(k)}, y_1^{(k)},\ldots, y_c^{(k)} \}$ is a basis of $\widehat{W}_{max}^{k,\beta}\big(V(\omega_p) \big)$ satisfying the following. 
\begin{enumerate}
\item $u_i^{(k)} \in \mathrm{span}\{(k+\beta)\wedge w_{j_1}^{(k)} \wedge \ldots \wedge w_{j_{p-1}}^{(k)} \ | \ j_r \neq 1, \ j_1 < \ldots < j_{p-1} \} \cap L_{\mathcal{H}}^{k,\beta}\big(V(\omega_p)\big) \ \forall \ 1 \leqslant i \leqslant a$.
\item $x_i^{(k)} \in \mathrm{span}\{w_{i_1}^{(k)} \wedge \ldots \wedge w_{i_{p}}^{(k)} \ | \ i_r \neq 1, \ i_1 < \ldots < i_{p} \} \cap L_{\mathcal{H}}^{k,\beta}\big(V(\omega_p)\big) \ \forall \ 1 \leqslant i \leqslant b$. 
\item For each $1 \leqslant i \leqslant c, \ y_i^{(k)}$ is an element of $L_{\mathcal{H}}^{k,\beta}\big(V(\omega_p)\big)$, having at least one of the summands of the form $(k+\beta) \wedge w_{1}^{(k)}\wedge w_{i_1}^{(k)} \wedge \ldots w_{i_{p-2}}^{(k)}$, where $i_r \neq 1$ and $i_1 < \ldots < i_{p-2}$.
\item $z_i^{(k)} \in \mathrm{span}\{w_1^{(k)} \wedge w_{i_1}^{(k)} \wedge \ldots \wedge w_{i_{p-1}}^{(k)} \ | \ i_r \neq 1, \ i_1 < \ldots <i_{p-1} \} \cap L_{\mathcal{H}}^{k,\beta}\big(V(\omega_p)\big) \ \forall \ 1 \leqslant i \leqslant d$.    
\end{enumerate}
\elmma 

\brmk
\
\begin{enumerate}
\item The basis elements in Lemma \ref{Basis} can be a linear combination of monomials. In particular, $y_i^{(k)}$ cannot be a single monomial, as in that case, $y_i^{(k)}$ will not lie in $L_{\mathcal{H}}^{k,\beta}\big(V(\omega_p)\big)$. 
\item In Lemma \ref{Basis}, it is possible that $\{x_1^{(k)}, \ldots, x_b^{(k)} \} = \emptyset$ or $\{y_1^{(k)}, \ldots, y_c^{(k)} \} = \emptyset$ in some specific cases, as we shall see later in Lemma \ref{Equal}. However, we always have $a, d \geqslant 1$.
\end{enumerate}
\ermk

\bppsn \label{Proposition}
Let $M = \oplus_{k \in \mathbb{Z}^N}M_k \otimes \mathbb{C}t^k$ be an $\mathcal{H}_N$-submodule of $L_{\mathcal{H}}^{\beta}\big(V(\omega_p)\big)$, where $1 \leqslant p \leqslant n$. If there exists $k \in \mathbb{Z}^N$ with $k+\beta \neq 0$ such that $M_k \cap \mathrm{span} \{w_1^{(k)} \wedge w_{i_1}^{(k)} \wedge \ldots \wedge w_{i_{p-1}}^{(k)} \ | \ i_j \neq 1 \} \neq (0)$, then $M = L_{\mathcal{H}}^{\beta}\big(V(\omega_p)\big)$.  
\eppsn

\begin{proof}
Pick $0 \neq v \in M_k \cap \mathrm{span} \{w_1^{(k)} \wedge w_{i_1}^{(k)} \wedge \ldots \wedge w_{i_{p-1}}^{(k)} \ | \ i_j \neq 1 \}$ and any $r \in \mathbb{Z}^N \setminus \{0 \}$. Consider 
\begin{align*}
h_{-r}\big(h_{r+s}(v \otimes t^k)\big) = -(\overline{r}|k+r+s+\beta)(\overline{r} + \overline{s}|k + \beta)v \otimes t^{k+s} + (\overline{r} + \overline{s}|k + \beta)(r\overline{r}^T)v \otimes t^{k+s} \\
- (\overline{r}|k+r+s+\beta)(r+s)(\overline{r+s})^Tv \otimes t^{k+s} + (r\overline{r}^T) (r+s)(\overline{r+s})^Tv \otimes t^{k+s} \\
= -(\overline{r}|k+s+\beta)(\overline{r}|k + \beta)v \otimes t^{k+s} -
(\overline{r}|k+s+\beta)(\overline{s}|k + \beta)v \otimes t^{k+s} + (\overline{r}|k + \beta)(r\overline{r}^T)v \otimes t^{k+s} \\
+ (\overline{s}|k + \beta)(r\overline{r}^T)v \otimes t^{k+s} -  (\overline{r}|k+s+\beta)(r\overline{r}^T)v \otimes t^{k+s} - (\overline{r}|k+s+\beta)(s\overline{s}^T)v \otimes t^{k+s} \\
 - (\overline{r}|k+s+\beta)(r\overline{s}^T + s\overline{r}^T)v \otimes t^{k+s}  
 + (r\overline{r}^T)^2v \otimes t^{k+s} - 
 (r\overline{r}^T)(s\overline{s}^T)v \otimes t^{k+s} \\
 + (r\overline{r}^T)(r\overline{s}^T + s\overline{r}^T)v \otimes t^{k+s} \in M_{k+s} \otimes \mathbb{C}t^{k+s}  \ \forall \ s \in \mathbb{Z}^N \setminus \{0 \}. 
\end{align*}
Now replace $r$ by $qr$ for $q \in \mathbb{Z}$ in the above equation. Henceforth, collecting  the coefficients of $q^2$ and using Lemma \ref{L2}, we obtain
\begin{align*}
(\overline{r}|k+s+\beta)(\overline{r}|k + \beta)v - (\overline{s}|k + \beta)(r\overline{r}^T)v - (r\overline{r}^T)(s\overline{s}^T)v \\ + (\overline{r}|k+s+\beta)(r\overline{s}^T + s\overline{r}^T)v \in M_{k+s} \ \forall \ s \in \mathbb{Z}^N \setminus \{0 \}.  
\end{align*}
Finally, putting $r = w_1^{(k)}$ in the previous equation gives
\begin{align*}
\big(\overline{w_1^{(k)}}|k+s+\beta \big)\big(\overline{w_1^{(k)}}|k + \beta \big)v - \big(\overline{s}|k + \beta \big) \big(w_1^{(k)}\overline{w_1^{(k)}}^T \big)v \\ - \big(w_1^{(k)}\overline{w_1^{(k)}}^T \big)\big(s\overline{s}^T \big)v  + \big(\overline{w_1^{(k)}}|k+s+\beta \big) \big(w_1^{(k)}\overline{s}^T + s\overline{w_1^{(k)}}^T \big)v \\
 = \big(\overline{w_1^{(k)}}|k+\beta \big)^2v - \big(\overline{s}|w_1^{(k)}\big)\big(\overline{w_1^{(k)}}|k+\beta \big)v + 
\big(\overline{s}|w_1^{(k)}\big)^2v 
+ \big(\overline{w_1^{(k)}}|k+\beta \big)\big(\overline{s}|w_1^{(k)}\big)v + \big(\overline{w_1^{(k)}}|s \big)\big(\overline{s}|w_1^{(k)}\big)v \\
 = \big(\overline{w_1^{(k)}}|k+\beta \big)^2v \in M_{k+s} \ \forall \ s \in \mathbb{Z}^N \setminus \{0 \} \ldots \ldots (**), \     
\end{align*}
where we have used the following relations, which are straightforward to deduce.
\begin{align*}
\big(w_1^{(k)}\overline{w_1^{(k)}}^T \big)v = 0, \,\,\ \big(w_1^{(k)}\overline{w_1^{(k)}}^T \big)\big(s\overline{s}^T \big)v = -\big(\overline{s}|w_1^{(k)}\big)^2v, \,\,\ \big(w_1^{(k)}\overline{s}^T + s\overline{w_1^{(k)}}^T \big)v = \big(\overline{s}|w_1^{(k)}\big)v. 
\end{align*}
Then (**) implies that $v \in \cap_{s \in \mathbb{Z}^N}M_{k+s}$ and thus the proposition directly follows by providing a similar argument as in the proof of Theorem \ref{Criterion}. 
\end{proof}

\noindent \textbf{Proof of Theorem \ref{Irreducible}}. By Lemma \ref{JH}, we need to show that there does not exist any proper $\mathcal{H}_N$-submodule of $L^{\beta}\big(V(\omega_p)\big)$ that properly contains $\widehat{W}^{\beta}_{max}\big(V(\omega_p)\big)$. We shall prove the result only for $\beta \notin \mathbb{Z}^N$, as the complementary case can be handled similarly. Let $M = \oplus_{k \in \mathbb{Z}^N}M_k \otimes \mathbb{C}t^k$ be an $\mathcal{H}_N$-module such that $\widehat{W}^{\beta}_{max}\big(V(\omega_p)\big) \subsetneqq M \subseteq L^{\beta}_{\mathcal{H}}\big(V(\omega_p)\big)$. Then by Lemma \ref{Inclusion} and Lemma \ref{Basis}, there exists $0 \neq v \in \mathrm{span} \{w_1^{(k)} \wedge w_{i_1}^{(k)} \wedge \ldots \wedge w_{i_{p-1}}^{(k)} \ | \ i_j \neq 1 \} \cap M_k$ for some $k \in \mathbb{Z}^N$. This implies that $M =L_{\mathcal{H}}^{\beta}\big(V(\omega_p)\big)$ by Proposition \ref{Proposition} and thus we are done.     

\smallskip

We now provide necessary and sufficient conditions for $L_{\mathcal{H}}^{\beta}\big(V(\omega_0)\big)$ to be irreducible as well as its Loewy series when it is not irreducible, in the following corollary.
 
\bcrlre \label{Cor} 
Let $\beta \in \mathbb{C}^N$.
\
\begin{enumerate}
\item $L_{\mathcal{H}}^{\beta}\big(V(\omega_0)\big)$ is irreducible over $\mathcal{H}_N$ if $\beta \notin \mathbb{Z}^N$, where $\omega_0 = 0$.
\item If $\beta \in \mathbb{Z}^N$, then both $V(\omega_0) \otimes \mathbb{C}t^{-\beta}$ and $\bigoplus_{k \in \mathbb{Z}^N, \ k \neq -\beta}\big( V(\omega_0) \otimes \mathbb{C}t^k\big)$ are irreducible $\mathcal{H}_N$-submodules of $L_{\mathcal{H}}^{\beta}\big(V(\omega_0)\big)$.
\item If $\beta \in \mathbb{Z}^N$, then $L_{\mathcal{H}}^{\beta}\big(V(\omega_0)\big)/V(\omega_0) \otimes \mathbb{C}t^{-\beta}$ is irreducible over $\mathcal{H}_N$.
\end{enumerate}
\ecrlre

\begin{proof}
First note that $L_{\mathcal{H}}^{\beta}\big(V(\omega_0)\big) \cong A_N$ for any $\beta \in \mathbb{C}^N$. Now define
	\begin{eqnarray*}
	\pi_0^{\beta} : L_{\mathcal{H}}^{\beta}\big(V(\omega_0)\big) & \longrightarrow & L_{\mathcal{H}}^{\beta}\big(V(\omega_1)\big) 
	\\
	 t^k &\longmapsto & (k+\beta) \otimes t^k, \ k \in \mathbb{Z}^N. 
\end{eqnarray*}
Then $\pi_0^{\beta}$ is an $\mathcal{H}_N$-module homomorphism and consequently (1) follows, thanks to Theorem \ref{Irreducible}. It is clear that $V(\omega_0) \otimes \mathbb{C}t^{-\beta}$ is irreducible over $\mathcal{H}_N$, for $\beta \in \mathbb{Z}^N$. Moreover, by restricting $\pi_0^{\beta}$ to $\bigoplus_{k \in \mathbb{Z}^N, \ k \neq -\beta}\big( V(\omega_0) \otimes \mathbb{C}t^k\big)$, we also obtain (2). Finally, (3) is immediate by just observing that $\mathrm{Ker}\pi_0^{\beta} = V(\omega_0) \otimes \mathbb{C}t^{-\beta}$ and $\mathrm{Im}\pi_0^{\beta} = \bigoplus_{k \in \mathbb{Z}^N, \ k \neq -\beta}\big( V(\omega_0) \otimes \mathbb{C}t^k\big)$, for $\beta \in \mathbb{Z}^N$.
\end{proof}

\brmk \label{Isom}
$\widehat{W}_{min}^{0}\big(V(\omega_1)\big) \cong \mathcal{H}_N$ and $W_{min}^{0}\big(V(\omega_1)\big) \cong A_N/\mathbb{C} \cong \mathcal{H}_N^{\prime}$ (see Remark \ref{Simple}).
\ermk

In the following two propositions, we investigate the successive quotients for the inclusions recorded in Lemma \ref{Inclusion}, which helps us to understand how the nature of these quotients varies when we restrict our attention to $L_{\mathcal{H}}^{\beta}\big(V(\omega_p)\big)$ (also see Lemma \ref{Similar} and Lemma \ref{JH}).
   
\bppsn \label{Equality}
Let $1 \leqslant p < N$ and $\beta \in \mathbb{C}^N$.
\
\begin{enumerate} 
	\item $W^{\beta}(\Lambda^p\mathbb{C}^N)/W_{min}^{\beta}(\Lambda^p\mathbb{C}^N) \cong W_{min}^{\beta}(\Lambda^{p-1}\mathbb{C}^N), \ p \geqslant 2$.
	\item If $\beta \in \mathbb{Z}^N$, then $\widehat{W}^{\beta}(\Lambda^p\mathbb{C}^N)/\widehat{W}_{min}^{\beta}(\Lambda^p\mathbb{C}^N) \cong W_{min}^{\beta}(\Lambda^{p-1}\mathbb{C}^N), \ p \geqslant 2$.
    \item $\widehat{W}^{\beta}_{max}(\Lambda^p\mathbb{C}^N)/\widehat{W}^{\beta}(\Lambda^p\mathbb{C}^N) \cong W_{min}^{\beta}(\Lambda^{p+1}\mathbb{C}^N)$.
	\item If $\beta \in \mathbb{Z}^N$, then $\widehat{W}_{max}^{\beta}(\Lambda^p\mathbb{C}^N)/\widehat{W}^{\beta}(\Lambda^p\mathbb{C}^N) \cong W_{min}^{\beta}(\Lambda^{p+1}\mathbb{C}^N)$. 
\end{enumerate}
\eppsn
\begin{proof} (1) Let us denote the restriction of $T_p^{\beta}$ to $W^{\beta}(\Lambda^p\mathbb{C}^N)$ by $\widetilde{T}_p^{\beta}$. Then 
\begin{align*}
\widetilde{T}_p^{\beta}\big(k+\beta \wedge v_1 \wedge \ldots \wedge v_{p-1} \otimes t^k \big) = (k+\beta)(\overline{k+\beta})^T\big(v_1 \wedge v_2 \wedge \ldots \wedge v_{p-1}\big) \otimes t^k \ \forall \ k \in \mathbb{Z}^N.
\end{align*}	
\textbf{Claim.} $\text{Ker}\widetilde{T}_p^{\beta} = W_{min}^{\beta}(\Lambda^p\mathbb{C}^N)$. \\    
	For any $k \in \mathbb{Z}^N$, it is clear that $\text{span} \{(k+\beta) \wedge w_{i_1}^{(k)} \wedge \ldots \wedge w_{i_{p-1}}^{(k)} \otimes t^k \ | \ i_m \neq 1 \} \subseteq \text{Ker}\widetilde{T}_p^{\beta}|_{W^{k,\beta}(\Lambda^p\mathbb{C}^N) \otimes \mathbb{C}t^k}$. Now for $(k+ \beta) \wedge w_1^{(k)} \wedge w_{j_1}^{(k)} \wedge \ldots \wedge w_{j_{p-1}}^{(k)} \in F^{\beta}(\Lambda^{p}\mathbb{C}^N)$, we have
	\begin{align*}
	(k+\beta)(\overline{k+\beta})^T\big(w_1^{(k)} \wedge w_{j_1}^{(k)} \wedge \ldots \wedge w_{j_{p-1}}^{(k)} \big) \neq 0.
	\end{align*}
	 But again,
	 \begin{align*}
	 W_{min}^{k,\beta}(\Lambda^p\mathbb{C}^N) = \text{span} \{(k+\beta) \wedge w_{i_1}^{(k)} \wedge \ldots \wedge w_{i_{p-1}}^{(k)} \ | \ i_m \neq 1 \},
	\end{align*}
	which proves our claim.      
	It is also evident that $\text{Im}\widetilde{T}_p^{\beta} = W_{min}^{\beta}(\Lambda^{p-1}\mathbb{C}^N)$, from which (1) and subsequently (2) follows. \\
	(3) Let us consider $\widehat{\pi}_p^{\beta} := \pi_p^{\beta}|_{\widehat{W}_{max}^{\beta}(\Lambda^p\mathbb{C}^N)}$. Note that
    \begin{align*}
	\widehat{\pi}_p^{\beta}(w_1^{(k)} \wedge w_{i_1}^{(k)} \wedge \ldots \wedge w_{i_{p-1}}^{(k)}) \neq 0, \ \text{where} \ i_r \neq 1 \ \text{and} \\
	\widehat{W}_{max}^{k,\beta}(\Lambda^{p}\mathbb{C}^N) = \text{span}\{(k+\beta) \wedge w_{j_1}^{(k)}\wedge \ldots \wedge w_{j_{p-1}}^{(k)}, w_{i_1}^{(k)} \wedge \ldots \wedge w_{i_{p}}^{(k)} \ | \ i_r \neq 1 \}.
	\end{align*} 
	for any $k \in \mathbb{Z}^N$ satisfying $k+\beta \neq 0$. Then it is easy to see that $\mathrm{Im}\widehat{\pi}_p^{\beta} = W_{min}^{\beta}(\Lambda^{p+1}\mathbb{C}^N)$ and $\mathrm{Ker}\widehat{\pi}_p^{\beta} = W^{\beta}(\Lambda^{p}\mathbb{C}^N)$, which proves (3) as well as (4).    
\end{proof}

\bppsn \label{P5.10}
Let $\beta \in \mathbb{C}^N$.
\
\begin{enumerate}
\item $W_{int}^{\beta}(\Lambda^{p-1}\mathbb{C}^N)/W_{min}^{\beta}(\Lambda^{p-1}\mathbb{C}^N) \cong W_{min}^{\beta}(\Lambda^{p}\mathbb{C}^N), \ 1<p < N$. 
\item If $\beta \in \mathbb{Z}^N$, then $\widehat{W}_{int}^{\beta}(\Lambda^{p-1}\mathbb{C}^N)/\widehat{W}_{min}^{\beta}(\Lambda^{p-1}\mathbb{C}^N) \cong W_{min}^{\beta}(\Lambda^{p}\mathbb{C}^N), \ 1<p < N$. 
\item $W_{max}^{\beta}(\Lambda^{p-1}\mathbb{C}^N)/W_{int}^{\beta}(\Lambda^{p-1}\mathbb{C}^N) \cong W_{min}^{\beta}(\Lambda^{p-2}\mathbb{C}^N), \ 2<p \leqslant N$. 
\item If $\beta \in \mathbb{Z}^N$, then $\widehat{W}_{max}^{\beta}(\Lambda^{p-1}\mathbb{C}^N)/\widehat{W}_{int}^{\beta}(\Lambda^{p-1}\mathbb{C}^N) \cong W_{min}^{\beta}(\Lambda^{p-2}\mathbb{C}^N), \ 2<p \leqslant N$.  
\end{enumerate}
\eppsn
\begin{proof}
(1) We can check that
\begin{align*}
\text{span}\{(k+\beta) \wedge w_{i_1}^{(k)} \wedge \ldots \wedge w_{i_{p-2}}^{(k)}, \ w_{j_1}^{(k)} \wedge \ldots \wedge w_{j_{p-1}}^{(k)} | \ i_r, j_s  \neq 1 \} \subseteq W_{int}^{k,\beta}(\Lambda^{p-1}\mathbb{C}^N).
\end{align*}
 Also, note that neither $(k+\beta) \wedge w_1^{(k)} \wedge w_{i_1}^{(k)} \wedge \ldots \wedge w_{i_{p-3}}^{(k)} \ (i_r \neq 1)$ nor $w_1^{(k)} \wedge w_{j_1}^{(k)} \wedge \ldots \wedge w_{j_{p-2}}^{(k)} \ (j_s \neq 1)$ belongs to $W_{int}^{k,\beta}(\Lambda^{p-1}\mathbb{C}^N)$ and so we get
\begin{align*} 
 W_{int}^{k,\beta}(\Lambda^{p-1}\mathbb{C}^N) = \text{span}\{(k+\beta) \wedge w_{i_1}^{(k)} \wedge \ldots \wedge w_{i_{p-2}}^{(k)}, \ w_{j_1}^{(k)} \wedge \ldots \wedge w_{j_{p-1}}^{(k)} | \ i_r, j_s  \neq 1 \}. 
\end{align*}
Then it can be easily deduced that $\text{Im}\big(\pi_p^{\beta}|_{W_{int}^{\beta}(\Lambda^{p-1}\mathbb{C}^N)}\big) = W_{min}^{\beta}(\Lambda^{p}\mathbb{C}^N)$ and $\text{Ker}\big(\pi_p^{\beta}|_{W_{int}^{\beta}(\Lambda^{p-1}\mathbb{C}^N)}\big) = W_{min}^{\beta}(\Lambda^{p-1}\mathbb{C}^N)$, which establishes (1) and consequently (2).  \\
To conclude the proof of this lemma, it suffices to prove (3). It is clear that
\begin{align*}
W^{\beta}(\Lambda^{p-1}\mathbb{C}^N) + W^{\beta}_{int} (\Lambda^{p-1}\mathbb{C}^N) = \widehat{W}^{\beta}_{max}(\Lambda^{p-1}\mathbb{C}^N), \\ 
W^{\beta}(\Lambda^{p-1}\mathbb{C}^N) \cap W^{\beta}_{int} (\Lambda^{p-1}\mathbb{C}^N) = W^{\beta}_{min}(\Lambda^{p-1}\mathbb{C}^N).
\end{align*}
Then by Proposition \ref{Equality}, 
\begin{align*}
\widehat{W}^{\beta}_{max}(\Lambda^{p-1}\mathbb{C}^N)/\widehat{W}^{\beta}_{int}(\Lambda^{p-1}\mathbb{C}^N) \cong W^{\beta}(\Lambda^{p-1}\mathbb{C}^N)/W^{\beta}_{min}(\Lambda^{p-1}\mathbb{C}^N) \cong W_{min}^{\beta}(\Lambda^{p-2}\mathbb{C}^N)
\end{align*}
and therefore we are done.
\end{proof}

The next result gives us an alternate realization of the intermediate submodule of $L_{\mathcal{H}}^{\beta}\big(V(\omega_p)\big)$, which also plays a vital role in determining the composition series of $L_{\mathcal{H}}^{\beta}\big(V(\omega_n)\big)$. 
 
\blmma \label{Lemma}
Let $1 < p \leqslant n$ and $\beta \in \mathbb{C}^N$.
\begin{enumerate}
\item $W_{int}^{\beta}\big(V(\omega_{p})\big) = \mathrm{Ker}\big(\widetilde{\theta}_{p+1}^{\beta}\circ \pi_{p}^{\beta}|_{W_{max}^{\beta}\big(V(\omega_{p})\big)}\big) = \mathrm{Ker}\big(\widetilde{\theta}_{p+1}^{\beta} \circ \pi_{p}^{\beta}|_{L^{\beta}\big(V(\omega_{p})\big)}\big)$, where $\widetilde{\theta}_{p}^{\beta}$ is the $\mathcal{H}_N$-module map defined in \S \ref{Intermediate}.  
\item $\widehat{W}_{max}^{\beta}\big(V(\omega_{p})\big)/\widehat{W}_{int}^{\beta}\big(V(\omega_{p})\big) \cong W_{min}^{\beta}\big(V(\omega_{p-1})\big)$.
\item If $\beta \in \mathbb{Z}^N$, then $\widehat{W}_{max}^{\beta}\big(V(\omega_{p})\big)/\widehat{W}_{int}^{\beta}\big(V(\omega_{p})\big) \cong W_{min}^{\beta}\big(V(\omega_{p-1})\big)$.
\end{enumerate}
\elmma

\begin{proof}
Fix any $k \in \mathbb{Z}^N$ satisfying $k+\beta \neq 0$. Then
\begin{align*}
(1) \ W_{int}^{k,\beta}\big(V(\omega_p)\big) = \mathrm{span}\{(k+\beta) \wedge w_{i_1}^{(k)} \wedge \ldots \wedge w_{i_{p-1}}^{(k)}, \ w_{j_1}^{(k)} \wedge \ldots \wedge w_{j_{p}}^{(k)} | \ i_r, j_s  \neq 1 \} \cap L^{k,\beta}\big(V(\omega_p)\big),
\end{align*}
which implies that 
\begin{align*}
W_{min}^{\beta}\big(V(\omega_{p})\big) \subseteq \mathrm{Ker}\big(\widetilde{\theta}_{p+1}^{\beta}\circ \pi_{p}^{\beta}|_{W_{max}^{\beta}\big(V(\omega_{p})\big)}\big), \\   
\mathrm{Im}\big(\pi_p^{\beta}|_{W_{max}^{k,\beta}\big(V(\omega_{p})\big)}\big) \subseteq \mathrm{span}\{ (k+\beta) \wedge w_{j_1}^{(k)} \wedge \ldots \wedge w_{j_{p}}^{(k)} | \  j_s  \neq 1 \}.
\end{align*}
Now since $\theta_{p+1}^{\beta}\big((k+\beta) \wedge w_{j_1}^{(k)} \wedge \ldots \wedge w_{j_{p}}^{(k)}\big) = (k+\beta) \wedge \big(\theta_{p}^{\beta}(w_{j_1}^{(k)} \wedge \ldots \wedge w_{j_{p}}^{(k)})\big) \ (j_s \neq 1)$, we obtain 
\begin{align*}
W_{int}^{\beta}\big(V(\omega_{p})\big) \subseteq \mathrm{Ker}\big(\widetilde{\theta}_{p+1}^{\beta}\circ \pi_{p}^{\beta}|_{W_{max}^{\beta}\big(V(\omega_{p})\big)}\big). 
\end{align*}
Finally, we can easily check that if $v$ is a non-zero element of    	
$L_{\mathcal{H}}^{k,\beta}\big(V(\omega_p)\big)$ of the form $y_i^{(k)}$ or $z_i^{(k)}$ (see Lemma \ref{Basis}) or their linear combination, then $\big(\widetilde{\theta}_{p+1}^{\beta}\circ \pi_{p}^{\beta}\big)(v \otimes t^k) \neq 0$ and hence we are done. 
(2) Put $\widehat{T}_p^{\beta}:= T_p^{\beta}|_{\widehat{W}_{max}^{\beta}\big(V(\omega_{p})\big)}$.
Now by Lemma \ref{Basis}, if $0 \neq v \in \widehat{W}_{max}^{k,\beta}\big(V(\omega_{p})\big)$ contains $w_1^{(k)}$ in one of its summands, then $k+\beta$ must also occur in that same summand and in that case, $\widehat{T}_p^{\beta}(v \otimes t^k) \neq 0$.  Also, it is evident from our above description of $\widehat{W}_{int}^{k,\beta}\big(V(\omega_p)\big)$ that $\widehat{T}_p^{\beta}|_{\widehat{W}_{int}^{k,\beta}\big(V(\omega_p)\big) \otimes \mathbb{C}t^k} = 0$, which yields that $\mathrm{Ker}\widehat{T}_p^{\beta} = \widehat{W}_{int}^{\beta}\big(V(\omega_p)\big)$. Next, we establish that 
$\mathrm{Im}\widehat{T}_p^{\beta} = W_{min}^{\beta}\big(V(\omega_{p-1})\big)$. Pick
\begin{align*} 
w= (k+\beta \wedge w_1^{(k)}) \wedge (\sum_{i_1 < \ldots < i_{p-2}} w_{i_1}^{(k)} \wedge \ldots \wedge w_{i_{p-2}}^{(k)}) + \sum_{j_1 < \ldots < j_{p}} w_{j_1}^{(k)}\wedge \ldots \wedge w_{j_p}^{(k)} \\
+ \sum_{l_1 < \ldots < l_{p-1}} (k + \beta) \wedge w_{l_1}^{(k)}\wedge \ldots \wedge w_{l_{p-1}}^{(k)}  \in W_{max}^{k,\beta}\big(V(\omega_{p})\big), 
\end{align*}
where $i_q, j_r, l_s \neq 1$. Then $\widehat{T}_p^{\beta}(w \otimes t^k) = -(k + \beta) \wedge (\sum_{i_1 < \ldots < i_{p-2}} w_{i_1}^{(k)} \wedge \ldots \wedge w_{i_{p-2}}^{(k)}) \otimes t^k$, with
\begin{align*}
\theta_{p-1}\big((k + \beta) \wedge (\sum_{i_1 < \ldots < i_{p-2}} w_{i_1}^{(k)} \wedge \ldots \wedge w_{i_{p-2}}^{(k)})\big) = (k+\beta) \wedge \theta_{p-2} \big(\sum_{i_1 < \ldots < i_{p-2}} w_{i_1}^{(k)} \wedge \ldots \wedge w_{i_{p-2}}^{(k)}\big) = 0,
\end{align*} 
since we have $\theta_p(w) = 0$. This implies that $(0) \neq \mathrm{Im}\widehat{T}_p^{\beta} \subseteq L_{\mathcal{H}}^{\beta}\big(V(\omega_{p-1})\big) \cap W_{min}^{\beta}(\Lambda^{p-1}\mathbb{C}^N)$. The desired result now follows from Lemma \ref{Similar} and Theorem \ref{Irreducible}, whereas a similar argument gives (4).
	
\end{proof}

\blmma \label{L5.12}
Let $1\leqslant p<n$ and $\beta \in \mathbb{C}^N$.
\
\begin{enumerate}
\item $W_{int}^{\beta}\big(V(\omega_{p})\big)/W_{min}^{\beta}\big(V(\omega_{p})\big) \cong W_{min}^{\beta}\big(V(\omega_{p+1})\big)$. 
\item If $\beta \in \mathbb{Z}^N$, then $\widehat{W}_{int}^{\beta}\big(V(\omega_{p})\big)/\widehat{W}_{min}^{\beta}\big(V(\omega_{p})\big) \cong W_{min}^{\beta}\big(V(\omega_{p+1})\big)$.  
\end{enumerate}
\elmma
\begin{proof}
Setting $\widetilde{\pi}_p^{\beta} := \pi_p^{\beta}|_{W_{int}^{\beta}\big(V(\omega_{p})\big)}$ and using (1) of Lemma \ref{Lemma}, we obtain 
\begin{align*}
\theta_{p+1}(k+\beta \wedge v) = (k+\beta)\wedge \theta_p(v) = 0 \ \forall \ v \in W_{int}^{\beta}\big(V(\omega_{p})\big) 
\end{align*}
and as a result, we have $(0) \neq \mathrm{Im}\widetilde{\pi}_p^{\beta} \subseteq L_{\mathcal{H}}^{\beta}\big(V(\omega_{p+1})\big)$. But again, $\mathrm{Im}\widetilde{\pi}_p^{\beta} \subseteq W_{min}^{\beta}(\Lambda^{p+1}\mathbb{C}^N)$, which implies that $\mathrm{Im}\widetilde{\pi}_p^{\beta} = W_{min}^{\beta}\big(V(\omega_{p+1})\big)$, by Lemma \ref{Similar} and Theorem \ref{Irreducible}. Finally an analogous argument also gives $\mathrm{Ker}\widetilde{\pi}_p^{\beta} = W_{min}^{\beta}\big(V(\omega_{p})\big)$, which proves (1) and consequently (2).    	
\end{proof}

The following lemma conveys that the minimal and maximal modules  can coincide with the intermediate module in some specific cases.

\blmma \label{Equal}
Let $\beta \in \mathbb{C}^N$.
\
\begin{enumerate}
\item $\widehat{W}_{max}^{\beta}\big(V(\omega_{1})\big) =\widehat{W}_{int}^{\beta}\big(V(\omega_{1})\big)$. 
\item $W_{int}^{\beta}\big(V(\omega_{n})\big) = W_{min}^{\beta}\big(V(\omega_{n})\big)$.
\end{enumerate}
\elmma
\begin{proof}
Fix any $k \in \mathbb{Z}^N$ with $k+\beta \neq 0$. \\
(1) $T_2^{\beta}(k+\beta \wedge w_1^{(k)}) = k + \beta$ and $T_2^{\beta}(w_i^{(k)} \wedge w_1^{(k)}) = w_i^{(k)} \ \forall \ 2 \leqslant i \leqslant N-1$, which implies that 
\begin{align*}
\widehat{W}_{max}^{k,\beta}\big(V(\omega_{1})\big) = \{v \in \mathbb{C}^N \ | \ (\overline{k+ \beta}|v) = 0 \} = \text{span}\{k+\beta, w_i^{(k)} \ | \ 2 \leqslant i < N \} \\
 \subseteq \widehat{W}_{int}^{k,\beta}\big(V(\omega_{1})\big) \subseteq  \widehat{W}_{max}^{k,\beta}\big(V(\omega_{1})\big) \ (\text{by Lemma} \ \ref{Inclusion}). 
\end{align*}
(2) Using Lemma \ref{Lemma} and Remark \ref{isomorphism}, we obtain
\begin{align*}
W_{int}^{k,\beta}\big(V(\omega_{n})\big) \otimes \mathbb{C}t^k = \mathrm{Ker}\big(\widetilde{\theta}_{n+1}^{\beta}|_{F^{\beta}(\Lambda^{n+1}\mathbb{C}^N) \otimes \mathbb{C}t^k} \circ \pi_{n}^{\beta}|_{L^{k,\beta}(V(\omega_{n})) \otimes \mathbb{C}t^k }\big) \\
 = \mathrm{Ker}\big(\pi_{n}^{\beta}|_{L^{k,\beta}(V(\omega_{n})) \otimes \mathbb{C}t^k}\big) 
 =   W_{min}^{k,\beta}\big(V(\omega_{n})\big) \otimes \mathbb{C}t^k.   
\end{align*} 
\end{proof}

\noindent Summarizing some of the results that we have proved in this section together with Lemma \ref{Inclusion}, we readily obtain the Loewy series for $L_{\mathcal{H}}^{\beta}\big(V(\omega_p)\big)$ for $1 \leqslant p \leqslant n$, whereas the Loewy series for $L_{\mathcal{H}}^{\beta}\big(V(\omega_0)\big)$ follows from Corollary \ref{Cor}.  

\bthm \label{Composition}
\
\begin{enumerate}
\item If $1 < p < n$ and $\beta \notin \mathbb{Z}^N$, then a Loewy filtration of $L_{\mathcal{H}}^{\beta}\big(V(\omega_p)\big)$ is given by 
\begin{align*}
(0) \subseteq W_{min}^{\beta}\big(V(\omega_{p})\big) \subseteq W_{int}^{\beta}\big(V(\omega_{p})\big) \subseteq W_{max}^{\beta}\big(V(\omega_{p})\big) \subseteq L_{\mathcal{H}}^{\beta}\big(V(\omega_p)\big),
\end{align*} 
where the Loewy layers of $L_{\mathcal{H}}^{\beta}\big(V(\omega_p)\big)$ are $W_{min}^{\beta}\big(V(\omega_{p-1})\big), \ W_{min}^{\beta}\big(V(\omega_{p})\big), \ W_{min}^{\beta}\big(V(\omega_{p+1})\big)$, with multiplicities
$[L_{\mathcal{H}}^{\beta}\big(V(\omega_p)\big) : W_{min}^{\beta}\big(V(\omega_{p-1})\big)] = [L_{\mathcal{H}}^{\beta}\big(V(\omega_p)\big) : W_{min}^{\beta}\big(V(\omega_{p+1})\big)] = 1$ and $[L_{\mathcal{H}}^{\beta}\big(V(\omega_p)\big) : W_{min}^{\beta}\big(V(\omega_{p})\big)] = 2$.  
\item If $1 < p < n$ and $\beta \in \mathbb{Z}^N$, then a Loewy filtration of $L_{\mathcal{H}}^{\beta}\big(V(\omega_p)\big)$ is given by 
\begin{align*}
(0) \subseteq W_{min}^{\beta}\big(V(\omega_{p})\big) \subseteq W_{int}^{\beta}\big(V(\omega_{p})\big) \subseteq W_{max}^{\beta}\big(V(\omega_{p})\big) \subseteq L_{\mathcal{H}}^{\beta}\big(V(\omega_p)\big), 
\end{align*}
where the Loewy layers of $L_{\mathcal{H}}^{\beta}\big(V(\omega_p)\big)$ are $\widehat{W}_{min}^{\beta}\big(V(\omega_{p-1})\big), \ W_{min}^{\beta}\big(V(\omega_{p})\big), \ W_{min}^{\beta}\big(V(\omega_{p+1})\big)$, with multiplicities 
$[L_{\mathcal{H}}^{\beta}\big(V(\omega_p)\big) : \widehat{W}_{min}^{\beta}\big(V(\omega_{p-1})\big)] = [L_{\mathcal{H}}^{\beta}\big(V(\omega_p)\big) : W_{min}^{\beta}\big(V(\omega_{p+1})\big)] = 1$ and
 $[L_{\mathcal{H}}^{\beta}\big(V(\omega_p)\big) : W_{min}^{\beta}\big(V(\omega_{p})\big)] = 2$.
\item If $\beta \notin \mathbb{Z}^N$, then a Loewy filtration of $L_{\mathcal{H}}^{\beta}\big(V(\omega_1)\big)$ is given by 
\begin{align*}
(0) \subseteq W_{min}^{\beta}\big(V(\omega_{1})\big) \subseteq W_{max}^{\beta}\big(V(\omega_{1})\big) \subseteq L_{\mathcal{H}}^{\beta}\big(V(\omega_1)\big),  
\end{align*} 
where the Loewy layers of $L_{\mathcal{H}}^{\beta}\big(V(\omega_1)\big)$ are given by $W_{min}^{\beta}\big(V(\omega_{1})\big)$ and $W_{min}^{\beta}\big(V(\omega_{2})\big)$, with multiplicities
$[L_{\mathcal{H}}^{\beta}\big(V(\omega_1)\big) : W_{min}^{\beta}\big(V(\omega_{2})\big)] = 1$ and $[L_{\mathcal{H}}^{\beta}\big(V(\omega_1)\big) : W_{min}^{\beta}\big(V(\omega_{1})\big)] = 2$.
\item If $\beta \in \mathbb{Z}^N$, then a Loewy filtration of $L_{\mathcal{H}}^{\beta}\big(V(\omega_1)\big)$ is given by 
\begin{align*}
(0) \subseteq W_{min}^{\beta}\big(V(\omega_{1})\big) \subseteq W_{max}^{\beta}\big(V(\omega_{1})\big) \subseteq L_{\mathcal{H}}^{\beta}\big(V(\omega_1)\big), 
\end{align*} 
where the Loewy layers of $L_{\mathcal{H}}^{\beta}\big(V(\omega_1)\big)$ are given by $W_{min}^{\beta}\big(V(\omega_{1})\big)$ and $\widehat{W}_{min}^{\beta}\big(V(\omega_{2})\big)$, with multiplicities
$[L_{\mathcal{H}}^{\beta}\big(V(\omega_1)\big) : \widehat{W}_{min}^{\beta}\big(V(\omega_{2})\big)] = 1$ and $[L_{\mathcal{H}}^{\beta}\big(V(\omega_1)\big) : W_{min}^{\beta}\big(V(\omega_{1})\big)] = 2$. 
\item If $\beta \notin \mathbb{Z}^N$, then a Loewy filtration of $L_{\mathcal{H}}^{\beta}\big(V(\omega_n)\big)$ is given by
\begin{align*}
(0) \subseteq W_{min}^{\beta}\big(V(\omega_{n})\big) \subseteq W_{max}^{\beta}\big(V(\omega_{n})\big) \subseteq L_{\mathcal{H}}^{\beta}\big(V(\omega_n)\big), \ where
\end{align*} 
the Loewy layers of $L_{\mathcal{H}}^{\beta}\big(V(\omega_n)\big)$ are $W_{min}^{\beta}\big(V(\omega_{n})\big)$ and $W_{min}^{\beta}\big(V(\omega_{n-1})\big)$, with multiplicities
$[L_{\mathcal{H}}^{\beta}\big(V(\omega_n)\big) : W_{min}^{\beta}\big(V(\omega_{n-1})\big)] = 1$ and $[L_{\mathcal{H}}^{\beta}\big(V(\omega_n)\big) : W_{min}^{\beta}\big(V(\omega_{n})\big)] = 2$.
\item If $\beta \in \mathbb{Z}^N$, then a Loewy filtration of $L_{\mathcal{H}}^{\beta}\big(V(\omega_1)\big)$ is given by
\begin{align*}
(0) \subseteq W_{min}^{\beta}\big(V(\omega_{n})\big) \subseteq W_{max}^{\beta}\big(V(\omega_{n})\big) \subseteq L_{\mathcal{H}}^{\beta}\big(V(\omega_n)\big), \ where 
\end{align*} 
the Loewy layers of $L_{\mathcal{H}}^{\beta}\big(V(\omega_n)\big)$ are $W_{min}^{\beta}\big(V(\omega_{n})\big)$ and $\widehat{W}_{min}^{\beta}\big(V(\omega_{n-1})\big)$, with multiplicities
$[L_{\mathcal{H}}^{\beta}\big(V(\omega_n)\big) : \widehat{W}_{min}^{\beta}\big(V(\omega_{n-1})\big)] = 1$ and $[L_{\mathcal{H}}^{\beta}\big(V(\omega_n)\big) : W_{min}^{\beta}\big(V(\omega_{n})\big)] = 2$. 
\end{enumerate}
\ethm

\brmk \label{S_2}
If $n=p=1$, then Lemma \ref{Equal} reveals that the minimal, intermediate and maximal submodules of $L_{\mathcal{H}}^{\beta}\big(V(\omega_1)\big) = F^{\beta}_{\mathcal{S}}\big(V(\delta_1)\big)$ all coincide, which conveys that the \textit{Loewy series} for the exceptional module $L_{\mathcal{H}}^{\beta}\big(V(\omega_1)\big)$ over $\mathcal{H}_2 (\cong \mathcal{S}_2)$ has \textit{Loewy length 2}. 
\ermk

\section{Classification of submodules of exceptional Shen--Larsson modules over $\mathcal{H}_N$} \label{S6}
In this section, we completely classify all possible submodules of $L_{\mathcal{H}}^{\beta}\big(V(\omega_p)\big)$ for $0 \leqslant p \leqslant n$. To this end, we first introduce a class of operators in $\mathfrak{sp}_N$, which is absolutely pivotal for our cause.

\subsection{Rank-reducing operators of $\mathfrak{sp}_{N}$} \label{H-invariant}
For every $k, r,s \in \mathbb{Z}^N$ and $\beta \in \mathbb{C}^N$, put 
\begin{align*}
\mathcal{T}_{k,r,s}^{\mathcal{H},\beta} = (\overline{k+\beta}|r)s - (\overline{k+\beta}|s)r. 
\end{align*}
Let us now define the corresponding \textit{rank-reducing operators} of $\mathfrak{sp}_{N}$ by setting
\begin{align*}
\Omega_{k,r,s}^{\mathcal{H},\beta} = \mathcal{T}_{k,r,s}^{\mathcal{H},\beta}\overline{\mathcal{T}_{k,r,s}^{\mathcal{H},\beta}}^T. 
\end{align*}

\blmma \label{Invariance}
Let $M = \oplus_{k \in \mathbb{Z}^N}(M_k \otimes \mathbb{C}t^k)$ and $M^{\prime} = \oplus_{k \in \mathbb{Z}^N}(M^{\prime}_k \otimes \mathbb{C}t^k)$ be $\mathcal{H}_N$-submodules of $L_{\mathcal{H}}^{\beta}(V)$ and $L_{\mathcal{H}}^{\beta}(V^{\prime})$ respectively, where $V$ and $V^{\prime}$ are non-trivial $\mathfrak{sp}_N$-modules. Then:
\begin{enumerate}
\item $M_k$ is invariant under $\Omega_{k,r,s}^{\mathcal{H},\beta}$ for each $k,r,s \in \mathbb{Z}^N$ and $\beta \in \mathbb{C}^N$. 
\item If $\phi: M \longrightarrow M^{\prime}$ is an $\mathcal{H}_N$-module map, then
$\phi\big(\Omega_{k,r,s}^{\mathcal{H},\beta}.(v \otimes t^m)\big) = \Omega_{k,r,s}^{\mathcal{H},\beta}.\phi(v \otimes t^m) \ \forall \ m \in \mathbb{Z}^N$.
\end{enumerate}
\elmma

\begin{proof}
(1) For any $v \in M_k$, consider 
\begin{align*}
h_{-s}h_{-r}\big(h_{r+s}(v \otimes t^k)\big) \\
= (\overline{r+s}|k+\beta)(\overline{r}|k+s+ \beta)\big[(\overline{s}|k+\beta)(v \otimes t^k) - (s\overline{s}^T)v \otimes t^k\big] - (\overline{r+s}|k+\beta)\big[(\overline{s}|k+ \beta)(r\overline{r}^T)v \otimes t^k \\  
-(s\overline{s}^T)(r\overline{r}^T)v \otimes t^k \big]+ (\overline{r}|k+ s + \beta)\big[(\overline{s}|k + \beta)(r+s) (\overline{r+s})^Tv \otimes t^k - (s\overline{s}^T)(r+s)(\overline{r+s})^Tv \otimes t^k\big] \\
-(\overline{s}|k+\beta)\big[(r\overline{r}^T)(r+s)(\overline{r+s})^Tv \otimes t^k - (s\overline{s}^T)(r\overline{r}^T)(r+s)(\overline{r+s})^Tv \otimes t^k\big] \\
 = - \Omega_{k,r,s}^{\mathcal{H},\beta}.(v \otimes t^k) + \ldots \ldots \ldots \in M_k \otimes \mathbb{C}t^k.
\end{align*}
Replacing $r$ by $pr$ and $s$ by $qs$ for $p, q \in \mathbb{Z}$ in the above equation, the result now follows by simply comparing the coefficients of $p^2q^2$ and then applying Lemma \ref{L2}. \\
(2) For all $r, s \in \mathbb{Z}^N$ and $v \in M_k$, we have
\begin{align*}
\phi\bigg((h_{-s}h_{-r}\big(h_{r+s}(v \otimes t^k)\big)\bigg) = \bigg(h_{-s}h_{-r}\big(h_{r+s}(\phi(v \otimes t^k)\big)\bigg). 
\end{align*}
Now expand both sides using the action of $\mathcal{H}_N$ and then replace $r$ by $pr$, $s$ by $qs$ for $p, q \in \mathbb{Z}$. The result subsequently follows by just comparing the coefficients of $p^2q^2$ and invoking Lemma \ref{L2}. 
\end{proof}

\bppsn \label{Prop}
Let $M = \oplus_{k \in \mathbb{Z}^N}(M_k \otimes \mathbb{C}t^k)$ be an $\mathcal{H}_N$-submodule of $L_{\mathcal{H}}^{\beta}(V)$, where $V$ is a non-trivial $\mathfrak{sp}_N$-module. 
\begin{enumerate}
\item If $(\overline{k+\beta})_i \neq 0$ for some $1 \leqslant i \leqslant N$, then $\{\mathcal{T}_{k,r,s}^{\mathcal{H},\beta} \ | \ r, s \in \mathbb{Z}^N \}$ spans $W_{max}^{k,\beta}\big(V(\omega_1)\big)$. Moreover, $\{\mathcal{T}_{k,e_i,e_j}^{\mathcal{H},\beta} \ | \ 1 \leqslant j \neq i \leqslant N \}$ is a basis of $W_{max}^{k,\beta}\big(V(\omega_1)\big)$. 
\item $(w\overline{w}^T).v \in M_k \ \forall \ v \in M_k$ and $w \in \mathbb{C}^N$ satisfying $(\overline{k+\beta}|w) = 0$.
\item $M_k$ is invariant under $u_1\overline{u_2}^T + u_2\overline{u_1}^T \ \forall \ u_1, u_2  \in \mathbb{C}^N$ satisfying $(\overline{k+\beta}|u_1) = 0 = (\overline{k+\beta}|u_2)$.
\end{enumerate} 
\eppsn

\begin{proof}
(1) If $(\overline{k+\beta})_i \neq 0$ for some $1 \leqslant i \leqslant N$, then it can be shown using elementary arguments that $\{\mathcal{T}_{k,e_i,e_j}^{\mathcal{H},\beta} \ | \ 1 \leqslant j \neq i \leqslant N \}$ is a linearly independent subset of $W_{max}^{k,\beta}\big(V(\omega_1)\big)$, which proves the assertion, as we already know that $W_{max}^{k,\beta}\big(V(\omega_1)\big)$ is $(N-1)$-dimensional. \\
(2) Without loss of generality, let $(\overline{k+\beta})_1 \neq 0$. Clearly $W_{max}^{k,\beta}\big(V(\omega_1)\big) = \{ v \in \mathbb{C}^N \ | \ (\overline{k+\beta}|v) = 0 \}$. By (1), there exist $a_2, \ldots, a_N \in \mathbb{C}$ such that $w = \sum_{i=2}^{N}a_i\mathcal{T}_{k,e_1,e_i}^{\mathcal{H},\beta}$. Then for any $v \in M_k$, we have 
\begin{align*}
(w\overline{w}^T).(v \otimes t^k) = \sum_{i=2}^{N}a_i^2\Omega_{k,e_1,e_i}^{\mathcal{H},\beta}.(v \otimes t^k) +  2\sum_{1 \leqslant i <j \leqslant N} a_ia_j \big(\mathcal{T}_{k,e_1,e_i}^{\mathcal{H},\beta}\overline{\mathcal{T}_{k,e_1,e_j}^{\mathcal{H},\beta}}^T + \mathcal{T}_{k,e_1,e_j}^{\mathcal{H},\beta}\overline{\mathcal{T}_{k,e_1,e_i}^{\mathcal{H},\beta}}^T \big).(v \otimes t^k) \\
= \sum_{i=2}^{N}a_i^2\Omega_{k,e_1,e_i}^{\mathcal{H},\beta}.(v \otimes t^k) + 2 \sum_{1 \leqslant i <j \leqslant N} a_ia_j \big(\Omega_{k,e_i,e_i+e_j}^{\mathcal{H},\beta} - \Omega_{k,e_1,e_i}^{\mathcal{H},\beta} - \Omega_{k,e_1,e_j}^{\mathcal{H},\beta}\big).(v \otimes t^k) 
\end{align*}
and thus we are done using Lemma \ref{Invariance}. \\
(3) is a direct consequence of (2), as for every $v \in M_k$, we have
\begin{align*}
(u_1\overline{u_2}^T + u_2\overline{u_1}^T).(v \otimes t^k) 
= \big((u_1 + u_2)(\overline{u_1+u_2})^T - (u_1\overline{u_1}^T) - (u_2\overline{u_2}^T) \big).(v \otimes t^k) \in M_k \otimes \mathbb{C}t^k.
\end{align*}
\end{proof}

\subsection{Indecomposable and completely reducible $\mathcal{H}_N$-modules} Recall that for any $k \in \mathbb{Z}^N$ with $k+\beta \neq 0$, we have a symplectic basis $\{k+\beta, w_1^{(k)}, w_2^{(k)}, \ldots, w_{N-1}^{(k)} \}$ of $\mathbb{C}^N$, which further contains a symplectic basis $\{w_i^{(k)} \}_{i=2}^{N-1}$ of $\mathbb{C}^{N-2}$. This gives rise to a Lie algebra of type $C_{n-1}$, with respect to this basis (see Remark \ref{basis}). We shall denote this symplectic Lie algebra by $\mathfrak{sp}_{N-2}^{(k)}$ and its Cartan subalgebra by $\mathfrak{h}_{\mathfrak{sp}_{N-2}}^{(k)}$, which is the abelian ad-diagonalizable subalgebra spanned by $\{w_2^{(k)}\overline{w_{3}^{(k)}}^T + w_{3}^{(k)}\overline{w_{2}^{(k)}}^T, \ w_{4}^{(k)}\overline{w_{5}^{(k)}}^T + w_{5}^{(k)}\overline{w_{4}^{(k)}}^T, \ldots , w_{N-2}^{(k)}\overline{w_{N-1}^{(k)}}^T + w_{N-1}^{(k)}\overline{w_{N-2}^{(k)}}^T \}$. It is again clear that $\mathfrak{sp}_{N-2}^{(k)}$ acts on $\mathbb{C}^{N-2}$, where $\mathbb{C}^{N-2} \cong \bigoplus_{i=2}^{N-1}\mathbb{C}w_i^{(k)}$ under the usual identification map preserving the skew-symmetric bilinear form. Set $H_i^{(k)} = w_{2i}^{(k)}\overline{w_{2i+1}^{(k)}}^T + w_{2i+1}^{(k)}\overline{w_{2i}^{(k)}}^T$ and $\epsilon_i^{(k)} \in (\mathfrak{h}_{\mathfrak{sp}_{N-2}}^{(k)})^*$ such that $\epsilon_i^{(k)}(H_j^{(k)}) = \delta_{j,i} \ \forall \ 1 \leqslant i \leqslant n-1$. This induces the weight space decomposition of $\mathbb{C}^{N-2}$ with respect to $\mathfrak{h}_{\mathfrak{sp}_{N-2}}^{(k)}$, given by  
$\mathbb{C}^{N-2} = \bigoplus_{i=2}^{n} (\mathbb{C}^{N-2})_{\pm \epsilon_i^{(k)}}$, where $(\mathbb{C}^{N-2})_{\epsilon_i^{(k)}} = \mathbb{C}w_{2i-1}^{(k)}$ and $(\mathbb{C}^{N-2})_{-\epsilon_i^{(k)}} = \mathbb{C}w_{2i-2}^{(k)} \ \forall \ 2 \leqslant i \leqslant n$. Also observe that $\mathfrak{sp}_{N-2}^{(k)}$ acts trivially on both $k+\beta$ and $w_1^{(k)}$. 

\blmma \label{Restricted}
Let $M = \oplus_{k \in \mathbb{Z}^N}(M_k \otimes \mathbb{C}t^k)$ be a non-trivial $\mathcal{H}_N$-submodule of $L_{\mathcal{H}}^{\beta}(V)$, where $V$ is an $\mathfrak{sp}_N$-module. Then $\mathfrak{sp}_{N-2}^{(k)}$ acts on $M_k$ for any arbitrary $k \in \mathbb{Z}^N$ satisfying $k+\beta \neq 0$, where $\mathfrak{sp}_{N-2}^{(k)} \cong \mathfrak{sp}\big(\{\overline{k+\beta} \}^{\perp}/\mathbb{C}(k + \beta)\big)$.
\elmma

\begin{proof}
First note that $\Omega_{k,r,s}^{\mathcal{H},\beta}$ acts trivially on the vector $k+\beta$. The lemma then follows by applying Proposition \ref{Prop}, since $(\overline{k+\beta}|w_i^{(k)}) = 0 \ \forall \ 2 \leqslant i \leqslant N-1$ and $\mathfrak{sp}_{N-2}^{(k)} = \text{span} \{r\overline{r}^T \ | \ r \in \mathbb{Z}^{N-2} \}$, where we identify $\mathbb{Z}^{N-2}$ with $\bigoplus_{i=2}^{N-1}\mathbb{Z}w_i^{(k)}$.  
\end{proof}
 
\bthm \label{Uniqueness}
If $1 \leqslant p \leqslant n$, then $W_{min}^{\beta}\big(V(\omega_p)\big)$ is the unique non-trivial irreducible submodule as well as the unique irreducible quotient of $L^{\beta}\big(V(\omega_p)\big)$.
\ethm

\begin{proof}
Let $M$ be an irreducible $\mathcal{H}_N$-submodule of $L_{\mathcal{H}}^{\beta}\big(V(\omega_p)\big)$. In view of Theorem \ref{Irreducible}, it suffices to prove that $W_{min}^{\beta}\big(V(\omega_p)\big) \cap M \neq (0)$. Now since $V(\omega_p)$ is not the trivial module, we can find $k \in \mathbb{Z}^N$ with $k+\beta \neq 0$ such that $M_k \neq (0)$. Let $0 \neq v  \in M_k$ be arbitrary. By Lemma \ref{Basis}, there exist $p_i, q_i, r_i, s_i \in \mathbb{C}$ such that $v = \sum_{i=1}^{a} p_iu_i^{(k)} + \sum_{i=1}^{b} q_ix_i^{(k)} + \sum_{i=1}^{c} r_iy_i^{(k)} + \sum_{i=1}^{d} s_iz_i^{(k)}$, where $\{u_1^{(k)}, \ldots, u_a^{(k)}, x_1^{(k)}, \ldots, x_b^{(k)}, y_1^{(k)}, \ldots, y_c^{(k)}, z_1^{(k)}, \ldots, z_d^{(k)} \}$ is a linearly independent subset of $V(\omega_p)$.\\
\noindent \textit{Case 1.} $s_i \neq 0$ for some $1 \leqslant i \leqslant d$. \\
By Proposition \ref{Prop}, we have $0 \neq (k+\beta)(\overline{k+\beta})^Tv  = \sum_{i=1}^{d} s_i(k+\beta)(\overline{k+\beta})^Tz_i^{(k)}\in M_k \cap W_{min}^{k,\beta}\big(V(\omega_p)\big)$. 

\noindent \textit{Case 2.} $s_i = 0 \ \forall \ 1 \leqslant i \leqslant d$ and $p_i \neq 0$ for some $1 \leqslant i \leqslant a$. \\
For $p=1, \ v = b_0(k + \beta) + \sum_{j}a_jw_{i_j}^{(k)}$, where $i_j \neq 1, \ b_0 \in \mathbb{C}^{\times}, \ a_j \in \mathbb{C}$ and the $\mathfrak{h}_{\mathfrak{sp}_{N-2}}^{(k)}$-weights of $w_{i_j}^{(k)}$ are all different, which implies that $\sum_{j}a_jw_{i_j}^{(k)} \in M_k$ and hence $b_0(k + \beta) \in M_k \cap W_{min}^{k,\beta}\big(V(\omega_p)\big)$. So let us assume that $2 \leqslant p \leqslant n$. In this case, it is easy to see that exactly $(p-1)$ elements of $\mathbb{C}^N$ occurring in the exterior product of any summand of $u_i^{(k)}$ have non-zero $\mathfrak{h}_{\mathfrak{sp}_{N-2}}^{(k)}$-weights, whereas each element of $\mathbb{C}^N$  in the exterior product of every summand of $x_i^{(k)}$ has non-zero $\mathfrak{h}_{\mathfrak{sp}_{N-2}}^{(k)}$-weight. Again, any summand of $y_i^{(k)}$ either has exactly $(p-2)$ elements of $\mathbb{C}^N$  with non-zero $\mathfrak{h}_{\mathfrak{sp}_{N-2}}^{(k)}$-weights in its exterior product or all the elements of $\mathbb{C}^N$ occuring in the exterior product of a summand of $y_i^{(k)}$  have non-zero $\mathfrak{h}_{\mathfrak{sp}_{N-2}}^{(k)}$-weights (see the discussion before Lemma \ref{Restricted}). It can be now easily deduced using elementary arguments that the $\mathfrak{h}_{\mathfrak{sp}_{N-2}}^{(k)}$-weight of $u_i^{(k)}$ differs from the $\mathfrak{h}_{\mathfrak{sp}_{N-2}}^{(k)}$-weights of both $x_i^{(k)}$ and $y_i^{(k)}$, from which we finally get $0 \neq \sum_{i=1}^{a} p_iu_i^{(k)} \in M_k \cap W_{min}^{k,\beta}\big(V(\omega_p)\big)$.  

\noindent \textit{Case 3.} $p_i = 0 \ \forall \ 1 \leqslant i \leqslant a$ and $s_i = 0 \ \forall \ 1 \leqslant i \leqslant d$. \\
For $p=1$, we have $v = \sum_{j}a_jw_{i_j}^{(k)}$, where $i_j \neq 1$ and $a_j \in \mathbb{C}$ for all $j$. Without loss of generality, let $i_1 = \mathrm{min} \{i_j\}_j$ and consider $s = i_1-1$ or $s = i_1+1$,  depending on whether $(\overline{w_{i_1-1}^{(k)}}|w_{i_1}^{(k)}) \neq 0$ or $(\overline{w_{i_1}^{(k)}}|w_{i_1+1}^{(k)}) \neq 0$. This unravels that $0 \neq \big((k + \beta)\overline{w_s^{(k)}}^T + w_s^{(k)}(\overline{k +\beta})^T\big)v \in M_k\cap W_{min}^{k,\beta}\big(V(\omega_p)\big)$ by Lemma \ref{Similar} and Proposition \ref{Prop}, as required. So let us suppose that $2 \leqslant p \leqslant n$. In this case, $v = \sum_{i=1}^{b} q_ix_i^{(k)} + \sum_{i=1}^{c} r_iy_i^{(k)}$. Consider the minimum of all the indices of the elements $w_{i_j}^{(k)}$ of $\mathbb{C}^N$, occuring in the exterior product of every summand of $v$ that does not contain $k+\beta$. For example, if $v = w_2^{(k)}\wedge w_4^{(k)} +(k+\beta \wedge w_1^{(k)} - w_2^{(k)} \wedge w_3^{(k)})$, then the indices present in the exterior products are $2,3,4$ and so their minimum will be $2$. Denoting this minimum by $i_l$, let us take $m= i_l+1$or $m= i_l-1$, depending on whether$(\overline{w_{i_1-1}^{(k)}}|w_{i_1}^{(k)}) \neq 0$ or $(\overline{w_{i_1}^{(k)}}|w_{i_1+1}^{(k)}) \neq 0$. If $r_i=0 \ \forall \ i$ or $w_m$ does not occur in the exterior product of any summand containing $k+\beta$, then it is easy to see that $\big((k + \beta)\overline{w_m^{(k)}}^T + w_m^{(k)}(\overline{k +\beta})^T\big)v$ is non-zero. On the other hand, if $r_i \neq 0$ for some $i$ and $w_m$ occurs in the exterior product of some summand containing $k+\beta$, then it is not difficult to deduce that $\big((k + \beta)\overline{w_m^{(k)}}^T + w_m^{(k)}(\overline{k +\beta})^T\big)v \neq 0$, using the fact that $y_i^{(k)}$ belongs to $V(\omega_p)$.
The result now follows verbatim as in the $p=1$ setting, since $\big((k + \beta)\overline{w_m^{(k)}}^T + w_m^{(k)}(\overline{k +\beta})^T\big)v \in M_k\cap W_{min}^{k,\beta}\big(V(\omega_p)\big)$.
 
Let $M^{\prime}$ be a proper non-zero $\mathcal{H}_N$-submodule of $L_{\mathcal{H}}^{\beta}\big(V(\omega_p)\big)$. To prove that $L_{\mathcal{H}}^{\beta}\big(V(\omega_p)\big)$ has a unique irreducible quotient, we shall show that $M^{\prime}$ is contained in $W_{max}\big(V(\omega_p)\big)$  if $\beta \notin \mathbb{Z}^N$ and in $\widehat{W}_{max}\big(V(\omega_p)\big)$ if $\beta \in \mathbb{Z}^N$. Pick $0 \neq v \in N_k$, where $k+\beta \neq 0$. By Lemma \ref{Basis}, there exist $p_i, q_i, r_i, s_i \in \mathbb{C}$ such that $v = \sum_{i=1}^{a} p_iu_i^{(k)} + \sum_{i=1}^{b} q_ix_i^{(k)} + \sum_{i=1}^{c} r_iy_i^{(k)} + \sum_{i=1}^{d} s_iz_i^{(k)}$. Moreover, as $W_{min}^{\beta}\big(V(\omega_p)\big) \subseteq M^{\prime}$, we are reduced to $v = \sum_{i=1}^{b} q_ix_i^{(k)} + \sum_{i=1}^{c} r_iy_i^{(k)} + \sum_{i=1}^{d} s_iz_i^{(k)}$. \\
\textbf{Claim.} $s_i = 0 \ \forall \ 1 \leqslant i \leqslant d$. \\
If not, then there exists $1 \leqslant i \leqslant d$ with $s_i \neq 0$. For $p =1$, it is evident that $v = \sum_{j}a_jw_{i_j}^{(k)} + c_0w_1^{(k)}$, where $i_j \neq 1$ and $c_0 \in \mathbb{C}^{\times}, a_j \in \mathbb{C}$ for all $j$.  
Now since the $\mathfrak{h}_{\mathfrak{sp}_{N-2}}^{(k)}$-weights of $w_{i_j}^{(k)}$ differ from each other, we have $\sum_{j}a_jw_{i_j}^{(k)} \in M^{\prime}_k$ and thus $0 \neq c_0w_1^{(k)} \in M^{\prime}_k$, whence it follows that $M^{\prime} = L_{\mathcal{H}}^{\beta}\big(V(\omega_1)\big)$ by Proposition \ref{Proposition}. This contradiction establishes our claim for $p=1$. For $2 \leqslant p \leqslant n$, consider $v = \sum_{i=1}^{b} q_ix_i^{(k)} + \sum_{i=1}^{c} r_iy_i^{(k)} + \sum_{i=1}^{d} s_iz_i^{(k)} \in M^{\prime}_k$. It is also clear that exactly $(p-1)$ elements of $\mathbb{C}^N$ occuring in the exterior product of any summand of $z_i^{(k)}$ have non-zero $\mathfrak{h}_{\mathfrak{sp}_{N-2}}^{(k)}$-weights. Then it is easy to deduce that $0 \neq \sum_{i=1}^{d} s_iz_i^{(k)} \in M^{\prime}_k$, by means of $\mathfrak{h}_{\mathfrak{sp}_{N-2}}^{(k)}$-weight arguments, as in \textit{Case 2}. This finally yields that $M^{\prime} = L_{\mathcal{H}}^{\beta}\big(V(\omega_p)\big)$ by Proposition \ref{Proposition}, which is a contradiction. Hence the claim.   
Consequently $L_{\mathcal{H}}^{\beta}\big(V(\omega_p)\big)$ admits a unique irreducible $\mathcal{H}_N$-quotient, given by $W_{min}\big(V(\omega_p)\big)$.
\end{proof}

\bcrlre
\
\begin{enumerate}
\item $L_{\mathcal{H}}^{\beta}\big(V(\omega_p)\big)$ is an indecomposable $\mathcal{H}_N$-module for each $1 \leqslant p \leqslant n$. 
\item Both $W_{int}^{\beta}\big(V(\omega_p)\big)$ and $W_{max}^{\beta}\big(V(\omega_p)\big)$ are indecomposable $\mathcal{H}_N$-modules for each $1 \leqslant p \leqslant n$. 
\end{enumerate}
\ecrlre

\begin{proof}
(1) and (2) follow by proceeding similarly as in Theorem \ref{Uniqueness}.
\end{proof}


The next lemma is crucial for showing that upto trivial modules, there is exactly one submodule sitting between the minimal and maximal submodules for $1 < p < n$.
\blmma \label{L6.6}
Let $M$ be an $\mathcal{H}_N$-submodule of $W_{max}^{\beta}\big(V(\omega_p)\big)$, where $1 < p < n$. Then:
\begin{enumerate}
\item For every $k, r, s \in \mathbb{Z}^N, \ \Omega_{k,r,s}^{\mathcal{H},\beta}.v  \in M_k \cap W_{int}^{k,\beta}\big(V(\omega_p)\big) \ \forall \ v \in M_k$.
\item If $\beta \notin \mathbb{Z}^N$ and $W_{min}^{\beta}\big(V(\omega_p)\big)$ is a proper subset of $M$, then $M \cap W_{int}^{\beta}\big(V(\omega_p)\big) \neq W_{min}^{\beta}\big(V(\omega_p)\big)$.  
\item  If $\beta \in \mathbb{Z}^N$ and $\widehat{W}_{min}^{\beta}\big(V(\omega_p)\big)$ is a proper subset of $M$, then $M \cap \widehat{W}_{int}^{\beta}\big(V(\omega_p)\big) \neq \widehat{W}_{min}^{\beta}\big(V(\omega_p)\big)$.   
\end{enumerate}
\elmma
\begin{proof}
(1) For any $v \in M_k$, it is easy to see that
 \begin{equation*}
  \theta_{p+1}^{\beta}\big((k+\beta)\wedge \Omega_{k,r,s}^{\mathcal{H},\beta}.v\big) =
    \begin{cases}
      (k+\beta)(\overline{k+\beta})^T(\Omega_{k,r,s}^{\mathcal{H},\beta}.v \big) + \theta_{p}^{\beta}\big(\Omega_{k,r,s}^{\mathcal{H},\beta}.v \big)(k+\beta), & p=2 \\
      (k+\beta)(\overline{k+\beta})^T(\Omega_{k,r,s}^{\mathcal{H},\beta}.v \big) 
      + (k+\beta) \wedge \theta_{p}^{\beta}\big(\Omega_{k,r,s}^{\mathcal{H},\beta}.v\big), & 2 < p < n.
    \end{cases}       
\end{equation*}
Then (1) follows using Lemma \ref{Lemma} and Lemma \ref{Invariance}, as $v \in W_{max}^{k,\beta}\big(V(\omega_p)\big)$.\\ 
To complete the proof of the lemma, it is enough to establish (2), as (3) can be handled in a more or less similar manner. Let us fix any $k \in \mathbb{Z}^N$ with $M_k \neq (0)$. We proceed by considering 2 cases.\\
\textit{Case 1.} $p = 2$. \\
If $M + W_{int}^{\beta}\big(V(\omega_p)\big) = M$ or $M + W_{int}^{\beta}\big(V(\omega_p)\big) = W_{int}^{\beta}\big(V(\omega_p)\big)$, then the result follows from Theorem \ref{Composition}. So it suffices to consider $M + W_{int}^{\beta}\big(V(\omega_p)\big) = W_{max}^{\beta}\big(V(\omega_p)\big)$, which implies that $\{k+\beta \wedge w_1^{(k)} - w_{2i-2}^{(k)} \wedge w_{2i-1}^{(k)} \ | \ 2 \leqslant i \leqslant N \} \subseteq M_k + W_{int}^{k,\beta}\big(V(\omega_p)\big)$. Therefore, we can find $v = \sum_{i=1}^{a}p_iu_i^{(k)} + \sum_{i=1}^{b}q_ix_i^{(k)} + \sum_{i=2}^{n}r_{2i-2,2i-1}(k+\beta \wedge w_1^{(k)} - w_{2i-2}^{(k)} \wedge w_{2i-1}^{(k)}) \in M_k$, where $r_{2i-2,2i-1} = 0$ for some $2 \leqslant i \leqslant N$, with either $\sum_{i=2}^{n}r_{2i-2,2i-1}(k+\beta \wedge w_1^{(k)} - w_{2i-2}^{(k)} \wedge w_{2i-1}^{(k)}) \neq 0$ or $\sum_{i=1}^{b}q_ix_i^{(k)} \neq 0$ (see Lemma \ref{Basis}). If $\sum_{i=2}^{n}r_{2i-2,2i-1}(k+\beta \wedge w_1^{(k)} - w_{2i-2}^{(k)} \wedge w_{2i-1}^{(k)}) = 0$, then $\sum_{i=1}^{b}q_ix_i^{(k)} \in M_k$, thanks to Theorem \ref{Uniqueness}, which gives us the desired result. Now suppose that $\sum_{i=2}^{n}r_{2i-2,2i-1}(k+\beta \wedge w_1^{(k)} - w_{2i-2}^{(k)} \wedge w_{2i-1}^{(k)}) \neq 0$ and $r_{2j-2,2j-1} = 0$ for some $2 \leqslant j \leqslant N$. Without loss of generality, let us take $r_{2,3} \neq 0$. Consequently, if $w_{2j-2}^{(k)} \wedge w_{2j-1}^{(k)} \neq x_i^{(k)}$ for any $1 \leqslant i \leqslant b$, then $(w_{2j-1}^{(k)}\overline{w_2^{(k)}}^T + w_2^{(k)}\overline{w_{2j-1}^{(k)}}^T)v = -r_{2,3}(w_2^{(k)} \wedge w_{2j-1}^{(k)}) + \ldots$ is a non-zero element of $M_k$, which proves the lemma. On the other hand, if $w_{2j-2}^{(k)} \wedge w_{2j-1}^{(k)}= x_i^{(k)}$ for some $1 \leqslant i \leqslant b$ , then $w_{2m-2}^{(k)} \wedge w_{2m-1}^{(k)}$ must also appear in one of the summands of $\sum_{i=1}^{b}q_ix_i^{(k)}$ for some $2 \leqslant m \leqslant N$ satisfying $r_{2m-2,2m-1} = 0$, which yields that $0 \neq (w_{2m-2}^{(k)}\overline{w_{2j-2}^{(k)}}^T + w_{2j-2}^{(k)}\overline{w_{2m-2}^{(k)}}^T)v \notin W_{min}^{k,\beta}\big(V(\omega_p)\big)$ and thus we are done using Proposition \ref{Prop}, Theorem \ref{Uniqueness} and (1). \\      
\textit{Case 2.} $2 < p <n$. \\
Let, if possible, $M \cap W_{int}^{\beta}\big(V(\omega_p)\big) = W_{min}^{\beta}\big(V(\omega_p)\big)$. Consider the $\mathcal{H}_N$-module map $T_p^{\beta}|_M$. From Lemma \ref{Lemma}, it is evident that $\mathrm{Ker}(T_p^{\beta}|_M) = M \cap W_{int}^{\beta}\big(V(\omega_p)\big)$, which implies that $T_p^{\beta}|_M$ is a  non-zero map, due to our hypothesis. But again, since $W_{min}^{\beta}\big(V(\omega_{p-1})\big)$ is irreducible, we  obtain $\mathrm{Im}(T_p^{\beta}|_M) = W_{min}^{\beta}\big(V(\omega_{p-1})\big)$, as $\mathrm{Im}(T_p^{\beta}|_{W_{max}^{\beta}\big(V(\omega_p)\big)}) = W_{min}^{\beta}\big(V(\omega_{p-1})\big)$ by Lemma \ref{Lemma}.\\ 
\textbf{Claim.} There exists $v \in M_k$ such that $\Omega_{k,r,s}^{\mathcal{H},\beta}.v \notin W_{min}^{k,\beta}\big(V(\omega_p)\big)$ for some $r,s \in \mathbb{Z}^N$. \\
If not, then $\Omega_{k,r,s}^{\mathcal{H},\beta}M_k \subseteq W_{min}^{k,\beta}\big(V(\omega_p)\big) \ \forall \ r,s \in \mathbb{Z}^N$. Now since $0 \neq v^{\prime} = (k+\beta) \wedge w_2^{(k)} \wedge w_4^{(k)} \wedge \ldots \wedge w_{2(p-3)}^{(k)} \wedge w_{2(p-2)}^{(k)} \in W_{min}^{\beta}\big(V(\omega_{p-1})\big)$, there exists $x \in M_k$ such that $T_p(x \otimes t^k) =  v^{\prime} \otimes t^k$. This implies that  $T_p\bigg(\big(w_{2p-3}^{(k)}\overline{w_{2p-3}^{(k)}}^T\big)(x \otimes t^k)\bigg) = 0$. But on the other hand, we have $\big(w_{2p-3}^{(k)}\overline{w_{2p-3}^{(k)}}^T\big)T_p(x \otimes t^k) \neq 0$, which contradicts Lemma \ref{Invariance}. Hence the claim. \\ 
The assertion now follows from Lemma \ref{Invariance}, Theorem \ref{Uniqueness} and (1).
\end{proof}

\bthm \label{Main}
Let $W$ be a non-zero proper $\mathcal{H}_N$-submodule of $L_{\mathcal{H}}^{\beta}\big(V(\omega_p)\big)$, where $0 \leqslant p \leqslant n$.
\begin{enumerate}
\item If $p=0$ and $\beta \in \mathbb{Z}^N$, then $W$ is isomorphic to either $\bigoplus_{k \in \mathbb{Z}^N, \ k \neq -\beta} \mathbb{C}t^k$ or $\mathbb{C}t^{-\beta}$. 
\item If $1 < p < n$ and $\beta \notin \mathbb{Z}^N$, then $W$ is isomorphic to either $W_{min}^{\beta}\big(V(\omega_p)\big)$ or $W_{int}^{\beta}\big(V(\omega_p)\big)$ or $W_{max}^{\beta}\big(V(\omega_p)\big)$.
\item If $1 < p < n$ and $\beta \in \mathbb{Z}^N$, then $W$ is isomorphic to either $W_{min}^{\beta}\big(V(\omega_p)\big) \oplus (W_1 \otimes \mathbb{C}t^{-\beta})$ or $W_{int}^{\beta}\big(V(\omega_p)\big) \oplus (W_2 \otimes \mathbb{C}t^{-\beta})$ or $W_{max}^{\beta}\big(V(\omega_p)\big)$, where $W_1$ and $W_2$ are arbitrary subspaces (including the zero subspace) of $V(\omega_p)$.
\item If $p \in \{1,n\}$ and $\beta \notin \mathbb{Z}^N$, then $W$ is isomorphic to either $W_{min}^{\beta}\big(V(\omega_p)\big)$ or $W_{max}^{\beta}\big(V(\omega_p)\big)$.
\item If $p \in \{1,n\}$ and $\beta \in \mathbb{Z}^N$, then $W$ is isomorphic to either $W_{min}^{\beta}\big(V(\omega_p)\big) \oplus (W_1 \otimes \mathbb{C}t^{-\beta})$ or  $W_{max}^{\beta}\big(V(\omega_p)\big)$, where $W_1$ is an arbitrary subspace (including the zero subspace) of $V(\omega_p)$. 
\end{enumerate}
\ethm

\begin{proof}
(1) is an immediate consequence of Corollary \ref{Cor}. \\
(2), (3), (4), (5) follow by simply combining Theorem \ref{Composition}, Theorem \ref{Uniqueness} and Lemma \ref{L6.6}.
\end{proof}

\bcrlre\label{C6.8}
\
\begin{enumerate}
\item For $1 \leqslant p \leqslant n$ and $\beta \notin \mathbb{Z}^{N}, \ \mathrm{Soc}L_{\mathcal{H}}^{\beta}\big(V(\omega_p)\big) \cong \mathrm{Head}L_{\mathcal{H}}^{\beta}\big(V(\omega_p)\big) = W_{min}^{\beta}\big(V(\omega_p)\big)$.
\item For $1 \leqslant p \leqslant n$ and $\beta \notin \mathbb{Z}^{N}, \ L_{\mathcal{H}}^{\beta}\big(V(\omega_p)\big)$ is uniserial.
\item For $\beta \in \mathbb{Z}^{N}, \  L_{\mathcal{H}}^{\beta}\big(V(\omega_0)\big)$ is completely reducible, admitting exactly two Loewy series with $\mathrm{Soc}L_{\mathcal{H}}^{\beta}\big(V(\omega_0)\big) \cong \mathrm{Head}L_{\mathcal{H}}^{\beta}\big(V(\omega_0)\big) = L_{\mathcal{H}}^{\beta}\big(V(\omega_0)\big)$.
\end{enumerate}
\ecrlre

\begin{proof}
Follows from Corollary \ref{Cor}, Theorem \ref{Composition} and Theorem \ref{Main}.
\end{proof}

\section{Differential and harmonic forms, homology, EALAs}\label{S7}
In this section, we first show that the submodules of the exceptional Shen--Larsson modules over $\mathcal{H}_N$ come from kernels and images of a suitable mixed complex of the modules of differential forms, which are also known as the modules of tensor fields on a torus \cite{F,R,GS}. We then provide a Hodge-theoretic interpretation of our results and compute various (co-)homologies concerning these modules before concluding this section by finally relating the differential $1$-forms to the construction of Hamiltonian extended affine Lie algebras.

\subsection{Mixed Complex arising from Modules of Tensor Fields} \label{SS7.1}
The module of tensor fields $\Omega^p(\beta) := F^{\beta}(\Lambda^p\mathbb{C}^N), \ 1 \leqslant p \leqslant N$, can be interpreted as $p$-differential forms via the correspondence
\begin{align*}
e_{i_1} \wedge \ldots \wedge e_{i_p} \otimes t^r \longmapsto t^rt_{i_1}^{-1}dt_{i_1} \wedge \ldots \wedge t_{i_p}^{-1}dt_{i_p}, \ r \in \mathbb{Z}^N.
\end{align*}
It is well-known that these modules of differential forms also admit an associative $A_N$-action \cite{YB,R3}. Moreover, they also form two complexes given by 
\begin{align*}
\Omega^0(\beta) \xrightarrow{\pi_0^{\beta}} \Omega^1(\beta) \xrightarrow{\pi_1^{\beta}} \ldots \ldots \xrightarrow{\pi_{N-1}^{\beta}} \Omega^N(\beta) \\
\Omega^N(\beta) \xrightarrow{T_N^{\beta}} \Omega^{N-1}(\beta) \xrightarrow{T_{N-1}^{\beta}} \ldots \ldots \xrightarrow{T_{1}^{\beta}} \Omega^0(\beta), 
\end{align*}
where $\Omega^0(\beta) = A_N$. Here $\pi_p^{\beta}$ and $T_p^{\beta}$ are $\mathcal{H}_N$-module complexes (but not $A_N$-module complexes) and so their kernels and images naturally give rise to $\mathcal{H}_{N}$-submodules of $\Omega^p(\beta)$. Then $\big(\Omega^p(\beta), \pi_p^{\beta}\ \big)$ is a \textit{de Rham type complex} and $\big(\Omega^p(\beta), T_p^{\beta}\ \big)$ is a \textit{Koszul type complex}.

Moreover, it is easy to see that $T_{p+1}^{\beta} \circ \pi_{p}^{\beta} + \pi_{p-1}^{\beta} \circ T_p^{\beta} = 0 \ \forall \ 1 \leqslant p \leqslant N-1$, which implies that $\{ \pi^{\beta}, T^{\beta} \} = 0$, where $\{\cdot \}$ denotes the \textit{anti-commutator} of $\pi^{\beta}$ and $T^{\beta}$. This shows that $\big(\Omega(\beta),T^{\beta},\pi^{\beta})$ has the structure of a \textit{mixed complex}, in the sense of Kassel \cite{Ka}. 

Following \cite{KNS}, we now define the \textit{cyclic complex} associated with the above mixed complex by 
\begin{align*}
\Omega_{\text{cyc}}^{p}(\beta) := \bigoplus_{i \geqslant 0}\Omega^{p-2i}(\beta), \ \text{where we take} \ \Omega^{k}(\beta) := 0 \ \forall \ k<0.
\end{align*}
with the differential $d_p^{\beta} := \pi^{\beta} + T^{\beta}$. Then the homology of the mixed complex is given by
\begin{align*}
\mathbb{HC}^{p}(\beta) = \dfrac{\mathrm{Ker}(d_{p}^{\beta}: \Omega_{\text{cyc}}^{p}(\beta) \longrightarrow \Omega_{\text{cyc}}^{p+1}(\beta))}{\mathrm{Im}(d_p^{\beta}: \Omega_{\text{cyc}}^{p-1}(\beta) \longrightarrow \Omega_{\text{cyc}}^{p}(\beta))}
\end{align*}

\bppsn \label{de Rham}
Let $\beta \in \mathbb{C}^N$.
\
\begin{enumerate}
\item If $\beta \notin \mathbb{Z}^N$, then the de Rham type cohomology $\mathbb{H}^{p}_{\mathrm{dR}}(\beta)$ of $\{\Omega^p(\beta), \pi_p^{\beta}\}$ as well as the Koszul type homology $\mathbb{H}^{p}_{\mathrm{Kos}}(\beta)$ of $\{\Omega^p(\beta),  T_p^{\beta} \}$ vanish.  
\item If $\beta \in \mathbb{Z}^N$, then the de Rham type cohomology $\mathbb{H}^{p}_{\mathrm{dR}}(\beta)$ of $\{\Omega^p(\beta), \pi_p^{\beta}\}$ as well as the Koszul type homology $\mathbb{H}^{p}_{\mathrm{Kos}}(\beta)$ of $\{\Omega^p(\beta),  T_p^{\beta} \}$ are given by the trivial module $\mathbb{C}^{{N \choose p}}$.
\item $\mathbb{HC}^{p}(\beta) = (0)$ for $\beta \notin \mathbb{Z}^N$, otherwise it is the trivial module $\mathbb{C}^{{N \choose p}}$.
\end{enumerate}
\eppsn

\begin{proof}
For (1) and (2), it suffices to show that $\mathrm{Ker}T_p^{k, \beta} = W_{int}^{k,\beta}(\Lambda^p\mathbb{C}^N)$ for any $k \in \mathbb{Z}^N$ satisfying $k+\beta \neq 0$. To prove this, first observe that
\begin{align*}
T_p^{k,\beta}\big((w_{i_1}^{(k)} \wedge \ldots \wedge w_{i_p}^{(k)}) \otimes t^k \big) =T_p^{k,\beta}\big((k +\beta \wedge w_{j_1}^{(k)} \wedge \ldots \wedge w_{j_{p-1}}^{(k)}) \otimes t^k \big) = 0, 
\end{align*}
where $i_r, j_s \neq 1$.
But again, we have $T_p^{k,\beta}\big((w_1^{(k)} \wedge w_{i_1}^{(k)} \wedge \ldots \wedge w_{i_{p-1}}^{(k)}) \otimes t^k \big) \neq 0 \ (i_r \neq 1)$, along with $W_{int}^{k,\beta}(\Lambda^p\mathbb{C}^N) = \mathrm{span} \{w_{i_1}^{(k)} \wedge \ldots \wedge w_{i_p}^{(k)}, \ k +\beta \wedge w_{j_1}^{(k)} \wedge \ldots \wedge w_{j_{p-1}}^{(k)} \ | \ i_r, j_s \neq 1 \}$ and so we are done.   \\
(3) follows by a similar computation. 
\end{proof}

\brmk \label{R7.2}
The minimal, intermediate and maximal submodules of $L_{\mathcal{H}}^{\beta}\big(V(\omega_p)\big)$ are obtained by restricting the kernels of $\pi_p^{\beta}, \ T_p^{\beta}$ and $T_{p+1}^{\beta} \circ \pi_p^{\beta}$ respectively to $L_{\mathcal{H}}^{\beta}\big(V(\omega_p)\big)$.
\ermk   

\begin{subsection}{A Hodge-type Interpretation of the Mixed Complex}\label{Hodge}
The mixed complex structure ($\Omega(\beta),\pi^{\beta},T^{\beta}$) admits a natural interpretation analogous to the ($d,\delta$)-formalism that appears in symplectic Hodge theory. In other words, the pair ($\pi^{\beta},T^{\beta}$) formally resembles the pair ($d,\delta$) consisting of the de Rham differential and a codifferential-type operator. However, unlike the classical geometric setting, the operators $\pi^{\beta}$ and $T^{\beta}$ are homogeneous lattice operators acting fiberwise on Fourier components rather than global differential operators on a manifold. 

The mixed complex also satisfies the identities
\begin{align*}
\mathrm{Ker}\pi_p^{\beta} \cap \mathrm{Im}T_{p+1}^{\beta} = W_{min}^{\beta}(\Lambda^p\mathbb{C}^N) = \mathrm{Im}(\pi_{p-1}^{\beta} \circ T_{p}^{\beta}) \  \forall \ 1 \leqslant p <N,
\end{align*}
which may be viewed as a formal algebraic version of the classical $d\delta$\textit{-lemma} in Hodge theory, where the intersection of \textit{closed} and \textit{co-exact} forms are described by mixed differential operators \cite{M}.

Motivated by the above analogy, we define the space of \textit{symplectically harmonic} $p$-forms as
\begin{align*}
\Omega_{\text{har}}^{p} (\beta) := \mathrm{Ker}\pi_{p}^{\beta} \cap \mathrm{Ker}T_{p}^{\beta}
\end{align*}
and the associated \textit{harmonic cohomology} is then given by
\begin{align}\label{HC}
\mathbb{H}_{\text{har}}^{p}(\beta) := \dfrac{\Omega_{\text{har}}^{p} (\beta)}{\Omega_{\text{har}}^{p}(\beta)\cap\mathrm{Im}(\pi_{p-1}^{\beta})}.
\end{align}
By Proposition \ref{de Rham}, it follows that $\mathbb{H}_{\text{har}}^{p}(\beta) = \mathbb{H}_{\text{dR}}^{p}(\beta)$, which is a characteristic trait of compact symplectic manifolds satisfying the $d\delta$-lemma \cite{Ma,M,Y}. Thus, the mixed complex exhibits a Hodge-type dichotomy analogous to the vanishing and nonvanishing of twisted cohomology on $\mathbb{T}^N$.

The submodules $L_{\mathcal{H}}^{\beta}\big(V(\omega_p)\big)$ of modules of tensor fields can be viewed as modules of \textit{primitive} $p$-form fields on $\mathbb{T}^N$. Although $\pi^{\beta}$ and $T^{\beta}$ do not preserve these primitive $p$-form fields individually, their harmonic intersection naturally interacts with it, in the sense that
  \begin{equation*}
  \Omega^{p}_{\text{har}}(\beta) \cap L_{\mathcal{H}}^{\beta}\big(V(\omega_p)\big) =
    \begin{cases}
      W_{min}^{\beta}\big(V(\omega_p)\big) & \text{if} \ \beta \notin \mathbb{Z}^N\\
      W_{min}^{\beta}\big(V(\omega_p)\big) \oplus \mathbb{C}^{{N \choose p}} & \text{if} \ \beta \in \mathbb{Z}^N,
    \end{cases}       
  \end{equation*}
and so this harmonic primitive sector is \textit{irreducible} when $\beta \notin \mathbb{Z}^N$ (by Theorem \ref{Irreducible}).
\subsection{Realization of Hamiltonian Lie algebras} \label{Realization}
In the space of differential $1$-forms $\Omega^1(\beta)$, 
\begin{align*}
W_{min}^{\beta}\big(V(\omega_1)\big) =  \mathrm{Im}\pi_0^{\beta} \subseteq \mathrm{Ker}T_1^{\beta} = \widehat{W}_{max}^{\beta}\big(V(\omega_1)\big) \ (\text{by Lemma} \ \ref{Equal}),
\end{align*}
which shows that every \textit{exact} $1$-form is automatically \textit{co-closed} in the present setup. 

Again, using Lemma \ref{JH}, Remark \ref{Isom} and Lemma \ref{Equal}, we obtain 
\begin{align*}
\mathcal{H}_N^{\prime} \cong \Omega^{0}(0)/\{\text{constants}\} \cong W_{min}^{0}\big(V(\omega_1)\big) = \{{\text{exact 1-forms}}\}  \cong \Omega^1(0)/\{{\text{co-closed forms}}\}, \\
 \mathcal{H}_N \cong \widehat{W}_{min}^{0}\big(V(\omega_1)\big) = \{{\text{closed 1-forms}}\} \cong \Omega^1(0)/\{{\text{co-exact forms}}\}.
\end{align*}
\end{subsection}


\subsection{Abelian extension of $\mathcal{H}_N$} \label{EALAs}
The differential $1$-forms $\Omega^1(0)$ are intimately connected to the construction of \textit{Hamiltonian extended affine Lie algebras} (HEALAs), which were recently introduced in \cite{R2}. The realization of HEALAs involves a certain abelian extension of $\mathcal{H}_N$ and this extension (modulo $\mathcal{H}_N$), in fact, comes from $\Omega^1(0)$. More specifically, we know that $\mathcal{H}_N$ acts on $\Omega^1(0)$ and 
\begin{align*}
(0) \subseteq W_{min}^{0}\big(V(\omega_{1})\big) \subseteq W_{max}^{0}\big(V(\omega_{1})\big) \subseteq  L_{\mathcal{H}}^{0}\big(V(\omega_1)\big) = \Omega^1(0) \ \text{(see Theorem \ref{Composition})}. 
\end{align*} 
In this case, it  is easy to see that $\pi_{0}^{0}(A_N) = W_{min}^{0}\big(V(\omega_{1})\big)$ is the (irreducible) submodule of \textit{exact} $1$-forms and $T_2^{0}\big(\Omega^2(\beta)\big) = W_{max}^{0}\big(V(\omega_{1})\big)$ is the submodule of \textit{co-exact} 1-forms over $\mathcal{H}_N$. Moreover, $L_{\mathcal{H}}^{0}\big(V(\omega_1)\big) / W_{max}^{0}\big(V(\omega_{1})\big)$ (which is isomorphic to $\mathcal{H}_N$) happens to be the abelian extension of $\mathcal{H}_N$ used in constructing HEALAs (see \cite{R2} for more details on HEALAs).

\section{Irreducibility of Shen--Larsson modules over $\mathcal{W}_N$ and $\mathcal{S}_N$} \label{S8}
In this section, we provide simpler conceptual proofs of some well-known irreducibilty results for Shen--Larsson modules over $\mathcal{W}_N$ and $\mathcal{S}_N$, by resorting to a similar strategy that we have already utilized in Theorem \ref{Criterion}. Throughout this section, let us fix any $N >1$.  

\blmma \label{LW}
Let $V$ be an irreducible $\mathfrak{gl}_N$-module (not necessarily a weight module).  Then:  
\begin{enumerate}
	\item $\mathfrak{gl}_N =  \mathrm{span}\{ru^T \ | \ r \in \mathbb{Z}^N,  \ u \in \mathbb{C}^N \}$. 
	\item $\mathcal{J}_{\mathcal{W}}(V) := \{v \in V \ | \ (ru^{T})^2v = (u|r)(ru^T)v \ \forall \ u \in \mathbb{C}^N, \ r \in \mathbb{Z}^N \}$ is a $\mathfrak{gl}_N$-submodule of $V$.
	\item $\mathcal{J}_{\mathcal{W}}(V) = V$ if and only if $V \cong V(\delta_k,k)$ for some $0 \leqslant k \leqslant N$, where $\delta_0 = \delta_N = 0$.
\end{enumerate}
\elmma

\begin{proof} 
	(1) follows by simply observing that $E_{i,i} - E_{j,j} = (e_i + e_j)(e_i - e_j)^T + e_ie_j^T - e_je_i^T$ and $E_{i,j} = e_ie_j^T \ \forall \ 1 \leqslant i \neq j \leqslant N$ and $E_{i,i} = e_ie_i^T \ \forall \ 1 \leqslant i \leqslant N$.\\
	(2) can be deduced using Lemma \ref{L2}, by providing a similar argument as  in Lemma \ref{L3}. \\
	(3) If $V \cong V(\delta_k,k)$ for some $0 \leqslant k \leqslant N$, then a direct checking readily yields $\mathcal{J}_{\mathcal{W}}(V) = V$. \\
	Conversely, assume that $\mathcal{J}_{\mathcal{W}}(V) = V$. The irreducibility of $V$ thereby implies that $V$ is countable-dimensional and so by Dixmier's lemma, the central elements (which are scalar multiples of $\mathrm{Id}$) act by fixed scalars on $V$. Consequently, $V$ remains irreducible over $\mathfrak{sl}_N$. Now putting $r=e_i$ and $u=e_j$ for any $1 \leqslant i \neq j \leqslant N$, we get $(E_{ij})^2v = 0 \ \forall \ v \in V$. So we can directly invoke Lemma \ref{L1} to conclude that $V$ is finite-dimensional. Then by the representation theory of $\mathfrak{gl}_N$, we have: \\
	(i) $V \cong V(\overline{\lambda}) \cong V(\lambda, \overline{\lambda}(\mathrm{Id}))$ as $\mathfrak{gl}_N$-modules for some $\overline{\lambda} \in P_{\mathfrak{gl}_N}^{+}$, where $\lambda = \overline{\lambda}|_{\mathfrak{h}_{\mathfrak{sl}_N}^*} \in P_{\mathfrak{sl}_N}^{+}$. \\ 
	(ii) There exists $0 \neq v_{\overline{\lambda}} \in V(\overline{\lambda})$ such that $E_{i,i}v_{\overline{\lambda}} = \overline{\lambda}(E_{i,i})v_{\overline{\lambda}} \ \forall \ 1 \leqslant i \leqslant N$. \\
	Now taking $r=u=e_i$ for any $1 \leqslant i \leqslant N$, we obtain $(E_{i,i})^2v = E_{i,i}v \ \forall \ v \in V(\overline{\lambda})$. This implies that $(\overline{\lambda}(E_{i,i})^2 - \overline{\lambda}(E_{i,i}))v_{\overline{\lambda}} = 0$, which gives $\overline{\lambda}(E_{i,i}) \in \{0, 1 \} \ \forall \ 1 \leqslant i \leqslant N$. But again since $\overline{\lambda}(E_{i,i}) - \overline{\lambda}(E_{i+1,i+1}) \in \mathbb{Z}_{\geqslant 0} \ \forall \ 1 \leqslant i \leqslant N-1$, it is clear that $\overline{\lambda} = \overline{\delta_k}$ for some $0 \leqslant k \leqslant N$. The result now immediately follows by observing that $\overline{\delta_k}(\mathrm{Id}) = k \ \forall \ 0 \leqslant k \leqslant N$.    
\end{proof}

\bcrlre \label{CW}
Let $V$ be any irreducible $\mathfrak{gl}_{N}$-module. Then:
\begin{enumerate}
\item The rank-one operators of $\mathfrak{gl}_{N}$ are given by $\mathcal{R}_1(\mathfrak{gl}_{N}) := \{0 \neq ru^T \ | \ u, r \in \mathbb{C}^{N} \}$.
\item $\mathfrak{gl}_{N}$ is spanned by $\mathcal{R}_1(\mathfrak{gl}_{N})$.
\item $x^2-\mathrm{tr}(x)x  \in \mathrm{Ann}_{\mathrm{U}(\mathfrak{gl}_{N})}V$ for all $x \in \mathcal{R}_1(\mathfrak{gl}_{N})$ if and only if $V$ is the trivial or a fundamental representation of $\mathfrak{gl}_{N}$.
\end{enumerate}
\ecrlre

\bthm \label{TW}
For any irreducible $\mathfrak{gl}_N$-module $V, \ F^{\beta}_{\mathcal{W}}(V)$ is irreducible over $\mathcal{W}_N$ if $V \ncong V(\delta_k,k)$ for any $0 \leqslant k \leqslant N$. 
\ethm

\begin{proof}
 Let $V \ncong V(\delta_k,k)$ for any $0 \leqslant k \leqslant N$ and $W$ be a non-zero $\mathcal{W}_N$-submodule of $F^{\beta}_{\mathcal{W}}(V)$. As in Theorem \ref{Criterion}, put $\widetilde{W} = \bigcap_{m \in \mathbb{Z}^N} W_m$ and apply Lemma \ref{LW} to obtain $r \in \mathbb{Z}^N,\ u \in \mathbb{C}^N$ such that $(ru^T)^2v \neq (u|r)(ru^T)v$ for every non-zero $v \in V$. Fix $0 \neq w \in W_m$ for some $m \in \mathbb{Z}^N$ and consider
\begin{align*}
D(u,-r)D(u,r+s)(w \otimes t^m) = [(u|r)(ru^T) -(ru^T)^2]w \otimes t^{m+s}+ \ldots \ldots \in W_{m+s} \otimes \mathbb{C}t^{m+s}.
\end{align*} 
Replacing $r$ by $pr$ and $u$ by $qu$ for $p,q \in \mathbb{Z}$ in the above equation and comparing the coefficients of $p^2q^2$, we can then apply Lemma \ref{L2} to infer that $[(ru^T)^2 -(u|r)(ru^T)]w \in W_{m+s} \ \forall \ s \in \mathbb{Z}^N$. This shows that $\widetilde{W} \neq (0)$. Again for any $m, s \in \mathbb{Z}^N, \ x \in \mathbb{C}^N$ and $v \in \widetilde{W}$, we have
\begin{align*}
D(x,m)(v \otimes t^{s-m}) = (x|s-m+\beta)v \otimes t^s + (mx^T)v \otimes t^s,
\end{align*}
which gives $(mx^T)v \in \widetilde{W}$. Applying Lemma \ref{LW}, we then get $\widetilde{W}= V$ and so we are done.	
\end{proof}

\blmma \label{LS}
Let $V$ be an irreducible $\mathfrak{sl}_N$-module (not necessarily a weight module).  Then:  
\begin{enumerate}
	\item $\mathfrak{sl}_N =  \mathrm{span}\{ru^T \ | \ r \in \mathbb{Z}^N, \ u \in \mathbb{C}^N \ \text{with} \ (u|r)=0 \} $. 
	\item $\mathcal{J}_{\mathcal{S}}(V) := \{v \in V \ | \ (ru^{T})^2v = 0 \ \forall \ u \in \mathbb{C}^N, \ r \in \mathbb{Z}^N \ with \ (u|r) = 0 \}$ is an $\mathfrak{sl}_N$-module.
	\item $\mathcal{J}_{\mathcal{S}}(V) = V$ if and only if $V \cong V(\delta_k)$ for some $0 \leqslant k \leqslant N-1$, where $\delta_0 = 0$.  
\end{enumerate}
\elmma

\begin{proof}
	(1) follows from Lemma \ref{LW}, whereas (2) can be deduced using Lemma \ref{L2}, by providing essentially an analogous argument as presented in Lemma \ref{L3}. \\
	(3) If $V \cong V(\delta_k)$ for some $0 \leqslant k \leqslant N-1$, then a direct checking readily yields $\mathcal{J}_{\mathcal{S}}(V) = V$. \\
	Conversely, assume that $\mathcal{J}_{\mathcal{S}}(V) = V$. Now putting $r=e_i$ and $u=e_j$ for any $1 \leqslant i \neq j \leqslant N$, we get $(E_{i,j})^2v = 0 \ \forall \ v \in V$. This shows that $V$ is finite-dimensional by Lemma \ref{L1}. Let $V \cong V(\lambda)$ for some $\lambda \in P_{\mathfrak{sl}_N}^{+} \setminus \{\delta_k \}_{k=0}^{N-1}$. Consequently, there exist $a_1, \ldots, a_{N-1} \in \mathbb{Z}_{\geqslant 0}$ such that $\lambda = \sum_{k=1}^{N-1} a_k \delta_k$ with $a_i + a_j \geqslant 2$ for some $1 \leqslant i \leqslant j < N$. Note that $\delta_i(E_{i,i} - E_{j+1,j+1}) = 1 = \delta_j(E_{i,i} - E_{j+1,j+1})$, which implies that $\lambda(E_{i,i} - E_{j+1,j+1}) \geqslant 2$. Now pick a non-zero $v_{\lambda} \in V(\lambda)_{\lambda}$ and consider the $S_{i,j+1}$-module $\text{span} \{(E_{j+1,i})^rv_{\lambda} \ | \ r \in \mathbb{Z}_{\geqslant 0}  \}$, where $S_{i,j+1} := \text{span} \{E_{i,j+1}, E_{j+1,i}, E_{i,i} - E_{j+1,j+1} \}$. Then by $\mathfrak{sl}_2$-theory, it immediately follows that $(E_{j+1,i})^2v_{\lambda} \neq 0$ and thus $v_{\lambda} \notin \mathcal{J}_{\mathcal{S}}(V(\lambda))$. 
\end{proof}

\bcrlre \label{CS}
Let $V$ be any irreducible $\mathfrak{sl}_{N}$-module. Then:
\begin{enumerate}
\item The rank-one operators of $\mathfrak{sl}_{N}$ are given by $\mathcal{R}_1(\mathfrak{sl}_{N}) := \{0 \neq ru^T \ | \ u, r \in \mathbb{C}^{N}, (u|r)=0 \}$.
\item $\mathfrak{sl}_{N}$ is spanned by $\mathcal{R}_1(\mathfrak{sl}_{N})$.
\item $x^2  \in \mathrm{Ann}_{\mathrm{U}(\mathfrak{sl}_{N})}V$ for all $x \in \mathcal{R}_1(\mathfrak{sl}_{N})$ if and only if $V$ is the trivial or a fundamental representation of $\mathfrak{sl}_{N}$.
\end{enumerate}
\ecrlre

\bthm \label{TS}
For any irreducible $\mathfrak{sl}_N$-module $V, \ F^{\beta}_{\mathcal{S}}(V)$ is irreducible over $\mathcal{S}_N$ if $V \ncong V(\delta_k)$ for any $0 \leqslant k \leqslant N-1$.
\ethm

\begin{proof}
Let $V \ncong V(\delta_k)$ for any $0 \leqslant k < N$ and $W$ be a non-zero $\mathcal{S}_N$-submodule of $F^{\beta}_{\mathcal{S}}(V)$. As in Theorem \ref{Criterion}, put $\widetilde{W} = \bigcap_{m \in \mathbb{Z}^N} W_m$ and then invoke Lemma \ref{LS} to obtain $r \in \mathbb{Z}^N,\ u \in \mathbb{C}^N$ with $(u|r)=0$ such that $(ru^T)^2v \neq  0$ for every non-zero $v \in V$. Fix $0 \neq w \in W_m$ for some $m \in \mathbb{Z}^N$. Consequently, proceeding exactly as in Theorem \ref{TW} yields $(ru^T)^2w \in W_{m+s} \ \forall \ s \in \mathbb{Z}^N$ with $(u|s)=0$. Let us now pick any $s \in \mathbb{Z}^N$ such that $(u|s) \neq 0$. Setting $u_1 = u$, extend $u_1$ to a basis $\{u_i \}_{i=1}^{N}$ of $\mathbb{C}^N$. Then $(u_i|r) \neq 0$ for some $2 \leqslant i \leqslant N$. So without loss of generality, we can take $(u_2|r) \neq 0$ and define $u_s^{\prime} := (u_1|s)u_2 - (u_2|r+s)u_1 \neq 0$. Note that $(u_s^{\prime}|r+s) =0$ and consider
\begin{align*}
D(u,-r)D(u_s^{\prime},r+s)(w \otimes t^m) = -(u_2|r)[(u|r)(ru^T) + (ru^T)^2](w \otimes t^{m+s}) + \ldots \in W_{m+s} \otimes \mathbb{C}t^{m+s}.
\end{align*} 
Now replacing $r$ by $pr$ and $u$ by $qu$ for $p,q \in \mathbb{Z}$ in the above equation and comparing the coefficients of $p^3q^2$, we can  apply Lemma \ref{L2} to obtain $(ru^T)^2w \in W_{m+s} \ \forall \ s \in \mathbb{Z}^N$ satisfying $(u|s) \neq 0$. This implies that $\widetilde{W} \neq (0)$ and so by providing the same argument as in Theorem \ref{TW}, we get $(mx^T)v \in \widetilde{W} \ \forall \ m \in \mathbb{Z}^N, \ x \in \mathbb{C}^N$ satisfying $(x|m)=0$. Applying Lemma \ref{LS} finally gives $\widetilde{W}= V$ and thus we obtain the desired result.	
\end{proof}

\subsection{Connection to minimal nilpotent orbits of semi-simple Lie algebras.} 


Let $\mathfrak{g}$ be a semi-simple Lie algebra and $G$ be the connected algebraic (adjoint) group with Lie algebra $\mathfrak{g}$. 
For $x \in \mathfrak{g}$, the orbit of $x$ (under the \textit{adjoint representation} $\mathrm{Ad}$ of $G$) is defined as
\begin{align*}
\mathcal{O}^{\mathfrak{g}}_x := \{\mathrm{Ad}_{\mathfrak{g}}(x) \ | \ g \in G \} = \{gxg^{-1} \ | \ g \in G \}.
\end{align*}
Then it is well-known that the \textit{minimal nilpotent orbit} of $\mathfrak{g}$, which we shall denote by $\mathcal{O}_{min}^{\mathfrak{g}}$, is given by the $G$-orbit of $e_\theta$, where $e_{\theta}$ is the highest root vector of $\mathfrak{g}$, i.e. $\mathcal{O}_{min}^{\mathfrak{g}} = G.e_{\theta}$.

\smallskip

The adjoint group of $\mathfrak{sl}_N$ is $\mathbb{PSL}_N =  \mathbb{SL}_N/ \mathrm{Z}(\mathbb{SL}_N)$, where $\mathbb{SL}_N$ is the special linear group and $\mathrm{Z}(\mathbb{SL}_N)$ denotes the center of $\mathbb{SL}_N$. Now 
the \textit{minimal nilpotent orbit} of $\mathfrak{sl}_N$ is given by the $\mathbb{SL}_N$-orbit of $E_{1,N}$, from which it can be easily deduced that $\mathcal{O}_{min}^{\mathfrak{sl}_N} = \mathcal{R}_1(\mathfrak{sl}_N)$, which spans $\mathfrak{sl}_N$.

Similarly, the adjoint group of $\mathfrak{sp}_{2n}$ is $\mathbb{PS}\mathrm{p}_{2n} = \mathbb{S}\mathrm{p}_{2n}/ \mathrm{Z}(\mathbb{S}\mathrm{p}_{2n})$, where $\mathbb{S}\mathrm{p}_{2n}$ denotes the symplectic group and $Z(\mathbb{S}\mathrm{p}_{2n})$ is the center of $\mathbb{S}\mathrm{p}_{2n}$. The \textit{minimal nilpotent orbit} of $\mathfrak{sp}_{2n}$ coincides with $\{X.E_{1,n+1} \ | \ X \in \mathbb{S}\mathrm{p}_{2n} \}$, which in turn yields that $\mathcal{O}_{min}^{\mathfrak{sp}_{2n}} =\mathcal{R}_1(\mathfrak{sp}_{2n})$, which spans $\mathfrak{sp}_{2n}$.

As a consequence of Lemma \ref{LS}, Lemma \ref{L3} and the above discussion, we obtain the following result, which may also be of independent interest.

\bcrlre \label{Nilpotent}
Let $\mathfrak{g}_{s} = \mathfrak{sl}_N$ or $\mathfrak{sp}_{2n}$ and V be any irreducible module over $\mathfrak{g}_{s}$ (need not be a weight module). Then $x^2 \in \mathrm{Ann}_{\mathrm{U}(\mathfrak{g}_{s})}V \ \forall \ x \in \mathcal{O}_{min}^{\mathfrak{g}_{s}}$ if and only if V is the trivial representation or a fundamental representation of $\mathfrak{g}_{s}$. 
\ecrlre

\brmk \label{Orbit}
\
\begin{enumerate}
\item As $\mathfrak{gl}_N$ is not semi-simple and contains matrices with \textit{non-zero} trace, we obtain an orbit decomposition, instead of a single orbit. More precisely, the general linear group $\mathbb{GL}_N$ acts on $\mathfrak{gl}_N$ under the adjoint action and $g(ru^T)g^{-1} = (gr)\big((g^{-1})^Tu\big)^T \ \forall \ g \in \mathbb{GL}_N, \ u,r \in \mathbb{C}^N$, which shows that $\mathbb{GL}_N$ acts on $\mathcal{R}_1(\mathfrak{gl}_N) := \{0 \neq ru^T \ | \ u, r \in \mathbb{C}^N \}$. Also, we clearly have $(gr|(g^{-1})^Tu) = \big((g^{-1})^Tu)\big)^Tgr = u^Tr = (u|r)$, from which it follows that
 \begin{align*}
\mathcal{R}_1(\mathfrak{gl}_N) = \bigsqcup_{\gamma \in \mathbb{C}} \mathrm{Orb}(\gamma), \ \text{where} \ \mathrm{Orb}(\gamma) = \{0 \neq ru^T \ | \ u, r \in \mathbb{C}^N, (u|r) = \gamma \},
\end{align*}
which is the fiber decomposition of the regular function $\mathrm{tr} : \mathcal{R}_1(\mathfrak{gl}_N) \longrightarrow \mathbb{C}$ over its values. 
\item Corollary \ref{Nilpotent} does not hold for $\mathfrak{so}_N$, i.e. for Lie algebras of type $B_n$ ($n \geqslant 3$) and $D_n$ ($n \geqslant 4$). Specifically, the second fundamental representation $V(\omega_2)$ of $\mathfrak{so}_N$ \textit{does not} satisfy the condition ``$x^2 \in \mathrm{Ann}_{\mathrm{U}(\mathfrak{so}_N)}V(\omega_2)$ for all $x \in \mathcal{O}_{min}^{\mathfrak{so}_N}$" . This is essentially because for types $B$ and $D$, the minimal nilpotent orbit consists of \textit{rank-two} skew-symmetric matrices rather than rank-one matrices. This highlights a special feature of types $A$ and $C$; the equality $\mathcal{O}_{min}^{\mathfrak{g}_{s}} = \mathcal{R}_1(\mathfrak{g}_{s})$ gives the rank-one quadratic condition its full rigidity in those cases.
\end{enumerate}
\ermk

\brmk
The irreducibility for the Shen--Larsson modules over both $\mathcal{W}_N$ and $\mathcal{S}_N$ has been studied extensively over the last three decades (see \cite{R1,GZ,LZ,LGW,T1} and the references therein).  In Theorem \ref{TW} and Theorem \ref{TS}, we recover these results uniformly by adapting the quadratic criterion approach of rank-one operators that we have also used for $\mathcal{H}_N$, which makes our proofs much more compact and streamlined than the ones already existing in the literature.
\ermk

\section{Classification of submodules of exceptional Shen--Larsson modules over $\mathcal{W}_N$}\label{S9}
In this section, we classify all possible $\mathcal{W}_N$-submodules of the exceptional Shen--Larsson modules $F^{\beta}(\Lambda^p\mathbb{C}^N), \ 0 \leqslant p \leqslant N$, by means of similar techniques that we have developed for $\mathcal{H}_N$, which subsequently reveals that all these modules are, in fact indecomposable. To this end, we shall work with the following orthogonal basis of $\mathbb{C}^N$ for any $k \in \mathbb{Z}^N$ satisfying $k+\beta \neq 0$, which will be extensively utilized throughout this section.

Extend $k+\beta$ to an orthogonal basis of $\mathbb{C}^N$, say $\{k+\beta, v_1^{(k)}, \ldots, v_{N-1}^{(k)} \}$ and denote the subspace of $\mathbb{C}^N$ consisting of vectors orthogonal to $k+\beta$ by $\{k+\beta \}^{\perp}$, i.e. $\{k+\beta \}^{\perp} = \oplus_{i=1}^{N-1}\mathbb{C}v_i^{(k)}$. 

\subsection{Rank-reducing operators of $\mathfrak{gl}_N$} \label{W-invariant}
For each $k \in \mathbb{Z}^N$ and $u, v, \beta \in \mathbb{C}^N$, put 
\begin{align*}
\mathcal{T}_{k,u,v}^{\mathcal{W},\beta} = (k+\beta|u)v - (k+\beta|v)u. 
\end{align*}
For any $r, s \in \mathbb{Z}^N$, let us now define the corresponding \textit{rank-reducing operators} of $\mathfrak{gl}_{N}$ by setting
\begin{align*}
\Omega_{r,s,u,v}^{\mathcal{W},k,\beta} = \mathcal{T}_{k,r,s}^{\mathcal{W},\beta}\big({\mathcal{T}_{k,u,v}^{\mathcal{W},\beta}}\big)^T. 
\end{align*}

\bppsn \label{Invariances}
Let $M = \oplus_{k \in \mathbb{Z}^N}(M_k \otimes \mathbb{C}t^k)$ be a non-trivial $\mathcal{W}_N$-submodule of $F^{\beta}_{\mathcal{W}}(V)$, where $V$ is a $\mathfrak{gl}_N$-module. Then:
\begin{enumerate}
\item $M_k$ is invariant under $\Omega_{r,s,u,v}^{\mathcal{W},k,\beta}$ for each $k,r,s \in \mathbb{Z}^N$ and $u, v, \beta \in \mathbb{C}^N$. 
\item If $(k+\beta)_i \neq 0$ for some $1 \leqslant i \leqslant N$, then $\{\mathcal{T}_{k,r,s}^{\mathcal{W},\beta} \ | \ r, s \in \mathbb{Z}^N \}$ spans $\{k+\beta \}^{\perp}$. Moreover, $\{\mathcal{T}_{k,e_i,e_j}^{\mathcal{W},\beta} \ | \ 1 \leqslant j \neq i \leqslant N \}$ is a basis of $\{k+\beta \}^{\perp}$. 
\end{enumerate}
\eppsn

\begin{proof}
(1) Choose $w \in \mathbb{C}^N$ such that $(r|w) \neq 0$. Now for any $v^{\prime} \in M_k$, consider
\begin{align*}
D(v,-(r+s))D(w,s)D(u,r).(v^{\prime} \otimes t^k) = - (r|w)[(k+\beta|u)rv^T - (k+\beta|v)ru^T].(v^{\prime} \otimes t^k) + \ldots
\end{align*}
Replacing $r$ by $pr, s$ by $qs$ and $u$ by $p^{\prime}u, v$ by $q^{\prime}v$ for $p, p^{\prime}, q, q^{\prime} \in \mathbb{Z}$ in the above equation and comparing the coeffecients of $p^2p^{\prime}q^{\prime}$, we can then apply Lemma \ref{L2} to get
\begin{align*}
[(k+\beta|u)rv^T - (k+\beta|v)ru^T]v^{\prime} \in W_k \ \forall \ r \in \mathbb{Z}^N \ \text{and} \ u, v \in \mathbb{C}^N. \\
\implies (k+\beta|r)[(k+\beta|u)sv^T - (k+\beta|v)su^T].v^{\prime} - (k+\beta|s)[(k+\beta|u)rv^T - (k+\beta|v)ru^T ].v^{\prime} \in W_k,
\end{align*}
which shows that $\Omega_{r,s,u,v}^{\mathcal{W},k,\beta}.v^{\prime} \in M_k$ and therefore (1) is proved. \\
(2) follows as $\{\mathcal{T}_{k,e_i,e_j}^{\mathcal{W},\beta} \ | \ 1 \leqslant j \neq i \leqslant N \}$ is a linearly independent subset of $\{k+\beta \}^{\perp}$, since $\{k+\beta \}^{\perp}$ is $(N-1)$-dimensional.
\end{proof}

\brmk \label{Rank-reducing}
We can construct similar \textit{rank-reducing operators} of $\mathfrak{sl}_N$, which enables us to show that each $\mathbb{Z}^N$-graded component of any submodule of a Shen--Larrson module $F_{\mathcal{S}}^{\beta}(V)$ over $\mathcal{S}_N$ is an $\mathfrak{sl}_{N-1}$-module. More precisely, for $u,v, \beta \in \mathbb{C}^N$ and $k,r,s \in \mathbb{Z}^N$ with $(u|r)=0=(v|s)$, define: 
\begin{align*}
\Omega_{r,s,u,v}^{\mathcal{S},k,\beta} := \big(\mathcal{P}_{\{k+\beta \}^{\perp}}[(v|k+\beta)(u|s)r + (u|k+\beta)(v|r)s] \big)\big((v|k+\beta)u-(u|k+\beta)v\big)^T,
\end{align*}
which spans $\mathfrak{sl}_{N-1}$ as a Lie algebra. Here $\mathcal{P}_{\{k+\beta \}^{\perp}}$ denotes the orthogonal projection onto $\{k+\beta \}^{\perp}$.
\ermk

\subsection{Indecomposable $\mathcal{W}_N$-modules} For any $k \in \mathbb{Z}^N$ satisfying $k+\beta \neq 0$, the orthogonal basis $\{k+\beta, v_1^{(k)}, v_2^{(k)}, \ldots, v_{N-1}^{(k)} \}$ of $\mathbb{C}^N$ again contains an orthogonal basis $\{v_i^{(k)} \}_{i=1}^{N-1}$ of $\mathbb{C}^{N-1}$. This gives rise to the Lie algebra $\mathrm{span}\{v_i^{(k)}{v_j^{(k)}}^T \ | \ 1 \leqslant i,j < N \} \cong \mathfrak{gl}_{N-1}$, which also consists of the subalgebra $\mathrm{span}\{v_i^{(k)}{v_j^{(k)}}^T \ | \ 1 \leqslant i \neq j < N \} \cong \mathfrak{sl}_{N-1}$. We shall denote these Lie algebras by $\mathfrak{gl}_{N-1}^{(k)}$ and $\mathfrak{sl}_{N-1}^{(k)}$ respectively. The corresponding Cartan subalgebras can be taken to be $\mathfrak{h}_{\mathfrak{gl}_{N-1}}^{(k)} := \mathrm{span}\{v_i^{(k)}{v_i^{(k)}}^T \ | \ 1 \leqslant i < N \}$ and $\mathfrak{h}_{\mathfrak{sl}_{N-1}}^{(k)} := \mathrm{span}\{v_i^{(k)}{v_i^{(k)}}^T - v_j^{(k)}{v_j^{(k)}}^T \ | \ 1 \leqslant i \neq j < N \}$ respectively. It is also clear that $\mathfrak{gl}_{N-1}^{(k)}$ acts on $\mathbb{C}^{N-1}$, where $\mathbb{C}^{N-1} \cong \bigoplus_{i=1}^{N-1}\mathbb{C}v_i^{(k)}$. Set $\overline{H_i^{(k)}} = v_i^{(k)}{v_i^{(k)}}^T$ and $\overline{\epsilon_i^{(k)}} \in (\mathfrak{h}_{\mathfrak{gl}_{N-1}}^{(k)})^*$ such that $\overline{\epsilon_i^{(k)}}\big(\overline{H_j^{(k)}}\big) = \delta_{j,i} \ \forall \ 1 \leqslant i < N$. Putting $\widehat{\epsilon_i^{(k)}} := \overline{\epsilon_i^{(k)}}|_{\mathfrak{h}_{\mathfrak{sl}_{N-1}}^{(k)}}$, this induces the weight space decomposition of $\mathbb{C}^{N-1}$ with respect to $\mathfrak{h}_{\mathfrak{sl}_{N-1}}^{(k)}$, given by  
$\mathbb{C}^{N-1} = \bigoplus_{i=1}^{N-1} (\mathbb{C}^{N-1})_{\widehat{\epsilon_i^{(k)}}}$, where $(\mathbb{C}^{N-1})_{\widehat{\epsilon_i^{(k)}}} = \mathbb{C}v_{i}^{(k)} \ \forall \ 1 \leqslant i < N$. Furthermore, observe that $\mathfrak{gl}_{N-1}^{(k)}$ acts trivially on $k+\beta$.

\blmma \label{Restrict}
Let $M = \oplus_{k \in \mathbb{Z}^N}(M_k \otimes \mathbb{C}t^k)$ be a non-trivial $\mathcal{W}_N$-submodule of $F^{\beta}_{\mathcal{W}}(V)$, where $V$ is a $\mathfrak{gl}_N$-module. Then $\mathfrak{gl}_{N-1}^{(k)}$ acts on $M_k$ for $k \in \mathbb{Z}^N$ satisfying $k+\beta \neq 0$, where $\mathfrak{gl}_{N-1}^{(k)} \cong \mathfrak{gl}(\{k+\beta \}^{\perp})$.
\elmma

\begin{proof}
The lemma follows by just applying Proposition \ref{Invariances}, since $v_i^{(k)} \in \{k+\beta \}^{\perp} \ \forall \ 1 \leqslant i < N$.
\end{proof}

\bthm\cite{ER,R1,GZ} \label{Recall}
\begin{enumerate}
\item If $p \in \{0,N\}$ and $\beta \notin \mathbb{Z}^N$, then $F_{\mathcal{W}}\big(V(\overline{\delta_p})\big)$ is irreducible.
\item If $p \in \{0,N\}$ and $\beta \in \mathbb{Z}^N$, then a Loewy series of $F_{\mathcal{W}}\big(V(\overline{\delta_p})\big)$ is given by 
\begin{align*}
(0) \subseteq V(\overline{\delta_p}) \otimes \mathbb{C}t^{-\beta} \subseteq V(\overline{\delta_p}) \otimes A_N, \\ \text{where the Loewy layers of} \ V(\overline{\delta_p}) \otimes A_N \ \text{are} \ W^{\beta}(\Lambda^1\mathbb{C}^N) \ \text{and} \ V(\overline{\delta_p}) \otimes \mathbb{C}t^{-\beta}.
\end{align*}

\item If $1 \leqslant p < N$ and $\beta \notin \mathbb{Z}^N$, then a Loewy series of $F_{\mathcal{W}}\big(V(\overline{\delta_0})\big)$ is given by
\begin{align*}
(0) \subseteq W^{\beta}(\Lambda^p\mathbb{C}^N) \subseteq F_{\mathcal{W}}\big(V(\overline{\delta_p})\big), \\ \text{where the Loewy layers of} \ F_{\mathcal{W}}\big(V(\overline{\delta_p})\big) \ \text{are} \ W^{\beta}(\Lambda^p\mathbb{C}^N) \ \text{and} \ W^{\beta}(\Lambda^{p+1}\mathbb{C}^N).
\end{align*}
\item If $1 \leqslant p < N$ and $\beta \in \mathbb{Z}^N$, then a Loewy series of $F_{\mathcal{W}}\big(V(\overline{\delta_p})\big)$ is given by 
\begin{align*}
(0) \subseteq \widehat{W}^{\beta}(\Lambda^p\mathbb{C}^N) \subseteq F_{\mathcal{W}}\big(V(\overline{\delta_p})\big), \\ \text{where the Loewy layers of} \ F_{\mathcal{W}}\big(V(\overline{\delta_p})\big) \ \text{are} \ W^{\beta}(\Lambda^p\mathbb{C}^N) \ \text{and} \ W^{\beta}(\Lambda^{p+1}\mathbb{C}^N), \\
\text{with} \ \widehat{W}^{\beta}(\Lambda^p\mathbb{C}^N)/W^{\beta}(\Lambda^p\mathbb{C}^N) \cong \mathbb{C}^{N \choose p}.
\end{align*} 
\item Let $M = \oplus_{k \in \mathbb{Z}^N}M_k \otimes \mathbb{C}t^k$ be a $\mathcal{W}_N$-submodule of $F^{\beta}_{\mathcal{W}}\big(V(\overline{\delta_p})\big)$, where $1 \leqslant p < N$. If there exists $k \in \mathbb{Z}^N$ with $k+\beta \neq 0$ such that $M_k \cap \mathrm{span}\{v_{i_1}^{(k)} \wedge \ldots \wedge v_{i_{p}}^{(k)} \ | \ 1 \leqslant i_j < N \} \neq (0)$, then $M = F^{\beta}_{\mathcal{W}}\big(V(\overline{\delta_p})\big)$.  
\end{enumerate}
\ethm

\blmma \label{Unique}
\
\begin{enumerate}
\item If $p \in \{0,N\}$ and $\beta \in \mathbb{Z}^N$, then $V(\overline{\delta_p}) \otimes \mathbb{C}t^{-\beta}$ is the unique submodule of $F^{\beta}_{\mathcal{W}}\big(V(\overline{\delta_p})\big)$. 
\item If $1 \leqslant p < N$, then $W^{\beta}(\Lambda^p\mathbb{C}^N)$ is the unique irreducible submodule of $F^{\beta}_{\mathcal{W}}\big(V(\overline{\delta_p})\big)$.
\end{enumerate}
\elmma

\begin{proof}
(1) Let $M = \oplus_{k \in \mathbb{Z}^N}M_k \otimes \mathbb{C}t^k$ be a non-zero proper $\mathcal{W}_N$-submodule of $F^{\beta}_{\mathcal{W}}\big(V(\overline{\delta_p})\big)$. In view of Theorem \ref{Recall}, we must have $M \subseteq \oplus_{k \in \mathbb{Z}^N, k \neq -\beta}V(\overline{\delta_p}) \otimes \mathbb{C}t^k \cong \oplus_{k \in \mathbb{Z}^N, k \neq -\beta} \mathbb{C}t^k$. So there exists $k \in \mathbb{Z}^N$ such that $k+\beta \neq 0$ and $M_k \neq (0)$. But this gives 
\begin{align*}
D(k+\beta,-k-\beta).t^k = -(k+\beta|k+\beta)t^{-\beta} \neq 0,
\end{align*} 
which is a contradiction and hence the assertion is proved. \\
(2) Let $W = \oplus_{k \in \mathbb{Z}^N}W_k \otimes \mathbb{C}t^k$ be an irreducible $\mathcal{W}_N$-submodule of $F^{\beta}_{\mathcal{W}}\big(V(\overline{\delta_p})\big)$ and $k \in \mathbb{Z}^N$ satisfying $k+\beta \neq 0$ and $W_k \neq (0)$. By Theorem \ref{Recall}, there exist $a_{i_1, \ldots,i_{p-1}} \in \mathbb{C}$ and $b_{j_1, \ldots,j_{p}} \in \mathbb{C}$ such that  
\begin{align*}
v = \sum_{(i_1, \ldots, i_{p-1})}a_{i_1, \ldots,i_{p-1}}(k+\beta) \wedge v_{i_1}^{(k)} \wedge \ldots \wedge v_{i_{p-1}}^{(k)} + \sum_{(j_1, \ldots, j_p)}b_{j_1, \ldots,j_{p}}v_{j_1}^{(k)} \wedge \ldots \wedge v_{j_{p}}^{(k)} \in W_k, 
\end{align*}
with $\sum_{(i_1, \ldots, i_{p-1})}a_{i_1, \ldots,i_{p-1}}(k+\beta) \wedge v_{i_1}^{(k)} \wedge \ldots \wedge v_{i_{p-1}}^{(k)} \neq 0$, where $1 \leqslant i_r, j_s < N$. By Lemma \ref{Restrict} and the discussion preceding this lemma, it now follows that the $\mathfrak{h}_{\mathfrak{sl}_{N-1}}^{(k)}$-weights of $(k+\beta) \wedge v_{i_1}^{(k)} \wedge \ldots \wedge v_{i_{p-1}}^{(k)}$ and $v_{j_1}^{(k)} \wedge \ldots \wedge v_{j_{p}}^{(k)}$ are different. This implies that $\sum_{(i_1, \ldots, i_{p-1})}a_{i_1, \ldots,i_{p-1}}(k+\beta) \wedge v_{i_1}^{(k)} \wedge \ldots \wedge v_{i_{p-1}}^{(k)} \in W_k$, which thereby proves the lemma, due to the irreducibility of $W^{\beta}(\Lambda^p\mathbb{C}^N)$ over $\mathcal{W}_N$.    
\end{proof}

\brmk \label{Indecomposable functor}
Proceeding analogously as in Theorem \ref{Recall} and Lemma \ref{Unique}, it can be shown that $F^{\beta}_{\mathcal{W}}\big(V(\overline{\delta_p})\big)$ is \textit{indecomposable} for all $1 \leqslant p < N$.
\ermk

\bthm\label{Classify}
\
Let $M$ be a non-zero proper $\mathcal{W}_N$-submodule of $F_{\mathcal{W}}^{\beta}\big(V(\overline{\delta_p})\big)$, where $0 \leqslant p \leqslant N$.
\begin{enumerate}
\item If $p \in \{0,N \}$ and $\beta \in \mathbb{Z}^N$, then $M$ is isomorphic to $V(\overline{\delta_p}) \otimes \mathbb{C}t^{-\beta}$. 
\item If $1 \leqslant p < N$ and $\beta \notin \mathbb{Z}^N$, then $M$ is isomorphic to $W^{\beta}(\Lambda^p\mathbb{C}^N)$.
\item If $1 \leqslant p < N$ and $\beta \in \mathbb{Z}^N$, then $M$ is isomorphic to either $W^{\beta}(\Lambda^p\mathbb{C}^N) \oplus W \otimes \mathbb{C}t^{-\beta}$, where $W$ is an arbitrary subspace (including the zero subspace) of $V(\overline{\delta_p})$.   
\end{enumerate}
\ethm

\begin{proof}
The proof follows from Remark \ref{Indecomposable functor}, as the Loewy series of $F^{\beta}_{\mathcal{W}}\big(V(\overline{\delta_p})\big)$ has length $2$, by Theorem \ref{Recall}. 
\end{proof}

\bcrlre \label{CW}
\
\begin{enumerate}
\item For $0 \leqslant p \leqslant N, \ F^{\beta}_{\mathcal{W}}\big(V(\overline{\delta_p})\big)$ is \textit{indecomposable} and has Loewy Length $2$ for all $\beta \in \mathbb{C}^N$.
\item If $p \in \{0,N \}$ and $\beta \in \mathbb{Z}^N$, then $F^{\beta}_{\mathcal{W}}\big(V(\overline{\delta_p})\big)$ is \textit{uniserial}.
\item For $1 \leqslant p < N, \ F^{\beta}_{\mathcal{W}}\big(V(\overline{\delta_p})\big)$ is uniserial if $\beta \notin \mathbb{Z}^N$.
\end{enumerate}
\ecrlre

\brmk
\
\begin{enumerate}
\item It is clear from Theorem \ref{Recall} and Theorem \ref{Classify} that $F_{\mathcal{W}}^{\beta}\big(V(\overline{\delta_p})\big)$ admits a \textit{unique} irreducible quotient for each $0 \leqslant p \leqslant N$.   
\item Recall that $\mathcal{H}_N = \text{span}\{D(Jr,r), d_i \ | \  r \in \mathbb{Z}^N, \ 1 \leqslant i \leqslant N  \}$, where $J$ is the skew-symmetric matrix mentioned in Remark \ref{Simple}. Now if we replace $J$ by any non-zero non-degenerate skew-symmetric matrix $P$, then we can form the corresponding derivation algebra 
$\mathcal{H}_N(P) := \text{span}\{D(Pr,r), d_i \ | \ r \in \mathbb{Z}^N, \ 1 \leqslant i \leqslant N \}$ and all our results again hold good for $\mathcal{H}_N(P)$. This derivation algebra appears in the construction of \textit{skew-symmetric extended affine Lie algebras}, which were recently introduced in \cite{CR}. It is also known that given any skew-symmetric matrix $P$ of order $N$, there exists $A \in \mathbb{GL}_{N}$ such that $J = APA^T$, but $\mathcal{H}_N$ and $\mathcal{H}_N(P)$ are not necessarily isomorphic. However, if all the entries of $A$ are known to be integers, then we have $\mathcal{H}_N \cong \mathcal{H}_N(P)$ (see \cite[Lemma 3.1]{CR} and \cite[Proposition 3.2]{CR}).
\end{enumerate}
\ermk

\brmk \label{RS}
As in the case of $\mathcal{W}$ and $\mathcal{H}$ types, the submodules of the exceptional Shen--Larsson modules over $\mathcal{S}_N$ can be classified using the rank-reducing operators $\Omega_{r,s,u,v}^{\mathcal{S},k,\beta}$ in Remark \ref{Rank-reducing} (this classification was also obtained in \cite{LGW} using different methods), from which it follows that:
\begin{enumerate}
\item For $1 \leqslant p < N$ and $\beta \in \mathbb{C}^N, \  F_{\mathcal{S}}^{\beta}\big(V(\delta_p)\big)$ is \textit{indecomposable} admitting a \textit{unique} irreducible submodule and a \textit{unique} irreducible quotient.
\item For $1 \leqslant p < N$ and $\beta \in \mathbb{C}^N, \ F_{\mathcal{S}}^{\beta}\big(V(\delta_p)\big)$ has \textit{Loewy length} $2$.
\item For $1 \leqslant p < N, \ F_{\mathcal{S}}^{\beta}\big(V(\delta_p)\big)$ is \textit{uniserial} if $\beta \notin \mathbb{Z}^N$.
\end{enumerate}
\ermk
 
\noindent \textbf{Acknowledgements.} Most of this work was done when the second author was a post-doc at Indian Institute of Science (IISc Bangalore), where he was supported by a post-doctoral fellowship from the National Board for Higher Mathematics (Ref. No. 0204/9/2024/R\&D-II/2965).

\end{document}